\documentclass[12pt]{amsart}
\usepackage{amsfonts,amsthm,amssymb}
\allowdisplaybreaks
\usepackage[mathscr]{eucal}
\usepackage{stmaryrd}
\usepackage{graphicx}
\usepackage{times}

\usepackage[curve,matrix,arrow,frame]{xy}

\usepackage[bookmarks=false, hyperfigures=false, a4paper]{hyperref}

\advance\textheight\topskip
\newcommand{\AT}{e^{-i\pi \frac{c}{12}}}

\newcommand{\Aconstplain}{\sqrt{\ffrac{p_+ p_-}{2}}}
\newcommand{\Aoverviiipp}{\ffrac{1}{\sqrt{\mathstrut2 p_+ p_-}}}
\newcommand{\iiA}{\sqrt{\mathstrut2 p_+ p_-\!}}
\newcommand{\iippA}{\sqrt{2}\,(p_+ p_-)^{3/2}}

\newcommand{\cZcft}{\cZ_{\mathrm{cft}}}

\newcommand{\CCplus}[1]{\mathscr{C}_{+,#1}}
\newcommand{\CCminus}[1]{\mathscr{C}_{-,#1}}
\newcommand{\CCpm}[1]{\mathscr{C}_{\pm,#1}}

\renewcommand{\tilde}{\widetilde}

\newcommand{\kk}{k}
\newcommand{\ep}{e_+}
\newcommand{\fp}{f_+}
\newcommand{\emi}{e_-}
\newcommand{\fm}{f_-}
\newcommand{\epm}{e_\pm}
\newcommand{\fpm}{f_\pm}
\newcommand{\q}{\mathfrak{q}}
\newcommand{\qp}{\mathfrak{q}_+}
\newcommand{\qm}{\mathfrak{q}_-}
\newcommand{\qpm}{\mathfrak{q}_\pm}
\newcommand{\dkk}{\kappa}

\newcommand{\idem}{\boldsymbol{e}}
\newcommand{\Idem}{\boldsymbol{e}}

\newcommand{\vSE}{\boldsymbol{v_{\phantom{h}}^{\scriptscriptstyle\searrow}}}
\newcommand{\vNW}{\boldsymbol{v_{\phantom{h}}^{\scriptscriptstyle\nwarrow}}}
\newcommand{\vSW}{\boldsymbol{v_{\phantom{h}}^{\scriptscriptstyle\swarrow}}}
\newcommand{\vNE}{\boldsymbol{v_{\phantom{h}}^{\scriptscriptstyle\nearrow}}}

\newcommand{\DnDn}{{\scriptstyle\boldsymbol{\downdownarrows}}}
\newcommand{\UpUp}{{\scriptstyle\boldsymbol{\upuparrows}}}

\newcommand{\NESW}{{\scriptscriptstyle\boldsymbol{\nearrow\!\!\!\swarrow}}}
\newcommand{\NWSE}{{\scriptscriptstyle\boldsymbol{\nwarrow\!\!\kern-1pt\searrow}}}

\newcommand{\Radbullet}{\boldsymbol{v_{\phantom{h}}^{\scriptscriptstyle\bullet}}}
\newcommand{\wUp}{\boldsymbol{w^{\uparrow}}\!\!}
\newcommand{\wLeft}{\boldsymbol{w^{\scriptstyle\leftarrow}}}
\newcommand{\wRight}{\boldsymbol{w^{\scriptstyle\rightarrow}}}
\newcommand{\wDown}{\boldsymbol{w^{\downarrow}}}
\newcommand{\Nilpbullet}{\boldsymbol{w^{\scriptstyle\bullet}}}
\newcommand{\vUp}{\boldsymbol{v^{\uparrow}}\!\!}
\newcommand{\vLeft}{\boldsymbol{v^{\scriptstyle\leftarrow}}}
\newcommand{\vRight}{\boldsymbol{v^{\scriptstyle\rightarrow}}}

\newcommand{\piUp}{\boldsymbol{\pi^{\uparrow}}\!\!}
\newcommand{\piLeft}{\boldsymbol{\pi^{\scriptstyle\leftarrow}}}
\newcommand{\piRight}{\boldsymbol{\pi^{\scriptstyle\rightarrow}}}
\newcommand{\piDown}{\boldsymbol{\pi^{\downarrow}}}
\newcommand{\nilp}{\boldsymbol{w}}

\newcommand{\cRadSE}[1]{{b_{\phantom{h}#1}^{\scriptscriptstyle\searrow}}}
\newcommand{\cRadNW}[1]{{b_{#1}^{\scriptscriptstyle\nwarrow}}}
\newcommand{\cRadSW}[1]{{b_{#1}^{\scriptscriptstyle\swarrow}}}
\newcommand{\cRadNE}[1]{{b_{#1}^{\scriptscriptstyle\nearrow}}}
\newcommand{\cRadbullet}[1]{{b_{#1}^{\scriptscriptstyle\bullet}}}
\newcommand{\cNilpUp}[1]{{c_{#1}^{\uparrow}}\!\!}
\newcommand{\cNilpLeft}[1]{{c_{#1}^{\scriptstyle\leftarrow}}}
\newcommand{\cNilpRight}[1]{{c_{#1}^{\scriptstyle\rightarrow}}}
\newcommand{\cNilpDown}[1]{{c_{#1}^{\downarrow}}}
\newcommand{\cNilpbullet}[1]{{c_{#1}^{\scriptscriptstyle\bullet}}}
\newcommand{\cRadUp}[1]{{b_{#1}^{\uparrow}}\!\!}

\newcommand{\cRadRight}[1]{{b_{#1}^{\scriptstyle\rightarrow}}}

\newcommand{\prBullet}[1]{\mathsf{#1}^{\scriptstyle\bullet}}
\newcommand{\prUp}[1]{\mathsf{#1}^{\scriptstyle\uparrow}}
\newcommand{\prDown}[1]{\mathsf{#1}^{\scriptstyle\downarrow}}
\newcommand{\prRight}[1]{\mathsf{#1}^{\scriptstyle\rightarrow}}
\newcommand{\prLeft}[1]{\mathsf{#1}^{\scriptstyle\leftarrow}}

\newcommand{\nprUp}[2]{\mathsf{#1}^{\scriptstyle\uparrow,#2}}
\newcommand{\nprDown}[2]{\mathsf{#1}^{\scriptstyle\downarrow,#2}}
\newcommand{\nprRight}[2]{\mathsf{#1}^{\scriptstyle\rightarrow,#2}}
\newcommand{\nprLeft}[2]{\mathsf{#1}^{\scriptstyle\leftarrow,#2}}

\newcommand{\alphaUp}[1]{{\alpha^{\scriptstyle\uparrow,#1}}}

\newcommand{\alphaDown}[1]{{\alpha^{\scriptstyle\downarrow,#1}}}

\newcommand{\betaUp}[1]{{\beta^{\scriptstyle\uparrow,#1}}}

\newcommand{\betaDown}[1]{{\beta^{\scriptstyle\downarrow,#1}}}

\newcommand{\sbullet}{{\scriptstyle\bullet}}
\newcommand{\suparrow}{{\scriptstyle\uparrow}}
\newcommand{\sdownarrow}{{\scriptstyle\downarrow}}
\newcommand{\srightarrow}{{\scriptstyle\rightarrow}}
\newcommand{\sleftarrow}{{\scriptstyle\leftarrow}}
\newcommand{\snearrow}{{\scriptscriptstyle\nearrow}}
\newcommand{\snwarrow}{{\scriptscriptstyle\nwarrow}}
\newcommand{\ssearrowA}{{\scriptscriptstyle\searrow}}
\newcommand{\sswarrowA}{{\scriptscriptstyle\swarrow}}

\newcommand{\rightarrows}{\suparrow,\sdownarrow,\srightarrow,\sleftarrow}
\newcommand{\skewarrows}{\snearrow,\snwarrow,\ssearrowA,\sswarrowA}

\newcommand{\repLyu}{\pi}
\newcommand{\repmin}{\repLyu_{0}}
\newcommand{\repplus}{\repLyu_{(+)}}
\newcommand{\repminus}{\repLyu_{(-)}}

\newcommand{\captionfont}[1]{\textit{\textbf{\small#1}}} 

\newcommand{\XXX}{\boldsymbol{\mathfrak{g}}_{p_+,p_-}}
\newcommand{\WWW}{\boldsymbol{\mathscr{W}}_{p_+,p_-}}

\newcommand{\TA}{\boldsymbol{\mathscr{C}}}

\newcommand{\XX}{\mathsf{X}}
\newcommand{\XXgen}{\mathsf{X}}
\newcommand{\tXX}{\mathsf{{\widetilde{X}}}}
\newcommand{\VV}{\mathsf{V}}

\newcommand{\PP}{\mathsf{P}}
\newcommand{\NN}{\mathsf{N}}
\newcommand{\ZZ}{\mathsf{Z}}

\renewcommand{\hat}{\widehat}

\newcommand{\bref}[1]{\textbf{\textup{\ref{#1}}}}

\newcommand{\acts}{{\rightharpoondown}}

\newcommand{\cZ}{\mathsf{Z}}

\newcommand{\Grring}{\mathsf{G}}

\newcommand{\id}{\mathrm{id}}

\renewcommand{\geq}{\,{\geqslant}\,}
\renewcommand{\leq}{\,{\leqslant}\,}

\renewcommand{\le}{\,{\leqslant}\,}

\newcommand{\coup}[2]{\langle#1,#2\rangle} 

\newcommand{\tensor}{\otimes}

\newcommand{\End}{\mathrm{End}}

\newcommand{\fusion}{%
  \mathop{{\otimes}\kern-7pt\raisebox{.6pt}{%
      \mbox{\footnotesize${\bullet}$}}}}

\newcommand{\algA}{\mathfrak{A}}
\newcommand{\algB}{\mathfrak{B}}
\newcommand{\algI}{\mathfrak{I}}

\newcommand{\ket}[1]{|#1\rangle}
\newcommand{\bra}[1]{\langle#1|}

\newcommand{\mfrac}[2]{\mbox{\small$\displaystyle\frac{#1}{#2}$}}
\newcommand{\ffrac}[2]{\mbox{\footnotesize$\displaystyle\frac{#1}{#2}$}}
\newcommand{\half}{%
  \mathchoice{\ffrac{1}{2}}{\frac{1}{2}}{\frac{1}{2}}{\frac{1}{2}}}

\newcommand{\qint}[1]{[#1]}
\newcommand{\qbin}[2]{\mathchoice%
  {{\qbinm{#1}{#2}}}{\qbinmm{#1}{#2}}%
  {\qbinmm{#1}{#2}}{\qbinmm{#1}{#2}}}
\newcommand{\qbinm}[2]{\mbox{\footnotesize$\displaystyle
    \genfrac{[}{]}{0pt}{}{#1}{#2}$}}
\newcommand{\qbinmm}[2]{\genfrac{[}{]}{0pt}{}{#1}{#2}}

\newcommand{\NLLambda}{\boldsymbol{\Omega}}
\newcommand{\Nllambda}{\boldsymbol{\omega}}

\newcommand{\Nbarllambda}{\bar{\boldsymbol{\omega}}}

\newcommand{\Radford}{\widehat{\pmb{\phi}}{}}
\newcommand{\RadfordProj}{\widehat{\pmb{\varkappa}}{}}

\newcommand{\Rmin}{R_{\mathrm{min}}}
\newcommand{\Rproj}{R_{\mathrm{proj}}}

\newcommand{\RLeft}{R_{\boxbslash}}
\newcommand{\RRLeft}{\mathscr{R}_{\boxbslash}}
\newcommand{\RRight}{R_{\boxslash}}
\newcommand{\RRRight}{\mathscr{R}_{\boxslash}}

\newcommand{\setI}{\mathcal{I}}
\newcommand{\setILeft}{\mathcal{I}_{\boxslash}}
\newcommand{\setIRight}{\mathcal{I}_{\boxbslash}}
\newcommand{\setR}{\mathcal{I}_1}
\newcommand{\setii}{\mathscr{I}_1}

\newcommand{\Drinfeld}{\boldsymbol{\chi}}
\newcommand{\DrinfeldNESW}{\Drinfeld^{\NESW}}
\newcommand{\DrinfeldNWSE}{\Drinfeld^{\NWSE}}
\newcommand{\DrinfeldDnDn}{\Drinfeld^{\DnDn}}
\newcommand{\RadfordNESW}{\Radford^{\NESW}{}}
\newcommand{\RadfordNWSE}{\Radford^{\NWSE}{}}
\newcommand{\RadfordUpUp}{\Radford^{\UpUp}{}}

\newcommand{\vvarphi}{\hat{\pmb{\varphi}}{}}
\newcommand{\ppsi}{\hat{\pmb{\psi}}}

\newcommand{\vvarphiNESW}{\vvarphi^{\NESW}}
\newcommand{\vvarphiNWSE}{\vvarphi^{\NWSE}}

\newcommand{\hrhoB}{\widehat{\pmb{\rho}}{}^{\,\boxbslash}{}}
\newcommand{\RadfordB}{\vvarphi^{\boxbslash}{}}

\newcommand{\hrhoS}{\widehat{\pmb{\rho}}{}^{\,\boxslash}{}}
\newcommand{\RadfordS}{\vvarphi^{\boxslash}{}}

\newcommand{\Ngamma}{\gamma}
\newcommand{\NDrinfeldDnDn}{\Drinfeld^{\DnDn}}
\newcommand{\NRadfordNESW}{\Radford^{\NESW}{}}
\newcommand{\NRadfordNWSE}{\Radford^{\NWSE}{}}
\newcommand{\NRadfordUpUp}{\Radford^{\UpUp}{}}

\newcommand{\Nppsi}{\hat{\pmb{\psi}}}

\newcommand{\NRadfordB}{\vvarphi^{\boxbslash}{}}
\newcommand{\NRadfordS}{\vvarphi^{\boxslash}{}}

\newcommand{\NvvarphiJoin}{\vvarphi^{{\displaystyle\times}}}

\newcommand{\NvvarphiNWSE}{\vvarphi^{\NWSE}}
\newcommand{\NvvarphiNESW}{\vvarphi^{\NESW}}

\newcommand{\NDrinfeldNWSE}{\Drinfeld^{\NWSE}}
\newcommand{\NDrinfeldNESW}{\Drinfeld^{\NESW}}

\newcommand{\NsigmaNESW}{\sigma^{\NESW}}
\newcommand{\NsigmaNWSE}{\sigma^{\NWSE}}
\newcommand{\NsigmaUpUp}{\sigma^{\UpUp}}

\newcommand{\Nttheta}{\boldsymbol{\vartheta}}

\newcommand{\xxi}{\boldsymbol{\xi}}

\newcommand{\modS}{\mathscr{S}}
\newcommand{\modT}{\mathscr{T}}

\newcommand{\modPr}{\mathscr{P}}
\newcommand{\repX}{\rep{X}}

\newcommand{\modQ}{\rep{Q}}
\newcommand{\repR}{\rep{R}}

\newcommand{\rep}{\mathscr}

\newcommand{\dd}{\partial}
\newcommand{\SLiiZ}{SL(2,\oZ)}
\newcommand{\oC}{\mathbb{C}}

\newcommand{\oZ}{\mathbb{Z}}

\newcommand{\one}{\boldsymbol{1}}
\newcommand{\tr}{\mathrm{Tr}^{\vphantom{y}}}
\newcommand{\Tr}{\mathrm{Tr}^{\vphantom{y}}}

\newcommand{\qTr}{\mathrm{qCh}}
\newcommand{\ad}{\mathrm{Ad}}

\newcommand{\Ch}{\mathsf{Ch}}

\newcommand{\pbw}{\boldsymbol{e}}

\newcommand{\pbwd}{\boldsymbol{m}}
\newcommand{\pbwdd}{\boldsymbol{n}}

\newcommand{\comodul}{{\boldsymbol{a}}}
\newcommand{\coint}{{\boldsymbol{\Lambda}}}
\newcommand{\rint}{{\boldsymbol{\lambda}}}
\newcommand{\balance}{{\boldsymbol{g}}}
\newcommand{\sqs}{{\boldsymbol{u}}}

\newcommand{\ribbon}{{\boldsymbol{v}}}
\newcommand{\ribbonplus}{\ribbon^*_{(+)}}
\newcommand{\ribbonminus}{\ribbon^*_{(-)}}
\newcommand{\ribbonpm}{\ribbon^*_{(\pm)}}

\newcommand{\cas}{\boldsymbol{C}}

\newcommand{\cheb}{U}

\numberwithin{equation}{section}
\makeatletter
\@addtoreset{equation}{section}
\@addtoreset{subsubsection}{section}

\def\@secnumfont{\bfseries}
\def\subsubsection{\@startsection{subsubsection}{3}%
  \z@{.5\linespacing\@plus.7\linespacing}{-.5em}%
  {\normalfont\bfseries}}
\def\paragraph{\@startsection{paragraph}{4}%
  \z@\z@{-\fontdimen2\font}%
  \normalfont\bfseries}
\def\subparagraph{\@startsection{subparagraph}{5}%
  \z@\z@{-\fontdimen2\font}%
  \normalfont\bfseries}

\makeatother

\swapnumbers
\newtheorem{Thm}[subsection]{Theorem}
\newtheorem{thm}[subsubsection]{Theorem}

\newtheorem{lemma}[subsubsection]{Lemma}

\newtheorem{prop}[subsubsection]{Proposition}

\theoremstyle{definition}

\newtheorem{rem}[subsubsection]{Remark}

\begin{document}

\title[Logarithmic $(p,q)$ CFTs and quantum groups]{%
  \vspace*{-4\baselineskip}
  \mbox{}\hfill\texttt{\small\lowercase{math}.QA/\lowercase{0606506}}
  \\[\baselineskip]
  Kazhdan--Lusztig-dual quantum group for logarithmic extensions of
  Virasoro minimal models}

\author[Feigin]{B.L.~Feigin}%

\address{\mbox{}\kern-\parindent blf: Landau Institute for Theoretical
  Physics \hfill\mbox{}\linebreak \texttt{feigin@mccme.ru}}

\author[Gainutdinov]{A.M.~Gainutdinov}%

\address{\mbox{}\kern-\parindent ams, amg, iyt: Lebedev Physics
  Institute \hfill\mbox{}\linebreak \texttt{ams@sci.lebedev.ru},
  \texttt{azot@mccme.ru}, \texttt{tipunin@td.lpi.ru}}

\author[Semikhatov]{A.M.~Semikhatov}%

\author[Tipunin]{I.Yu.~Tipunin}

\begin{abstract}
  We derive and study a quantum group
  $\boldsymbol{\mathfrak{g}}_{p,q}$ that is Kazhdan--Lusztig-dual to
  the $W$-algebra $\boldsymbol{\mathscr{W}}_{p,q}$ of the logarithmic
  $(p,q)$ conformal field theory model.  The algebra
  $\boldsymbol{\mathscr{W}}_{p,q}$ is generated by two currents
  $W^+(z)$ and $W^-(z)$ of dimension $(2p-1)(2q-1)$, and the
  energy--momentum tensor $T(z)$.  The two currents generate a
  vertex-operator ideal $\boldsymbol{\mathscr{R}}$ with the property
  that the quotient
  $\boldsymbol{\mathscr{W}}_{p,q}/\boldsymbol{\mathscr{R}}$ is the
  vertex-operator algebra of the $(p,q)$ Virasoro minimal model.  The
  number ($2 p q$) of irreducible $\boldsymbol{\mathfrak{g}}_{p,q}$
  rep\-resenta\-tions is the same as the number of irreducible
  $\boldsymbol{\mathscr{W}}_{p,q}$-rep\-resenta\-tions on which
  $\boldsymbol{\mathscr{R}}$ acts nontrivially.  We find the center of
  $\boldsymbol{\mathfrak{g}}_{p,q}$ and show that the modular group
  representation on it is equivalent to the modular group
  representation on the $\boldsymbol{\mathscr{W}}_{p,q}$ characters
  and ``pseudocharacters.''  The factorization of the
  $\boldsymbol{\mathfrak{g}}_{p,q}$ ribbon element leads to a
  factorization of the modular group representation on the center.  We
  also find the $\boldsymbol{\mathfrak{g}}_{p,q}$ Grothendieck ring,
  which is presumably the ``logarithmic'' fusion of the $(p,q)$ model.
\end{abstract}

\maketitle
\enlargethispage{\baselineskip}

\thispagestyle{empty}


\setcounter{tocdepth}{2}

\vspace*{-24pt}

\begin{footnotesize}\addtolength{\baselineskip}{-6pt}
  \tableofcontents
\end{footnotesize}

\section{Introduction}
We consider quantum group counterparts of several important structures
manifested in logarithmic $(p,q)$ conformal field theory
models\,---\,examples of conformal field theory with nonsemisimple
representation categories, which have been the subject of some
attention from various
standpoints~\cite{[KeLu],[FHST],[My],[HLZ],[FGST],[FG],[FGST2],[jF]}.
The logarithmic $(p,q)$ models can be quite interesting (for the
$(3,2)$ model, there is evidence of its relation to the percolation
problem~\cite{[Cardy],[Watts],[DuplantierSaleur]}), but their
attractiveness is matched or, rather, counterbalanced by the
difficulty of studying them directly, which arises primarily because
the symmetries of such models are not the Virasoro algebra but its
nonlinear extensions, some $W$-algebras (and in the $(3,2)$ model, for
example, the leading pole in the operator product expansion of two
currents generating the $W$-algebra is already of order~$28$).  This
highlights the importance of ``indirect'' methods, which concentrate
around the Kazhdan--Lusztig correspondence~\cite{[KL]}.

The origin of the Kazhdan--Lusztig correspondence for logarithmic
models can already be seen in their definition as \textit{kernels} of
screening operators acting in a free-field space (we recall that the
rational minimal models are defined as the \textit{cohomology} of
screenings~\cite{[F],[BMP]}).  The Kazhdan--Lusztig-dual quantum
group~$\boldsymbol{\mathfrak{g}}$ is derived from the action of
screenings that makes the space of states in conformal field theory
into a \textit{bimodule}\,---\,a module over the $W$-algebra and the
quantum group~$\boldsymbol{\mathfrak{g}}$.  In the $(p,1)$ logarithmic
models, this leads to the equivalence of braided quasitensor
categories of the $W$-algebra representations and of the quantum group
representations~\cite{[FGST],[FGST2]}.

Quantum group counterparts can be found not only for irreducible
$W$-algebra repre\-sentations$/$characters but also for a larger space
$\TA$ of torus amplitudes.  We recall that in the nonsemisimple
(logarithmic) case, $\TA$ is spanned not only by the characters of
irreducible representations but also by ``generalized characters'' (or
``pseudo\-characters'') associated with some
\textit{pseudotraces}~\cite{[My]}
(cf.~\cite{[F-mod],[F-fusion],[FG]}).  On the quantum group side,
similarly, irreducible $\boldsymbol{\mathfrak{g}}$-modules provide
only a subspace in the space $\Ch(\mathfrak{g})$ of $q$-characters
($\ad^*$-invariant functionals on the quantum group); spanning all
of~$\Ch(\mathfrak{g})$ by elements associated with representations
requires taking pseudotraces associated with some indecomposable
$\boldsymbol{\mathfrak{g}}$-modules.

It can be expected in some generality that
\begin{equation*}
  \TA=\Ch(\mathfrak{g}),
\end{equation*}
i.e., the space of torus amplitudes and the space of $q$-characters of
the Kazhdan--Lusztig-dual quantum group are isomorphic.\pagebreak[3]
Such an isomorphism becomes much more interesting if it extends to
some structures defined on the respective spaces.  The modular group
action is a structure of major importance on $\TA$.  Remarkably, it is
nicely matched by the general theory of the modular group action on
ribbon factorizable quantum groups~\cite{[Lyu],[LM],[Kerler]} (also
see~\cite{[L-center]} for an application of this theory to the
\textit{small} quantum group).  For factorizable quantum groups, there
is an isomorphism (actually, of associative
algebras~\cite{[Drinfeld]})
\begin{equation}\label{S-acts-0}
  \Ch(\mathfrak{g})\xrightarrow{\ \Drinfeld\ }\cZ,
\end{equation}
which allows replacing $\Ch=\Ch(\mathfrak{g})$ with a somewhat more
``tangible'' object, the quantum group center $\cZ$. \ For the $(p,q)$
models, the $\SLiiZ$ representation~$\repLyu$ constructed on the
center of the Kazhdan--Lusztig-dual quantum group in accordance with
the recipe in~\cite{[Lyu],[LM],[Kerler]} turns out to be equivalent to
the $\SLiiZ$ representation on the torus amplitudes, as we show in
what follows (this has been known for the $(p,1)$ logarithmic
models~\cite{[FGST]}).

To formulate the result, we use the notation $p=p_+$ and $q=p_-$ for a
fixed pair of coprime positive integers and set\footnote{The paper is
  somewhat overloaded with the ``$\pm$-type'' notation, but it allows
  saving some space by ``$\pm$'' formulas.}
\begin{gather}\label{small-q}
  \q=e^{\frac{i\pi}{2p_+p_-}},\qquad
  \qp=\q^{2p_-}=e^{\frac{i\pi}{p_+}},\qquad  
  \qm=\q^{2p_+}=e^{\frac{i\pi}{p_-}}.
\end{gather}

\textit{On the conformal field theory side}, the chiral algebra of the
logarithmic $(p_+,p_-)$ conformal field theory is a certain
``triplet'' $W$-algebra $\WWW$ defined in~\cite{[FGST3]}.  The
characters of its irreducible representations give rise to a
$\frac{1}{2}(3p_+\,{-}\,1)(3p_-\,{-}\,1)$-dimensional representation
$\cZcft$ of the modular group $\SLiiZ$, evaluated as~\cite{[FGST3]}
\begin{equation}\label{W-structure}
  \cZcft
  =\Rmin\oplus\Rproj\oplus\oC^2\tensor(\RLeft\oplus\RRight)\oplus
  \oC^3\tensor\Rmin,
\end{equation}
where $\oC^2$ is the standard two-dimensional representation, $\oC^3$
is its symmetrized square, $\Rmin$ is the
$\frac{1}{2}(p_+\,{-}\,1)(p_-\,{-}\,1)$-dimensional $\SLiiZ$
representation on the characters of the $(p_+,p_-)$ Virasoro minimal
model, and $\Rproj$, $\RLeft$, and $\RRight$ are certain $\SLiiZ$
representations of the respective dimensions
$\frac{1}{2}(p_+\,{+}\,1)(p_-\,{+}\,1)$,
$\frac{1}{2}(p_+\,{+}\,1)(p_-\,{-}\,1)$, and
$\frac{1}{2}(p_+\,{-}\,1)(p_-\,{+}\,1)$ (see~\cite{[FGST3]} for the
precise formulas).

\textit{On the quantum group side}, the quantum group $\XXX$ that is
Kazhdan--Lusztig-dual to the $(p_+,p_-)$ logarithmic conformal field
theory is the quotient
\begin{equation}\label{as-quotient}
  \XXX=\frac{\overline{\mathscr{U}}_{Q_{+}}s\ell(2)
    \tensor
    \overline{\mathscr{U}}_{Q_{-}}s\ell(2)}{(K_+^{p_+} - K_-^{p_-})}
\end{equation}
of the product of two restricted quantum groups
$\overline{\mathscr{U}}_{Q_{\pm}}s\ell(2)$ at the roots of unity
$Q_{\pm}=\q_{\pm}^{p_{\mp}}$ over the (Hopf) ideal generated by the
central element $K_+^{p_+} - K_-^{p_-}$.
\enlargethispage{\baselineskip}

\begin{Thm}\label{thm:q-structure}
  \mbox{}

  \begin{enumerate}
  \item The center $\cZ$ of $\XXX$ is
    $\frac{1}{2}(3p_+\!-\!1)(3p_-\!-\!1)$-dimen\-sional\textup{.}

  \item\label{claim:factorization} The $\SLiiZ$ representation
    $\repLyu$ on $\cZ$ factors into the product of three $\SLiiZ$
    representations,
    \begin{equation*}
      \repLyu(A)=\repmin(A)\,\repplus(A)\,\repminus(A),
      \qquad A\in\SLiiZ,
    \end{equation*}
    which pairwise commute \textup{(}i.e.,
    $\repmin(A)\,\repplus(A')=\repplus(A')\,\repmin(A)$, \
    $\repmin(A)$\linebreak[0]$\repminus(A')
    $\linebreak[0]${}=\repminus(A')$\linebreak[0]$\repmin(A)$, and
    $\repplus(A)\,\repminus(A')
    =\repminus(A')\,\repplus(A)$\textup{).}

  \item\label{claim:decomp}The center decomposes with respect to this
    product as
    \begin{equation*}
      \cZ=
      \Rproj\tensor 1 \tensor 1
      \oplus
      \RLeft\tensor\oC^2\tensor 1
      \oplus
      \RRight\tensor 1 \tensor\oC^2
      \oplus
      \Rmin\tensor\oC^2\tensor\oC^2,
    \end{equation*}
    where $\oC^2$ is the standard two-dimensional representation, $1$
    is the trivial representation, $\Rmin$ is the
    $\frac{1}{2}(p_+\!-\!1)(p_-\!-\!1)$-dimensional $\SLiiZ$
    representation on the characters of the $(p_+,p_-)$ Virasoro
    minimal model, and $\Rproj$, $\RLeft$, and $\RRight$ are those
    in~\eqref{W-structure}\textup{.}

  \item\label{claim:grring}The representation
    \begin{equation*}
      \Rproj
      \oplus
      \RLeft
      \oplus
      \RRight
      \oplus
      \Rmin
    \end{equation*}
    is the $\SLiiZ$ representation on the Grothendieck ring
    of~$\XXX$\textup{.}
    
  \item\label{claim:equiv} Moreover, the decomposition
    $\oC^2\tensor\oC^2=1\oplus\oC^3$ establishes an equivalence of the
    $\SLiiZ$ representation on $\cZ$ to the representation on the
    space of generalized $\WWW$-characters in~\eqref{W-structure}.
  \end{enumerate}
\end{Thm}

The $\XXX$ quantum group has $2 p_+ p_-$ irreducible representations,
which we label as $\XX^\pm_{r,r'}$ with $1\leq r\leq p_+$ and $1\leq
r'\leq p_-$.
\begin{Thm}\label{thm:Gr-ring}
  Multiplication in the $\XXX$ Grothendieck ring is given by
  \begin{equation}\label{the-Groth}
    \XX^{\alpha}_{r,r'}\XX^{\beta}_{s,s'}
    =\sum_{\substack{u=|r - s| + 1\\
        \mathrm{step}=2}}^{r + s - 1}
    \sum_{\substack{u'=|r' - s'| + 1\\
        \mathrm{step}=2}}^{r' + s' - 1}
    {\tXX}^{\alpha\beta}_{u,u'},
  \end{equation}
  where
  \begin{equation*}
    {\tXX}^{\alpha}_{r,r'} =
    \begin{cases}
      \XX^{\alpha}_{r,r'},& 
      \begin{array}[t]{l}
        1\leq r\leq p_+,\\[-6pt]
        1\leq r'\leq p_-,
      \end{array}
      \\
      \XX^{\alpha}_{2p_+ - r,r'} + 2\XX^{-\alpha}_{r - p_+, r'},&
      \begin{array}[t]{l}
        p_+\!+\!1\leq r\leq 2 p_+\!-\!1,\\[-6pt]
        1\leq r'\leq p_-,
      \end{array}
      \\
      \XX^{\alpha}_{r,2p_- - r'} + 2\XX^{-\alpha}_{r,r' - p_-},&
      \begin{array}[t]{l}
        1\leq r\leq p_+,\\[-6pt]
        p_-\!+\!1\leq r'\leq 2 p_-\!-\!1,
      \end{array}
      \\
      \mbox{}\kern-3pt\begin{aligned}[b]
        &\XX^{\alpha}_{2p_+ - r, 2 p_- - r'}
        + 2\XX^{-\alpha}_{2p_+ - r, r' - p_-}\\[-3pt]
        &{}+ 2\XX^{-\alpha}_{r - p_+, 2 p_- - r'}
        + 4\XX^{\alpha}_{r - p_+, r' - p_-},
      \end{aligned}&
      \begin{array}[t]{l}
        p_+\!+\!1\leq r\leq 2 p_+\!-\!1,\\[-6pt]
        p_-\!+\!1\leq r'\leq 2 p_-\!-\!1.
      \end{array}
    \end{cases}
  \end{equation*}
\end{Thm}
We also characterize this ring as a quotient of $\oC[x,y]$
(see~\bref{prop:Cheb}).  \ Its interpretation as a $\WWW$ fusion ring
is discussed in~\cite{[FGST3]}.

\medskip

To derive the Kazhdan--Lusztig-dual quantum group for the logarithmic
$(p_+,p_-)$ conformal field theory models, we follow the strategy
proposed in~\cite{[FGST]} for $(p,1)$ models.  The starting point is
the Hopf algebra $\mathscr{H}$ of screening operators in a free-field
space.  In the $(p,1)$ case, $\mathscr{H}$ is a Taft Hopf algebra at a
certain root of unity, and in the $(p_+,p_-)$ case, $\mathscr{H}$ is
constructed from two Taft algebras (taken at different roots of
unity).  We then pass to the Drinfeld double $D(\mathscr{H})$.  {}From
the conformal field theory standpoint, the double is a quantum group
generated by the screening and contour-removal operators, the latter
being just the dual to the former in Drinfeld's construction.  But the
doubling procedure also yields the dual $\dkk$ to the Cartan element
$\kk$ in $\mathscr{H}$, and it is to be ``eliminated'' as
$\dkk=\kk^{-1}$; this amounts to taking a quotient~$\bar D$
of~$D(\mathscr{H})$.

As any Drinfeld double, $D(\mathscr{H})$ is quasitriangular, i.e., has
a universal $R$-matrix $R\in D(\mathscr{H})\tensor D(\mathscr{H})$.
It follows that $\bar D$ is still quasitriangular, but the $M$-matrix
$\bar M = {}^t\bar R \bar R$ turns out to be an element not of $\bar
D\tensor\bar D$ but of $\boldsymbol{\mathfrak{g}}_{p_+,p_-}\tensor
\boldsymbol{\mathfrak{g}}_{p_+,p_-}$, where (for the $(p_+,p_-)$ model
specifically) $\boldsymbol{\mathfrak{g}}_{p_+,p_-}$ is the subalgebra
in $\bar D$ generated by $K\equiv\kk^2$ and the other $\bar D$
generators.  This $\boldsymbol{\mathfrak{g}}_{p_+,p_-}$ can be
alternatively described as in~\eqref{as-quotient}.

This Kazhdan-Lusztig-dual quantum group
$\boldsymbol{\mathfrak{g}}_{p_+,p_-}$ is \textit{not} quasitriangular,
but we have a pair $\boldsymbol{\mathfrak{g}}_{p_+,p_-}\subset\bar D$
of quantum groups, one of which has an $M$-matrix but not an
$R$-matrix, and the other has an $R$-matrix but not a (nondegenerate)
$M$-matrix.  Abusing the terminology, we still refer to the
$\boldsymbol{\mathfrak{g}}_{p_+,p_-}$ quantum group as factorizable
(although it, e.g., is not unimodular).  Such a ``deficient
factorizability'' suffices for constructing a modular group action on
the $\boldsymbol{\mathfrak{g}}_{p_+,p_-}$ center~$\cZ$.  In accordance
with~\cite{[Lyu],[LM],[Kerler]}, the action of
$\modS=\left(\begin{smallmatrix}
    0&1\\
    -1&0
  \end{smallmatrix}\right)\in\SLiiZ$ on $\cZ$ 
is constructed in terms of the Drinfeld and Radford maps $\Drinfeld$
and $\Radford$ and is given by composing one of these with the inverse
of the other,
\begin{equation}\label{3-diagram}
  \xymatrix@=16pt{
    {}&\Ch\ar[dl]_{\Drinfeld}
    \ar[dr]^{\Radford}
    \\
    \cZ
    &{}&\cZ\ar@/^/[ll]^{\modS^{-1}}
    }
\end{equation}
(the Radford map $\Radford$ is the inverse compared with the standard
definition).  Next, $\modT=\left(\begin{smallmatrix}
    1&1\\
    0&1
  \end{smallmatrix}\right)\in\SLiiZ$ is represented on $\cZ$
essentially by (multiplication with) the ribbon element,
\begin{equation}
  \label{eq:T-diagram}
  \cZ\xrightarrow{\quad\ribbon\quad}\cZ.
\end{equation}
Here, too, the properties of
$\boldsymbol{\mathfrak{g}}_{p_+,p_-}\subset\bar D$ are such that the
ribbon element~$\ribbon$ actually belongs to
$\boldsymbol{\mathfrak{g}}_{p_+,p_-}$ (not surprisingly, in view of
the relation between the $M$-matrix and the ribbon
element).

We prove Theorem~\bref{thm:q-structure} rather directly, by an
explicit calculation relying on the existence of two special bases in
the center, associated with the Drinfeld and Radford maps.  Some
``obvious'' central elements that are related by $\modS$ are given by
the Radford and Drinfeld images of the (``quantum'') traces over
irreducible representations.\footnote{From the standpoint of
  \textit{boundary} conformal field theory, the $\Radford$ map of the
  basis of characters and pseudocharacters in~$\Ch$ yields the
  Ishibashi states, i.e., states with a certain ``momentum.''  The
  $\Drinfeld$ map of the same elements in $\Ch$ yields the localized,
  i.e., Cardy states~\cite{[Ca89]}.}  Adding some pseudotraces yields
the ``Radford'' and ``Drinfeld'' bases in~$\cZ$, with $\modS$ mapping
between them.\footnote{The Drinfeld and Radford-map images of only the
  (traces over) irreducible representations do not span the entire
  center; this is an essential complication compared with the quantum
  $s\ell(2)$ in~\cite{[L-center]} and in~\cite{[FGST]}.}  The theorem
eventually follows from studying how the ribbon element acts on the
two bases.

The relation between the $\SLiiZ$ representations in
Claims~\ref{claim:decomp} and~\ref{claim:grring} can be stated in
terms of an automorphy factor, similarly to~\cite{[FHST]}.  The
factorization itself is related to the multiplicative Jordan
decomposition of the ribbon element but goes somewhat further: the
ribbon element factors as
\begin{equation*}
  \ribbon = 
  \Bigl(\sum_{(r,r')\in\setI}\!\!
    e^{2i\pi\Delta_{r,r'}}\Idem(r,r')\Bigr)
  \Bigl(\one
  + \ffrac{1}{p_+}\,
  \NDrinfeldNWSE(1, 1)\Bigr)
  \Bigl(\one
  + \ffrac{1}{p_-}\,
  \NDrinfeldNESW(1, 1)
  \Bigr).
\end{equation*}
The first factor is the semisimple part, constructed from the
primitive idempotents $\Idem(r,r')$ (weighted by exponentials of the
conformal dimensions of $\WWW$ primary fields in the $(p_+,p_-)$
logarithmic conformal field theory model); in the other two factors,
which are unipotent, $\DrinfeldNWSE(1, 1)$ and $\DrinfeldNESW(1, 1)$
are elements in the radical given by Drinfeld-map images of certain
pseudotraces associated with indecomposable $\XXX$ representations in
the full subcategory containing the trivial representation.

\medskip

The rest of the paper is not difficult to describe, despite its
length.  We recall how the quantum group is derived from conformal
field theory data in~\bref{sec:from} and describe its representations
in~\bref{sec:repr} (and Appendix~\ref{app:priojective}),
$q$-characters in~\bref{sec:all-q-chars}, and the center
in~\bref{sec:center}.  We then travel along the diagram
in~\eqref{3-diagram} from top down, considering the $\Radford$ arrow
in Sec.~\ref{sec:radford} and the $\Drinfeld$ arrow in
Sec.~\ref{sec:drinfeld}.  The Radford and Drinfeld maps $\Radford$ and
$\Drinfeld$ give two ``$\modS$-adapted'' bases in the center.  It then
remains to study how $\modT\in\SLiiZ$ acts on these bases.  This
amounts to studying the action of the ribbon element derived
in~\bref{XXX-ribbon}.  We prove Theorem~\bref{thm:q-structure} in
Sec.~\ref{sec:SL2Z}.

The basic ingredients for the $\Radford$ map (and its inverse), the
integral and cointegral, are found in~\bref{sec:int}; for the
$\Drinfeld$ map, the required $M$-matrix is derived from the Drinfeld
double in~\bref{sec:the-M}.  The analysis of the $q$-characters $\Ch$
and the center $\cZ$ relies on the study of representations, with an
important role played by projective and some other indecomposable
modules.  These are studied in Appendix~\ref{app:priojective}. In
Appendix~\ref{app:the-center}, we detail the construction of the
$\XXX$ center and collect some explicit calculations.

\subsection*{Notation} In $\oZ[q,q^{-1}]$, we consider the
$q$-integers and $q$-binomial coefficients
\begin{equation*}
  \qint{n}_q = \ffrac{q^{n}-q^{-n}}{q-q^{-1}},\quad
  {\qbin{m}{n}}_q
  =\ffrac{\qint{m}_q!}{\qint{n}_q!\,\qint{m-n}_q!},\quad
  \qint{n}_q! = \qint{1}_q\qint{2}_q\dots
  \qint{n}_q.
\end{equation*}
We then let
\begin{equation}
  \qint{m}_+=\qint{m}_{\qp^{p_-}},\quad
  \qint{m}_-=\qint{m}_{\qm^{p_+}},\quad
  {\qbin{m}{n}}_+={\qbin{m}{n}}_{\qp^{p_-}},\quad
  {\qbin{m}{n}}_-={\qbin{m}{n}}_{\qm^{p_+}}
\end{equation}
denote their respective specializations at $\q_{\pm}^{p_{\mp}}$
(see~\eqref{small-q}).

\begin{figure}[htb]
   \centering
   \includegraphics[bb=1.5in 8.4in 7in 10.6in, clip]{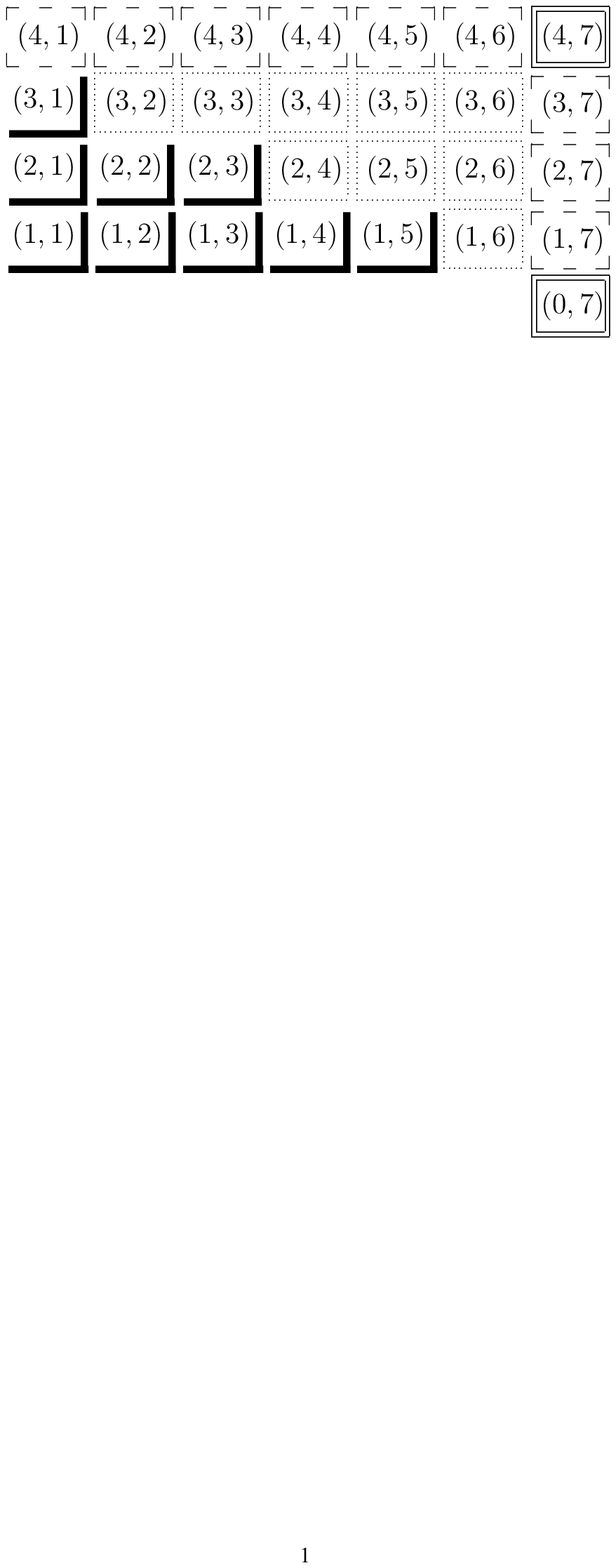}
  \caption[The four sets of pairs]{\captionfont{The four sets of
      pairs}.  \small The $\begin{xy}
      0*=<29pt,14pt>{(r,r')}*\frm<2pt>{,}
    \end{xy}$ \ boxes constitute the set $\setR$.  The column of
    $\begin{xy} 0*=<35pt,16pt>{(r,p_-)}*\frm{--}
    \end{xy}$ \ boxes is added to make $\setILeft$, and the row of
    $\begin{xy} 0*=<37pt,16pt>{(p_+,r')}*\frm{--}
    \end{xy}$ \ boxes to make $\setIRight$.  Adding both the row and
    the column and the two $\begin{xy}
      0*=<27pt,16pt>{(~\;,\;)}*\frm{=}
    \end{xy}$ \ boxes makes the set~$\setI$.}
  \label{fig:Kac-table}
\end{figure}

The fixed coprime pair $p_+$ and $p_-$ determines the four sets of
pairs
\begin{alignat}{2}
  \label{set-R}
  \setR&=
  \{(r, r')\mid 1\leq r\leq p_+ - 1,\  1\leq r'\leq p_- - 1,\ 
  p_- r + p_+ r' \leq p_+ p_-\}
  \kern-260pt\\
  \intertext{with $\vert\setR\vert=\half(p_+\!-\!1)(p_-\!-\!1)$,}
  \setILeft&=\setR\cup
  \{(r,p_-)\mid1\leq r\leq p_+\!-\!1\},&
  &\vert\setILeft\vert=\half(p_+\!-\!1)(p_-\!+\!1),
  \\
  \setIRight&=\setR\cup
  \{(p_+,r')\mid1\leq r'\leq p_-\!-\!1\},&
  \qquad
  &\vert\setIRight\vert=\half(p_+\!+\!1)(p_-\!-\!1),
\end{alignat}
and
\begin{multline}
  \setI=\setR
  \cup\{(r,p_-)\mid1\leq r\leq p_+\!-\!1\}
  \cup\{(p_+,r')\mid1\leq r'\leq p_-\!-\!1\}\\*
  {}\cup\{(p_+,p_-)\}\cup\{(0,p_-)\}
  \label{eq:setI}  
\end{multline}
with $|\setI|
=\half(p_+\!+\!1)(p_-\!+\!1)$.  The sets $\setR$, $\setILeft$, and
$\setIRight$ are subsets of the extended Kac table, a $p_+\times p_-$
rectangle, and $\setI$ additionally contains the box conventionally
labeled by $(0,p_-)$ (it can be considered a ```second copy'' of the
Steinberg box $(p_+,p_-)$).  For $p_+=4$ and $p_-=7$, these are
illustrated in Fig.~\ref{fig:Kac-table}.

\section{The $\smash{\XXX}$ quantum group, its representations,
  $q$-characters, and center}\label{sec:2} The quantum group
corresponding to the logarithmic $(p_+,p_-)$ conformal field theory
model, $\XXX$ at an appropriate root of unity, is introduced
in~\bref{sec:from}.  We then consider its representations
in~\bref{sec:repr}, describe the space of $q$-characters
in~\bref{sec:all-q-chars} and give the structure of the $\XXX$ center
in~\bref{sec:center}.

We refer the reader
to~\cite{[LSw],[Rad-antipode],[Drinfeld],[Kassel],[ChP]} for the
fundamental facts regarding Hopf algebra structures.

\subsection{$\smash{\XXX}$ from the double of the ``two-screening''
  quantum group~\cite{[FGST3]}}\label{sec:from} A ``gateway'' between
the direct study of the logarithmic $(p_+,p_-)$ model and its
quantum group counterpart is provided by the screening operators $\ep$
and $\fm$ and a ``Cartan'' operator $\kk$ constructed from the
free-field zero mode~\cite{[FGST3]}.  They act on the chiral sector of
the space of states in the logarithmic model such that the relations
\begin{equation}\label{B-relations}
  \begin{gathered}
    \ep^{p_+}=\fm^{p_-}=0,\quad\kk^{4p_+p_-}=\one,
    \quad\ep\fm=\fm\ep,\\
    \kk\ep\kk^{-1}=\qp\ep,\quad\kk\fm\kk^{-1}=\qm^{-1}\fm
  \end{gathered}
\end{equation}
are satisfied.  Let $\mathscr{H}$ denote the Hopf algebra with these
generators and relations; its comultiplication, antipode, and counit
(obtained by naturally combining those for the two Taft Hopf algebras
generated by $(\ep,\kk)$ and $(\fm,\kk)$) are given in~\cite{[FGST3]}.
The PBW basis in~$\mathscr{H}$ is
\begin{equation}\label{ejmn}
  \pbw_{jmn}=\kk^j\ep^m\fm^n,
  \qquad 0\leq j\leq 4p_+p_--1,
  \quad 0\leq m\leq p_+ -1,
  \quad 0\leq n\leq p_- -1.
\end{equation}

We next construct the Drinfeld double $D(\mathscr{H})$
of~$\mathscr{H}$.  This involves introducing generators $\dkk$, $\fp$,
and $\emi$, ``dual'' to $\kk$, $\ep$, and $\fm$, by the relations
\begin{equation}\label{basis}
  \begin{gathered}
    \coup{\dkk}{\pbw_{jmn}}=\delta_{m,0}\delta_{n,0}\q^{j},\\
    \coup{\fp}{\pbw_{jmn}}
    =-\delta_{m,1}\delta_{n,0}\ffrac{\qp^{j}}{\qp^{p_-}\!- \qp^{-p_-}},
    \quad
    \coup{\emi}{\pbw_{jmn}}
    =-\delta_{m,0}\delta_{n,1}\ffrac{\qm^{-j}}{\qm^{p_+}\!- \qm^{-p_+}}
  \end{gathered}
\end{equation}
Further, we eliminate $\dkk$, which is extraneous from the conformal
field theory standpoint, by taking the quotient
\begin{equation}\label{D-bar}
  \bar D(\mathscr{H})=D(\mathscr{H})\!\bigm/\!(\kk\dkk-\one)
\end{equation}
over the \textit{Hopf} ideal generated by the central element
$\kk\dkk-\one$.  We then let $\XXX$ be the subalgebra in $\bar
D(\mathscr{H})$ generated by (the images of) $\ep$, $\fp$, $\emi$,
$\fm$, and $K=k^2$.  It is shown in~\cite{[FGST3]} that this it can be
equivalently described as follows.

\begin{prop}\label{XXX-def}
  $\XXX$ is the Hopf algebra generated by $\ep$, $\fp$, $\emi$, $\fm$,
  and $K$ with the relations
  \begin{gather*}
    \epm^{p_\pm}=\fpm^{p_\pm}=0,\quad K^{2p_+p_-}=\one,\\
    K\epm K^{-1}=\qpm^2\epm,\quad K\fpm K^{-1}=\qpm^{-2}\fpm,\\
    \ep\emi=\emi\ep,\quad\fp\fm=\fm\fp,
    \quad\ep\fm=\fm\ep,\quad
    \emi\fp=\fp\emi,\\
    {[}\ep,\fp]
    =\ffrac{K^{p_-} - K^{-p_-}}{\qp^{p_-}\!- \qp^{-p_-}},\qquad
    {[}\emi,\fm]
    =\ffrac{K^{p_+} - K^{-p_+}}{\qm^{p_+}\!- \qm^{-p_+}}\\
    \intertext{and the Hopf algebra structure}
    \Delta(K)= K\tensor K,
    \quad
    \Delta(\ep)=\ep\tensor\one+ K^{p_-}\tensor\ep,
    \quad
    \Delta(\fm)=\fm\tensor\one + K^{-p_+}\tensor\fm,
    \\
    \Delta(\fp)=\fp\tensor K^{-p_-}+\one\tensor\fp,
    \quad
    \Delta(\emi)=\emi\tensor K^{p_+}+\one\tensor\emi,
    \\
    S(K)= K^{-1},\quad
    S(\ep) =-K^{-p_-}\ep, \quad S(\fm)=-K^{p_+}\fm,\\
    S(\fp)=-\fp K^{p_-},
    \quad S(\emi)=-\emi K^{-p_+},\\
    \epsilon(\epm)=\epsilon(\fpm)=0,\quad\epsilon(K)=1.
  \end{gather*}
\end{prop}

\subsubsection{} Furthermore, $\XXX$ can be described as the
quotient~\eqref{as-quotient} of the product of two restricted quantum
groups $\overline{\mathscr{U}}_{Q_{\pm}}s\ell(2)$ with
$Q_{\pm}=\q_{\pm}^{p_{\mp}}$ over the Hopf ideal generated by the
central element $K_+^{p_+} - K_-^{p_-}$.  By $\overline{\mathscr{U}}_q
s\ell(2) \equiv\overline{\mathscr{U}}_q\langle e,f,K;p\rangle$ at
$q^{2p}=1$, we mean the quantum group with the relations $e^{p}=
f^{p}=0$ and $K^{2 p}=\one$, the standard relations
\begin{gather*}
  K e  K^{-1}= q^2 e ,\quad K f K^{-1}=q^{-2} f,\quad
  [ e, f] = \ffrac{K - K^{-1}}{q\!- \q^{-1}},
  \\
  \intertext{and the Hopf algebra structure given by}
  \Delta(K)= K\tensor K,
  \quad
  \Delta( e)= e \tensor\one + K\tensor e,
  \quad
  \Delta(f)=f\tensor K^{-1} + \one\tensor f,
  \\
  S(K)= K^{-1},\quad
  S(e) = -K^{-1} e,\quad
  S(f)=- f K.  
\end{gather*}

Indeed, because $p_+$ and $p_-$ are coprime, the $\XXX$ algebra
in~\bref{XXX-def} is equivalently generated by $\ep$, $\fp$, $\emi$,
$\fm$, and the elements
\begin{equation}\label{Kpm}
  K_+=K^{p_-},\qquad K_-=K^{p_+},
\end{equation}
{}From the formulas in~\bref{XXX-def}, we then have
\begin{gather*}
  K_{\pm}e_{\pm}=\qpm^{2p_{\mp}}e_{\pm}K_{\pm},
  \quad
  K_{\pm}f_{\pm}=\qpm^{-2p_{\mp}}f_{\pm}K_{\pm},\quad
  [\epm,\fpm] =\ffrac{K_{\pm} - K_{\pm}^{-1}}{
    \q_{\pm}^{p_{\mp}}\!- \q_{\pm}^{-p_{\mp}}},
\end{gather*}
and it follows in addition that $K_{\pm}e_{\mp}=e_{\mp}K_{\pm}$ and
$K_{\pm}f_{\mp}=f_{\mp}K_{\pm}$.  The formulas for $\Delta$ and $S$
also ``separate,'' e.g., $S(\ep) =-K_+^{-1}\ep$.  The two
$\overline{\mathscr{U}}_{q}s\ell(2)$ algebras are therefore
$\overline{\mathscr{U}}_{\qp^{p_-}}\langle \ep,\fp,K_+;p_+\rangle$ and
$\overline{\mathscr{U}}_{\qm^{p_+}}\langle
\fm,\emi,K_-^{-1};p_-\rangle$.  Clearly, the relations
$K_+^{2p_+}=K_-^{2p_-}=K_+^{p_+}K_-^{p_-}=\one$ are satisfied by
$K_\pm$ in~\eqref{Kpm}.

\subsection{$\smash{\XXX}$-modules}\label{sec:repr} In what
follows, we need the irreducible and projective $\XXX$-modules; on the
way from the former to the latter, useful intermediate objects are the
Verma modules.

\subsubsection{Irreducible modules}\label{subsec:irrep} It is
easy to see that there are $2p_+ p_-$ irreducible finite-dimensional
$\XXX$-modules.  We label them as $\XX^{\pm}_{r,r'}$, where $1\leq
r\leq p_+$ and $1\leq r'\leq p_-$, with the highest-weight vector
$\ket{r,r'}^\pm$ of $\XX^{\pm}_{r,r'}$ defined by the relations
\begin{gather*}
  \ep\ket{r,r'}^\pm=\emi\ket{r,r'}^\pm=0,\quad
  K\ket{r,r'}^\pm=\pm\qp^{r-1}\qm^{r'-1}\ket{r,r'}^\pm.
\end{gather*}
It follows that $\dim\XX^{\pm}_{r,r'}=r r'$. \ $\XX^{+}_{1,1}$ is the
trivial module.  The module $\XX^{\alpha}_{r,r'}$ is linearly spanned
by elements $\ket{r,r',n,n'}^{\alpha}$, $0\leq n\leq r\,{-}\,1$ and
$0\leq n'\leq r'\,{-}\,1$ (with
$\ket{r,r',0,0}^{\pm}\equiv\ket{r,r'}^{\pm}$), with the $\XXX$-action
given by
\begin{align}
  K\, \ket{r,r',n,n'}^{\alpha}&=
  \alpha\qp^{r-1-2n}\,\qm^{r'-1-2n'}\ket{r,r',n,n'}^{\alpha},
  \notag\\
  \ep\, \ket{r,r',n,n'}^{\alpha}&= \alpha^{p_-}(-1)^{r'-1}
  \qint{n}_+\qint{r\!-\!n}_+\ket{r,r',n\!-\!1,n'}^{\alpha},
  \notag\\
  \emi\, \ket{r,r',n,n'}^{\alpha}&=\alpha^{p_+}(-1)^{r-1}
  \qint{n'}_-\qint{r'\!-\!n'}_-\ket{r,r',n,n'\!-\!1}^{\alpha},
  \label{on-vectors}\\
  \fp\, \ket{r,r',n,n'}^{\alpha}&=\ket{r,r',n\!+\!1,n'}^{\alpha},
  \notag\\
  \fm\, \ket{r,r',n,n'}^{\alpha}&=\ket{r,r',n,n'\!+\!1}^{\alpha},
  \notag
\end{align}
where we set $\ket{r,r',r,n'}^\alpha=\ket{r,r',n,r'}^\alpha
=\ket{r,r',-1,n'}^\alpha=\ket{r,r',n,-1}^\alpha=0$.

\subsubsection{Verma modules}\label{subsec:verma-mod}
There are $2p_+ p_-$ Verma modules, labeled as $\VV^{\pm}_{r,r'}$ with
$1\!\leq r\!\leq p_+$ and $1\!\leq r'\!\leq p_-$.  First, these are
the two Steinberg modules
\begin{equation*}
  \VV^{\pm}_{p_+,p_-}=\XX^{\pm}_{p_+,p_-}.
\end{equation*}
Next, for each $r=1,\dots,p_+\,{-}\,1$, $r'=1,\dots,p_-\,{-}\,1$, and
$a=\pm$, the Verma modules $\VV^{a}_{r,p_-}$ and $\VV^{a}_{p_+,r'}$
can be described as the respective extensions
\begin{align*}
  &0\to\XX^{-a}_{p_+-r,p_-}\to\VV^{a}_{r,p_-}\to\XX^{a}_{r,p_-}\to0,\\
  &0\to\XX^{-a}_{p_+,p_--r'}\to\VV^{a}_{p_+,r'}\to\XX^{a}_{p_+,r'}\to0.
\end{align*}
For consistency with the notation used in more complicated extensions
below, we also write these extensions as
\begin{equation*}
  \overset{\XX^{a}_{r,p_-}}{\bullet}{}\longrightarrow{}
  \overset{\!\!\!\XX^{-a}_{p_+-r,p_-}\!\!\!}{\bullet},\quad\quad\quad
  \overset{\XX^{a}_{p_+,r'}}{\bullet}{}\longrightarrow{}
  \overset{\!\!\!\XX^{-a}_{p_+,p_--r'}\!\!\!}{\bullet},
\end{equation*}
with the convention that arrows are directed toward
\textit{sub}modules.

Next, for $1\leq r\leq p_+-1$ and $1\leq r'\leq p_--1$, the Verma
module $\VV^{a}_{r,r'}$ can be described as a second extension with
the following structure of subquotients:\footnote{In diagrams of this
  type, it is understood that the $\XXX$-action on each irreducible
  representation is changed in agreement with the arrows connecting a
  given subquotient with others.  For example, the lowest-weight
  vector of $\XX^{a}_{r,r'}$ is annihilated by $\fm$ in the
  \textit{irreducible} representation, but is mapped by $\fm$ into the
  highest-weight vector in $\smash{\XX^{-a}_{r,p_--r'}}$
  in~\eqref{schem-Verma}.  With some more work, the remaining details
  of the $\XXX$-action can then be reconstructed.}
\begin{equation}\label{schem-Verma}
  \xymatrix@=12pt{
    &&\stackrel{\XX^{a}_{r,r'}}{\clubsuit}
    \ar@/^/[dl]_{\fm} \ar@/_/[dr]^{\fp}&\\
    &\stackrel{\XX^{-a}_{r,p_--r'}}
    {\diamondsuit}\ar@/^/[dr]_{\fp}&
    &\stackrel{\XX^{-a}_{p_+-r,r'}}{\heartsuit}
    \ar@/_/[dl]^{\fm}\\
    &&\stackrel{\XX^{a}_{p_+-r,p_--r'}}{\spadesuit}&
  }
\end{equation}
(a ``four-suit'' module).  We note that $\dim\VV^{a}_{r,r'}=p_+p_-$.

\subsubsection{Projective modules}\label{subsec:proj-mod} The $2p_+
p_-$ projective modules $\PP^{\pm}_{r,r'}$, with $1\leq r\leq p_+$ and
$1\leq r'\leq p_-$, are constructed in the standard way, as a
projective cover of each irreducible module.

\begin{prop}\label{prop:proj-struc}
  The subquotient structure of the projective $\XXX$-modules is as
  follows.
  
  \begin{enumerate}
  \item $\PP^{\pm}_{p_+,p_-}\simeq\XX_{p_+,p_-}^{\pm}$ are irreducible
    modules.
    
  \item For $1\leq r\leq p_+-1$, the projective module
    $\PP^{\pm}_{r,p_-}$ admits a filtration
    \begin{equation*}
      \NN_0^\pm\subset\NN_1^\pm\subset\PP^{\pm}_{r,p_-},
    \end{equation*}
    where $\NN_0^\pm\simeq\XX^{\pm}_{r,p_-}$,
    $\NN_1^\pm/\NN_0^\pm\simeq\XX^{\mp}_{p_+-r,p_-}
    \oplus\XX^{\mp}_{p_+-r,p_-}$, and
    $\PP^{\pm}_{r,p_-}/\NN_1^\pm\simeq\XX^{\pm}_{r,p_-}$.
    
  \item For $1\leq r'\leq p_--1$, the projective module
    $\PP^{\pm}_{p_+,r'}$ admits a filtration
    \begin{equation*}
      \NN_0^\pm\subset\NN_1^\pm\subset\PP^{\pm}_{p_+,r'},
    \end{equation*}
    where $\NN_0^\pm\simeq\XX^{\pm}_{p_+,r'}$,
    $\NN_1^\pm/\NN_0^\pm\simeq\XX^{\mp}_{p_+,p_--r'}
    \oplus\XX^{\mp}_{p_+,p_--r'}$, and
    $\PP^{\pm}_{p_+,r'}/\NN_1^\pm\simeq\XX^{\pm}_{p_+,r'}$.
    
  \item For $1\leq r\leq p_+-1$ and $1\leq r'\leq p_--1$, the
    projective module $\PP^{\pm}_{r,r'}$ admits a filtration
    \begin{equation*}
      \NN_0^\pm\subset\NN_1^\pm\subset\NN_2^\pm
      \subset\NN_3^\pm\subset\PP^{\pm}_{r,r'},
    \end{equation*}
    where $\NN_0^\pm \simeq\XX^{\pm}_{r,r'}$, \ $\NN_1^\pm/\NN_0^\pm
    \simeq 2\XX^{\mp}_{r,p_--r'}\oplus 2\XX^{\mp}_{p_+-r,r'}$, \ 
    $\NN_2^\pm/\NN_1^\pm \simeq 4\XX^{\pm}_{p_+-r,p_--r'}\oplus
    2\XX^{\pm}_{r,r'}$, \ $\NN_3^\pm/\NN_2^\pm \simeq
    2\XX^{\mp}_{r,p_--r'}\oplus 2\XX^{\mp}_{p_+-r,r'}$, and
    $\PP^{\pm}_{r,r'}/\NN_3^\pm \simeq\XX^{\pm}_{r,r'}$.
  \end{enumerate}
\end{prop}

\begin{figure}[tbp]
  \centering
  \includegraphics[bb=1.7in 6.8in 7.6in 9.6in, clip]{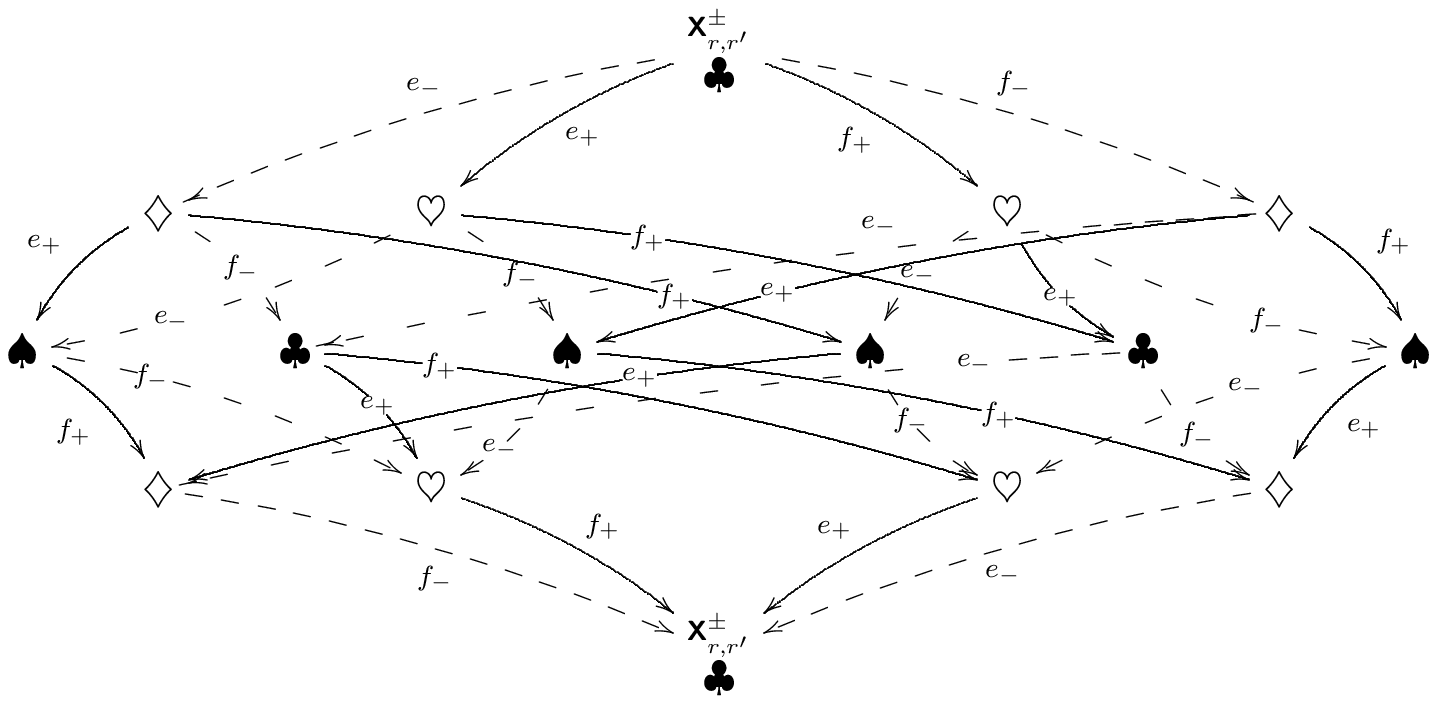}
    \caption[Subquotients of the projective
    module]{\captionfont{Subquotients of the projective module.}
      \small We use the notation
    \begin{alignat*}{2}
      \clubsuit&=\XX^{\pm}_{r,r'},\qquad
      &\diamondsuit&=\XX^{\mp}_{r,p_--r'},\\
      \heartsuit&=\XX^{\mp}_{p_+-r,r'},\qquad
      &\spadesuit&=\XX^{\pm}_{p_+-r,p_--r'}.
    \end{alignat*}
    Continuous lines show the action of the $\ep$ and $\fp$ generators
    and dashed lines show the action of the $\emi$ and $\fm$
    generators.}
  \label{app:proj-pic-plus}
\end{figure}

The structure of subquotients and the left $\XXX$-action on
$\PP^{\pm}_{r,r'}$ can be visualized with the aid of the diagram
in~Fig.~\ref{app:proj-pic-plus}.
Their definition in terms of bases, whence~\bref{prop:proj-struc}
follows, is given in~\bref{app:proj-mod}.

\subsubsection{The $\smash{\XXX}$ Grothendieck ring}\label{sec:Grring}
Closely following the strategy in~\cite{[FGST]}, we obtain
Theorem~\bref{thm:Gr-ring} describing the Grothendieck ring structure
of~$\XXX$.  We do not repeat the standard steps leading to this
statement, the derivation being totally similar to the one
in~\cite{[FGST]}; it relies on a property of tensor products with
Verma modules and simple explicit calculations for two-dimensional
representations.  The $\XXX$ Grothendieck ring is generated by
$\repX^+_{1,2}$ and $\repX^+_{2,1}$.  In~\bref{sec:Grring-U}, it is
further characterized as a polynomial quotient~ring.

\subsection{The space of $q$-characters for
  $\smash{\XXX}$}\label{sec:all-q-chars} We now use irreducible and
projective modules to introduce a basis in the space $\Ch=\Ch(\XXX)$
of $q$-characters of $\XXX$.  In subsequent sections, this basis is
mapped into the center, to produce two distinguished bases
there.\enlargethispage{\baselineskip}

\subsubsection{} For a Hopf algebra $A$, the space $\Ch=\Ch(A)$
of $q$-characters is defined as
\begin{multline}\label{Ch-def}
  \Ch(A)=\bigl\{\beta\in A^* \bigm| \ad^*_x(\beta)
  =\epsilon(x)\beta\quad \forall x\in A\bigr\}\\
  = \bigl\{\beta\in A^* \bigm| \beta(xy)=\beta\bigl(S^2(y)x\bigr)
  \quad \forall x,y\in A\bigr\},
\end{multline}
where the coadjoint action $\ad^*_a:A^*\to A^*$ is
$\ad^*_a(\beta)=\beta\bigl(\sum_{(a)} S(a')?a''\bigr)$, $a\in A$,
$\beta\in A^*$.

In what follows, we need the so-called \textit{balancing element}
$\balance\in A$ that satisfies~\cite{[Drinfeld]}
\begin{equation*}
  \Delta(\balance)=\balance\tensor\balance,\quad
  S^2(x)=\balance x\balance^{-1}
\end{equation*}
for all $x\in A$.  For $A=\XXX$, $\balance=K^{p_+ - p_-}$, as we
calculate in~\eqref{balancing} below.

\subsubsection{Irreducible representation traces}
\label{subsec:func-gamma}
The space of $q$-characters contains a homomorphic image of the
Grothendieck ring under the $q$-trace: for any
$A$-module~$\XXgen$,
\begin{equation}\label{qCh}
  \qTr_{\XXgen} \equiv \tr_{\XXgen}(\balance^{-1}?)\in\Ch(A),
\end{equation}
where $\balance$ is the balancing element~\cite{[Drinfeld]}.  For
$A=\XXX$, we thus have a $2p_+ p_-$-dimen\-sion\-al subspace in $\Ch$
spanned by $q$-traces over irreducible modules, i.e., by
\begin{align*}
  &\gamma^{\pm}(r,r'):x\mapsto\tr_{\XX^{\pm}_{r,r'}}(\balance^{-1}x),
  \quad
  1\leq r\leq p_+,\quad
  1\leq r'\leq p_-,\\
  &\gamma^{\pm}(r,r')\in\Ch.
\end{align*}

\subsubsection{Pseudotraces}\label{sec:pseudo}
The space of $q$-characters $\Ch(\XXX)$ is not spanned by $q$-traces
over irreducible modules; it also contains ``pseudotraces'' associated
with projective modules (in the Kazhdan--Lusztig context, they can be
considered quantum group counterparts of pseudotraces~\cite{[My]}).

The strategy for constructing the pseudotraces is as follows.  For any
indecomposable module $\mathbb{P}$ and a map
$\sigma:\mathbb{P}\to\mathbb{P}$, the functional
\begin{equation}
 \gamma{}:{}
  x\mapsto \Tr_{\mathbb{P}}
  (\balance^{-1} x \sigma)
\end{equation}
is a $q$-character if and only if (cf.~\eqref{Ch-def})
\begin{equation}\label{q-char-cond}
    0=\gamma(x y)- \gamma(S^2(y)x) \equiv
    \Tr_{\mathbb{P}}
    (\balance^{-1}x[y,\sigma]).
\end{equation}
We now find maps $\sigma$ satisfying~\eqref{q-char-cond}.

First, to choose $\mathbb{P}$, we note that the space of
$q$-characters $\Ch(\XXX)$ is the center of the dual algebra $\XXX^*$.
To describe $\Ch(\XXX)$, we can therefore use the isomorphism between
the center of $\XXX^*$ and the algebra of \textit{bimodule}
endomorphisms of the regular representation of $\XXX^*$.  As a
$\XXX^*$-bimodule, the regular representation is contragredient to the
representation described in~\bref{sec-RegDecompose}.  This means that
$q$-characters are endomorphisms of the modules that are
contragredient to projective $\XXX$-modules.  We therefore take
$\mathbb{P}$ to be the projective $\XXX$-module in one full
subcategory, that is, the direct sum
\begin{equation}\label{PP}
  \mathbb{P}_{r,r'}
  =\PP^{+}_{r,r'}\oplus\PP^{-}_{p_+-r,r'} 
  \oplus\PP^{-}_{r,p_--r'}\oplus\PP^{+}_{p_+-r,p_--r'}
\end{equation}
(for $(r,r')\in\setR$) of the four indecomposable projective
$\XXX$-modules $\PP^{\pm}_{r,r'}$ described in~\bref{app:proj-mod}.
We then consider the maps
\begin{equation*}
  \sigma:\mathbb{P}_{r,r'} \to\mathbb{P}_{r,r'}
\end{equation*}
defined by its action on the corresponding basis vectors (see
\eqref{app:projective-basis}): in each of the four modules
in~\eqref{PP}, $\sigma$ acts by zero on all elements
in~\eqref{app:projective-basis} except
\begin{equation*}
  \sigma:\prDown{b}_{n,n'}
  \mapsto\alpha^{\uparrow}\prUp{b}_{n,n'}
  + \alpha^{\downarrow}\prDown{b}_{n,n'} +
  \beta^{\uparrow}\prUp{t}_{n,n'} +
  \beta^{\downarrow}\prDown{t}_{n,n'},
\end{equation*}
with the coefficient $\alpha^{\uparrow}$, $\alpha^{\downarrow}$,
$\beta^{\uparrow}$, and $\beta^{\downarrow}$ chosen arbitrarily in
each module.  To distinguish the basis elements in each of the four
modules in~\eqref{PP}, we introduce the second superscript and write
\begin{equation}\label{Sigma-def}
  \sigma:\nprDown{b}{\bullet}_{n,n'}
  \mapsto\alphaUp{\bullet}\nprUp{b}{\bullet}_{n,n'}
  + \alphaDown{\bullet}\nprDown{b}{\bullet}_{n,n'} +
  \betaUp{\bullet}\nprUp{t}{\bullet}_{n,n'} +
  \betaDown{\bullet}\nprDown{t}{\bullet}_{n,n'},
\end{equation}
where $\sbullet \in \{\rightarrows\}$ and the symbols where $\sbullet$
``takes the value'' \,$\suparrow$\, correspond to the basis elements
in $\PP^{+}_{r,r'}$, the elements where $\sbullet=\sleftarrow$
correspond to the basis elements in $\PP^{-}_{r,p_--r'}$, those where
$\sbullet=\srightarrow$ to the basis elements in $\PP^{-}_{p_+-r,r'}$,
and those where $\sbullet=\sdownarrow$ to the basis elements in
$\PP^{+}_{p_+-r,p_--r'}$.  This $\sigma$ map thus depends on ($(r,r')$
and) sixteen coefficients but, abusing the notation, we do not write
these parameters explicitly.

For any such $\sigma$, we now define a functional on $\XXX$ as
\begin{equation*}
  \gamma(r,r'):
  x\mapsto \Tr_{\mathbb{P}_{r,r'}}
  (\balance^{-1} x \sigma).
\end{equation*}
A lengthy but straightforward calculation shows the following
proposition.
\begin{prop}\label{prop:gammaUpUp}
  For $(r,r')\in\setR$,
  \begin{equation}\label{gamma-general}
    \gamma(r,r')\in\Ch
  \end{equation}
  if and only if
  \begin{alignat*}{2}
    \alphaUp{\suparrow}  &= \alphaUp{\srightarrow}, &
    \alphaUp{\sdownarrow}  &=  \alphaUp{\sleftarrow},\\
    \betaDown{\suparrow}  &= \betaDown{\sleftarrow},& 
    \betaDown{\sdownarrow}  &= \betaDown{\srightarrow},\\
    \betaUp{\suparrow} &= \betaUp{\sleftarrow} &={}
    \betaUp{\srightarrow} &= \betaUp{\sdownarrow}.
  \end{alignat*}
\end{prop}
The proof is outlined in~\bref{app:proof-prop:gammaUpUp}.

\subsubsection{The $\gamma$ basis}\label{subsec:func-gammaUpUp}
In what follows, we assume the relations in~\bref{prop:gammaUpUp} to
be satisfied.  

If $\alphaUp{\suparrow} =\alphaUp{\sdownarrow}= \betaDown{\suparrow} =
\betaDown{\sdownarrow}=\betaUp{\suparrow}=0$, then
\begin{multline*}
  \gamma(r,r')=\alphaDown{\suparrow}\gamma^{+}(r,r') +
  \alphaDown{\sleftarrow}\gamma^{-}(r,p_-\!-\!r')\\*
  {}+ \alphaDown{\srightarrow}\gamma^{-}(p_+\!-\!r,r') +
  \alphaDown{\sdownarrow}\gamma^{+}(p_+\!-\!r,p_-\!-\!r')
\end{multline*}
is a linear combination of the $q$-traces over irredicible modules
introduced in~\bref{subsec:func-gamma}. To construct the other
linearly independent $q$-characters, we take one of the coefficients
$\alphaUp{\sbullet}$, $\betaDown{\sbullet}$, $\betaUp{\sbullet}$ to be
nonzero and thus define
\begin{multline}\label{last-azat-label}
  \gamma^{\NESW}(r,r')
  =\gamma(r,r')\quad\mbox{for}\quad
    \alphaUp{\sdownarrow} =
    \alphaDown{\bullet} = \betaDown{\bullet} =\betaUp{\bullet}=0,\\
    \alphaUp{\suparrow} = \ffrac{[r]_+}{\qp^{p_-}\!-\!\qp^{-p_-}},
    \qquad
    (r,r')\in\setILeft,
\end{multline}
\vspace{-\baselineskip}
\begin{multline}\label{gammaNESW-def}
  \gamma^{\NESW}(p_+\!-\!r,p_-\!-\!r')=\gamma(r,r')
  \quad\mbox{for}\quad
  \alphaUp{\suparrow}
    = \alphaDown{\bullet} = \betaDown{\bullet} =\betaUp{\bullet}=0,\\
    \alphaUp{\sdownarrow} =
    \ffrac{[p_+\!-\!r]_+}{\qp^{p_-}\!-\!\qp^{-p_-}},
    \qquad (r,r')\in\setILeft,
\end{multline}
\vspace{-\baselineskip}
\begin{multline}
  \gamma^{\NWSE}(r,r')=\gamma(r,r')
  \quad\mbox{for}\quad
  \alphaDown{\bullet}
  =\alphaUp{\bullet} = \betaDown{\sdownarrow} = \betaUp{\bullet}=0,\\
   \betaDown{\suparrow} = \ffrac{[r']_-}{\qm^{p_+}\!-\!\qm^{-p_+}},
   \qquad
  (r,r')\in\setIRight,
\end{multline}
\vspace{-\baselineskip}
\begin{multline}
  \gamma^{\NWSE}(p_+\!-\!r,p_-\!-\!r')=\gamma(r,r')
  \quad\mbox{for}\quad
  \alphaDown{\bullet}
  =\alphaUp{\bullet} = \betaDown{\suparrow} = \betaUp{\bullet}=0\\
  \betaDown{\sdownarrow} =
  \ffrac{[p_-\!-\!r']_-}{\qm^{p_+}\!-\!\qm^{-p_+}},
  \qquad
  (r,r')\in\setIRight,
\end{multline}
\vspace{-\baselineskip}
\begin{multline}
  \gamma^{\UpUp}(r,r')=\gamma(r,r')
  \quad\mbox{for}\quad
  \alphaDown{\bullet}
  =\alphaUp{\bullet} =\betaDown{\bullet}=0,\\
  \betaUp{\suparrow}
  = \alphaUp{\suparrow} \betaDown{\suparrow},\qquad
  (r,r')\in\setR.
\end{multline}
For convenience, we introduce the maps $\sigma^{\NESW}(r,r')$ (with
$(r,r')\in\setILeft$), $\sigma^{\NWSE}(r,r')$ (with
$(r,r')\in\setIRight$), and $\sigma^{\UpUp}(r,r')$ (with
$(r,r')\in\setR$) obtained by substituting the respective sets of
coefficients in~\eqref{Sigma-def}.

We thus have
 $ 2p_+ p_- + (p_+\!-\!1)p_- + p_+(p_-\!-\!1)
  +\half(p_+\!-\!1)(p_-\!-\!1)
  =\half(3p_+\!-\!1)(3p_-\!-\!1)$
$q$-characters in the space $\Ch$:
\begin{equation}\label{eq:gamma-basis}
  \mbox{}\kern-13pt\begin{gathered}
    \gamma^{+}(p_+,p_-),\qquad
    \underset{\strut\mbox{\footnotesize$1\leq r\leq p_+\!-\!1$}}{
      \gamma^{\NESW}(r,p_-),}\qquad
    \underset{\strut\mbox{\footnotesize$(r,r')\in\setR$}}{
      \gamma^{\UpUp}(r,r'),}\qquad
    \underset{\strut\mbox{\footnotesize$1\leq r'\leq p_-\!-\!1$}}{
      \gamma^{\NWSE}(p_+,r'),}\qquad
    \gamma^{-}(p_+,p_-),\\[4pt]
    \underset{\strut\mbox{\footnotesize$1\leq r\leq p_+\!-\!1$}}{
      \gamma^{+}(r,p_-),~~\gamma^{-}(p_+\!-\!r,p_-),}\quad
    \underset{\strut\mbox{\footnotesize$
        \begin{array}[t]{c}
          1\leq r\leq p_+\!-\!1,\\
          1\leq r'\leq p_-\!-\!1
        \end{array}$}}{
      \gamma^{\NESW}(r,r'),~~\gamma^{\NWSE}(r,r'),}\quad
    \underset{\strut\mbox{\footnotesize$1\leq r'\leq p_-\!-\!1$}}{
      \gamma^{+}(p_+, r'),~\gamma^{-}(p_+, p_-\!-\!r'),}\\[4pt]
    \underset{\strut\mbox{\footnotesize$(r,r')\in\setR$}}{
      \gamma^{+}(r,r'),~\gamma^{-}(p_+\!-\!r,r'),~
      \gamma^{-}(r,p_-\!-\!r'),~\gamma^{+}(p_+\!-\!r,p_-\!-\!r')}.
  \end{gathered}\kern-7pt
\end{equation}
There are $\half(3p_+\!-\!1)(3p_-\!-\!1)$ linearly independent
elements in $\Ch$ thus obtained.  In what follows, we also establish
that $\dim\cZ=\half(3p_+\!-\!1)(3p_-\!-\!1)$, with the consequence
that the $\gamma$ listed above are a basis in~$\Ch$.  We call it
\textit{the $\gamma$-basis} in what follows.

\subsection{The $\smash{\XXX}$ center}\label{sec:center}
We now describe the structure of the center of~$\XXX$.

\begin{thm}\label{prop-center}
  The center $\cZ$ of $\XXX$ is
  $\half(3p_+\,{-}\,1)(3p_-\,{-}\,1)$-dimensional and decomposes into
  a direct sum of associative algebras as
  \begin{equation}\label{center-decomposition}
    \cZ=\algI_{p_+,p_-}^{(1)}\oplus
    \algI_{0,p_-}^{(1)}\oplus
    \bigoplus_{r=1}^{p_+-1}\algB_{r,p_-}^{(3)}
    \oplus
    \bigoplus_{r'=1}^{p_--1}\algB_{p_+,r'}^{(3)}
    \oplus
    \bigoplus_{r,r'\in\setR}\algA_{r,r'}^{(9)},
  \end{equation}
  where the dimension of each algebra is shown as a superscript.
  Furthermore,
  \begin{itemize}
  \item $\algI_{p_+,p_-}^{(1)}$ is spanned by a single idempotent
    $\Idem(p_+,p_-)$,
    
  \item $\algI_{0,p_-}^{(1)}$ is spanned by a single idempotent
    $\Idem(0,p_-)$,
    
  \item each $\algB_{r,p_-}^{(3)}$ is spanned by the primitive
    idempotent $\Idem(r,p_-)$ and the radical elements $\vUp(r,p_-)$
    and $\vRight(r,p_-)$, on all of which $\Idem(r,p_-)$ acts as
    identity and the other products vanish,
    
  \item each $\algB_{p_+,r'}^{(3)}$ is spanned by the primitive
    idempotent $\Idem(p_+,r')$ and the radical elements $\vUp(p_+,r')$
    and $\vLeft(p_+,r')$, on all of which $\Idem(p_+,r')$ acts as
    identity and the other products vanish,
    
  \item each $\algA_{r,r'}^{(9)}$ is spanned by the primitive
    idempotent $\Idem(r,r')$ \textup{(}acting as identity on
    $\algA_{r,r'}^{(9)}$\textup{)} and the radical elements
    $\vNE(r,r')$, $\vSW\!(r,r')$, $\vNW(r,r')$, $\vSE(r,r')$,
    $\wUp\!(r,r')$, $\wRight(r,r')$, $\wDown(r,r')$, $\wLeft(r,r')$,
    which have the nonzero products\footnote{The arrow notation is
      partly justified by the ``momentum conservation law'' satisfied
      in~\eqref{rad-times}.}
   \begin{equation}\label{rad-times}
     \begin{alignedat}{2}
       \vNE(r,r')\vNW(r,r')&=\wUp(r,r'),&
       \quad
       \vNE(r,r')\vSE(r,r')&=\wRight(r,r'),
       \\
       \vSW(r,r')\vNW(r,r')&=\wLeft(r,r'),&
       \quad
       \vSW(r,r')\vSE(r,r')&=\wDown(r,r').
     \end{alignedat}
   \end{equation}
   Pairwise products between elements in the radical other than those
   in~\eqref{rad-times} are identically zero.
  \end{itemize}
\end{thm}
The theorem is proved in Appendix~\bref{app:the-center}.

The $\half(3p_+\,{-}\,1)(3p_-\,{-}\,1)$ elements $\Idem$,
$\boldsymbol{v}$, and $\boldsymbol{w}$ are called the
\textit{canonical central elements}.  A number of facts in what
follows are established by decomposing various central elements with
respect to the basis of the canonical central element.  But because
our next aim is to study the Radford map, it proves convenient to use
a different basis in the center, given by the Radford-map images of
the above $\gamma\in\Ch$ and actually obtained by certain mixing of
the $\Idem$, $\boldsymbol{v}$, and $\boldsymbol{w}$ within each
respective algebra in~\eqref{center-decomposition}.

\section{Radford map of  $\smash{\Ch(\XXX)}$}\label{sec:radford} The
$\SLiiZ$-action on the center to be considered in Sec.~\ref{sec:SL2Z}
involves the Radford and Drinfeld maps and a ribbon structure; the
Radford map $\Radford$ is studied in this section.  We first find the
necessary ingredients for $\Radford$ and $\Radford^{-1}$ (the
cointegral and integral) and then evaluate the Radford images of the
basis in $\Ch$ constructed in the previous section.

\subsection{Radford map} We recall that for a Hopf algebra $A$, a
\textit{right integral}~$\rint$ is a linear functional on $A$
satisfying
\begin{equation*}
  (\rint\tensor\id)\Delta(x)=\rint(x)\one
\end{equation*}
for all $x\,{\in}\, A$.  If such a functional exists, it is unique up
to multiplication by a nonzero constant.  The left--right\footnote{For
  a unimodular $A^*$, which is the case with~$A=\XXX$.}
\textit{cointegral}~$\coint$ is an element in $A$ such that
\begin{equation*}
  x\coint=\coint x =\epsilon(x)\coint,\quad\forall x\in A.
\end{equation*}
Whenever it exists, this element is unique up to multiplication by a
nonzero constant.  We also note that the cointegral gives an embedding
of the trivial representation of $A$ in the regular bimodule~$A$.  We
use the standard normalization $\rint(\coint)=1$.

Next, a \textit{comodulus}~$\comodul$ is an element in $A$ such that
\begin{equation*}
  (\id\tensor\rint)\Delta(x)
  =\rint(x)\comodul.
\end{equation*}
Whenever a square root of the comodulus exists, we
have~\cite{[Drinfeld]}
\begin{equation}\label{bal-comod}
  \balance^2=\comodul,
\end{equation}
where $\balance\in A$ is the balancing element.

The Radford map $\Radford:A^*\to A$ and its inverse
$\Radford^{-1}:A\to A^*$ are given by
\begin{equation}\label{radford-def}
  \Radford(\beta)
  =\sum_{(\coint)}\beta(\coint')\coint'',
  \qquad
  \Radford^{-1}(x)=\rint(S(x)?).
\end{equation}
Both $\Radford$ and $\Radford^{-1}$ intertwine the left actions of $A$
on $A$ and $A^*$, with the left $A$-module structures given by
$a\acts\beta=\beta(S(a)?)$ on $A^*$ and by the regular action on $A$,
and similarly for the right actions~\cite{[Swe],[Rad]}.  We use the
hat for notational consistency in what follows (and swap $\Radford$
and $\Radford^{-1}$ compared with the commonly used notation).

In particular, $\Radford$ gives an isomorphism of linear spaces
$\Radford:\Ch\to\cZ$. 

We now find the integral and cointegral and then evaluate $\Radford$
on the $q$-characters constructed in~\bref{sec:all-q-chars}.

\subsubsection{The integral, cointegral, and comodulus for
  $\XXX$}\label{sec:int} A simple calculation in $\XXX$ shows that the
right integral is given by
\begin{equation}\label{rint}
  \rint(\emi^{m'}\fp^{n} K^j \ep^{m}\fm^{n'})
  = \ffrac{\qp^{2p_-}\qm^{2p_+}}{\zeta}\,
  \delta_{m',p_- - 1}\delta_{n,p_+ - 1}\delta_{m,p_+ - 1}
  \delta_{n',p_- - 1} \delta_{j, p_+ - p_-}.
\end{equation}
Next, as is easy to verify, the left--right cointegral for~$\XXX$ is
given by the product of the cointegrals for the two
$\overline{\mathscr{U}}_{q} s\ell(2)$ in~\eqref{as-quotient},
\begin{equation}\label{coint}
  \coint
  =\zeta\,\fp^{p_+ - 1}\ep^{p_+ - 1}\fm^{p_- - 1}\emi^{p_- - 1}
  \!\sum_{n=0}^{2 p_+ p_- - 1}K^n
\end{equation}
(which is at the same time the product of the cointegrals for the two
$\overline{\mathscr{U}}_{q} s\ell(2)$ in~\eqref{as-quotient}, up to
normalization).  \textit{Factorizable ribbon} quantum groups offer a
``canonical'' normalization of the Radford map (i.e., of the integral
and cointegral) up a power of $i=\sqrt{-1}$, from the condition
$\modS^4|_{\cZ}=\id$ for $\modS$ in~\eqref{3-diagram},
and we therefore set
\begin{equation}\label{Rad-normalization}
  \zeta([p_+ - 1]_+!\,[p_- - 1]_-!)^2=\Aconstplain.
\end{equation}

We further calculate
\begin{equation*}
  (\id\tensor\rint)
  \Delta(\emi^{p_- - 1}\fp^{p_+ - 1} k^{2p_+ - 2 p_-}
  \ep^{p_+ - 1}\fm^{p_- - 1})
  =\ffrac{\qp^{2p_-}\qm^{2p_+}}{\zeta}\,k^{4 p_+ - 4 p_-}, 
\end{equation*}
which allows choosing the comodulus as $\comodul=K^{2 p_+ - 2 p_-}$.
In accordance with~\eqref{bal-comod}, we then have the balancing
element
\begin{equation}\label{balancing}
  \balance=K^{p_+ - p_-}.
\end{equation}

\subsection{Radford map of the $\gamma$ basis in $\Ch(\XXX)$}
Explicitly, the Radford map $\Radford:\Ch\to\cZ$ of the irreducible
representation characters is given by
\begin{align*}
  \XX^{\pm}_{r,r'}\mapsto{}&
  \Radford^{\pm}(r,r')\equiv\Radford(\qTr_{\XX^{\pm}_{r,r'}})
  =\sum_{(\coint)} \Tr_{\XX^{\pm}_{r,r'}}(K^{p_- - p_+}\coint')\,
  \coint'',\quad
  \begin{array}[t]{l}
    1\leq r\leq p_+,\\ 1\leq r'\leq p_-.
  \end{array}
\end{align*}
Clearly, $\Radford^{\,+}(1,1)=\coint$, in accordance with the fact
that $\coint$ furnishes an embedding of the trivial representation
$\XX^+_{1,1}$ into~$\XXX$.  Because the Radford map $\Grring_{2 p_+
  p_-}\to\ZZ$ is a morphism of representations, it follows that the
linear span of the $\Radford^{\pm}(r,r')$ is the annihilator of the
radical in the center (see formulas below
in~\bref{prop:radford-decompose}).

The elements $\gamma^{\NESW}(r,r')$ are mapped under $\Radford$ as
\begin{equation}\label{radfordNESW-prepare}
  \Radford(\Ngamma^{\NESW}(r,r'))=\NRadfordNESW(r,r')=
  \sum_{(\coint)}
  \Tr_
  {\mathbb{P}_{r,r'}}
  (K^{p_- - p_+}\coint'\NsigmaNESW(r,r'))\,\coint''
\end{equation}
for all $1\leq r\leq p_+\,{-}\,1$ and $1\leq r'\leq p_-$.

The elements $\gamma^{\NWSE}(r,r')$ are mapped under $\Radford$ as
\begin{equation}
  \label{radfordNWSE-prepare}
  \Radford(\Ngamma^{\NWSE}(r,r'))=\NRadfordNWSE(r,r')=
  \sum_{(\coint)}
  \Tr_
  {\mathbb{P}_{r,r'}}
  (K^{p_- - p_+}\coint'\NsigmaNWSE(r,r'))\,\coint'',
\end{equation}
for all $1\leq r\leq p_+$ and $1\leq r'\leq p_-\,{-}\,1$.

Finally, the elements $\gamma^{\UpUp}(r,r')$ with $(r,r')\in\setR$ are
mapped under $\Radford$ as
\begin{equation}
  \Radford(\Ngamma^{\UpUp}(r,r'))=\NRadfordUpUp(r,r')= \sum_{(\coint)}
  \Tr_{\mathbb{P}_{r,r'}}\!( K^{p_- -
  p_+}\coint'\NsigmaUpUp(r,r'))\coint''.\kern-4pt
\end{equation}

\begin{prop}\label{prop:radford-decompose}
  The Radford-map images
  \begin{equation*}
    \Radford^\bullet(r,r')=
    \Radford(\gamma^\bullet(r,r')),
    \quad \bullet=\UpUp,\NESW,\NWSE,+,-,
  \end{equation*}
  of the $\gamma$-basis in~\eqref{eq:gamma-basis} are as follows.
  \begin{enumerate}
  \item The Radford images of the irreducible representation traces
    $\gamma^{\pm}(r,r')$, $1\leq r\leq p_+$, $1\leq r'\leq p_-$, are
    given by
    \begin{align*}
      &\begin{aligned} \Radford^{+}(r,r')&= \Aoverviiipp\,
        \wUp(r, r'),\\
        \Radford^{-}(r, p_-\!-\!r')&= \Aoverviiipp\,
        \wLeft(r, r'),\\
        \Radford^{-}(p_+\!-\!r,r')&= \Aoverviiipp\,
        \wRight(r, r'),\\
        \Radford^{+}(p_+\!-\!r,p_-\!-\!r')&= \Aoverviiipp\, \wDown(r,
        r'),
      \end{aligned}\quad (r,r')\in\setR,\\
      &\begin{aligned} \Radford^{+}(p_+,r')&= (-1)^{p_- +
          r'}\ffrac{p_+}{2 p_-}\, \Aconstplain \,
        \vUp(p_+, r'),\\
        \Radford^{-}(p_+, p_-\!-\!r')&= (-1)^{p_- + r'} \ffrac{p_+}{2
          p_-}\, \Aconstplain\,
        \vLeft(p_+, r'),\\
      \end{aligned}\quad 1\leq r'\leq p_-\!-\!1,
      \\
      &\begin{aligned} \Radford^{+}(r,p_-)&= (-1)^{p_+ + r}
        \ffrac{p_-}{2 p_+}\, \Aconstplain \,
        \vUp(r, p_-),\\
        \Radford^{-}(p_+\!-\!r, p_-)&= (-1)^{p_+ + r}\ffrac{p_-}{2
          p_+}\, \Aconstplain \,
        \vRight(r, p_-),\\
      \end{aligned}\quad 1\leq r\leq p_+\!-\!1,\\
      &\begin{aligned} \Radford^{+}(p_+,p_-)&=
        \iippA \Idem(p_+, p_-),\\
        \Radford^{-}(p_+,p_-)&= (-1)^{p_+ + p_-} \iippA \Idem(0, p_-);
      \end{aligned}
    \end{align*}
    
  \item the Radford images of $\gamma^{\NESW}(r,r')$ with    
    $1\leq r \leq p_+\!-\!1$, $1 \leq r'\leq p_-$ and of
    $\gamma^{\NWSE}(r,r')$ with $1\leq r\leq p_+$, $1\leq r'\leq
    p_-\!-\!1$ are given by
    \begin{multline}\label{radfordNESW-decompose}
      \NRadfordNESW(r,r') = (-1)^{r'}\ffrac{\qp^{p_- r}\!+ \qp^{-p_-
          r}}{ \qp^{p_- r}\!- \qp^{-p_- r}} \,\bigl(\Radford^+(r, r')
      + \Radford^-(p_+\!-\!r, r')\bigr)
      \\*
      {}- \ffrac{1}{ \qp^{p_- r}\!- \qp^{-p_- r}}\,\ffrac{p_+}{2
        p_-}\, \Aconstplain\cdot
      \begin{cases}
        \vNE(r, r'),& (r,r')\in\setR,\\
        \vSW(p_+\!-\!r, p_-\!-\!r'),&(p_+\!-\!r, p_-\!-\!r')\in\setR,
      \end{cases}
    \end{multline}
    \vspace{-.4\baselineskip}
    \begin{equation}\label{radford-bdry-decomp1}
      \NRadfordNESW(r,p_-) = (-1)^{r + p_+ + 1}
      \ffrac{\iippA}{\qp^{p_- r}\!- \qp^{-p_- r}}\,
      \Idem(r, p_-)
      + (-1)^{p_-}\ffrac{\qp^{p_- r}\!+ \qp^{-p_- r}}{
        \qp^{p_- r}\!- \qp^{-p_- r}}\,
      \RadfordProj(r,p_-),
    \end{equation}
    \vspace{-.6\baselineskip}
    \begin{multline}\label{radfordNWSE-decompose}
      \NRadfordNWSE(r,r')= (-1)^{r}\ffrac{\qm^{p_+ r'}\!+ \qm^{-p_+
          r'}}{ \qm^{p_+ r'}\!-\!\qm^{-p_+ r'}}
      \,\bigl(\Radford^+(r, r') + \Radford^-(r, p_-\!-\!r')\bigr)\\*
      - \ffrac{1}{ \qm^{p_+ r'}\!-\!\qm^{-p_+ r'}}\, \ffrac{p_-}{2
        p_+}\, \Aconstplain\cdot
      \begin{cases}
        \vNW(r, r'),& (r,r')\in\setR,\\
        \vSE(p_+\!-\!r, p_-\!-\!r'),&(p_+\!-\!r, p_-\!-\!r')\in\setR,
      \end{cases}
    \end{multline}
    \vspace{-\baselineskip}
    \begin{multline}\label{radford-bdry-decomp2}  
      \NRadfordNWSE(p_+, r') =\\
      =(-1)^{r' + p_- + 1}
      \ffrac{\iippA}{\qm^{p_+ r'}\!\!-\!\qm^{-p_+ r'}}\,
      \Idem(p_+, r')
      +(-1)^{p_+}\,\ffrac{\qm^{p_+ r'}\!\!+\!\qm^{-p_+ r'}}{
        \qm^{p_+ r'}\!\!-\!\qm^{-p_+ r'}}\,
      \RadfordProj(p_+,r');
    \end{multline}
  \item and the Radford images of $\gamma^{\UpUp}(r,r')$,
    $(r,r')\in\setR$, are
    \begin{multline}\label{phiUpUp-decompose}
      \NRadfordUpUp(r,r') = \ffrac{\iippA}{ (\qp^{p_- r}\!-\!\qp^{-p_-
          r}) (\qm^{p_+ r'}\!-\!\qm^{-p_+ r'})}\,
      \Idem(r, r')\\
      +(-1)^{r'} \ffrac{\qp^{p_- r}\!+ \qp^{-p_- r}}{ \qp^{p_-
          r}\!-\!\qp^{-p_- r}} (\NRadfordNWSE(r, r') -
      (-1)^{p_+}\NRadfordNWSE(p_+\!-\!r, p_-\!-\!r'))
      \\
      \shoveright{ {}+(-1)^{r} \ffrac{\qm^{p_+ r'}\!+ \qm^{-p_+ r'}}{
          \qm^{p_+ r'}\!-\!\qm^{-p_+ r'}} (\NRadfordNESW(r, r')
        -(-1)^{p_-}\NRadfordNESW(p_+\!-\!r, p_-\!-\!r'))
        \quad}\\
      {}- (-1)^{r + r'}\ffrac{(\qp^{p_- r}\!+ \qp^{-p_- r}) (\qm^{p_+
          r'}\!+ \qm^{-p_+ r'})}{ (\qp^{p_- r}\!-\!\qp^{-p_- r})
        (\qm^{p_+ r'}\!-\!\qm^{-p_+ r'})}\, \RadfordProj(r,r'),
    \end{multline}
    where we use the notation
    \begin{equation*}
      \RadfordProj(r,r')
      = \Radford^{+}(r,r')
      + \Radford^{-}(p_+\!-\!r,r')
      + \Radford^{-}(r,p_-\!-\!r')
      + \Radford^{+}(p_+\!-\!r,p_-\!-\!r').
    \end{equation*}
  \end{enumerate}  
  The $\half(3p_+-1)(3p_--1)$ elements $\Radford^\bullet(r,r')$ are a
  basis in $\cZ$, and therefore the $\gamma^\bullet(r,r')$ are a basis
  in~$\Ch$.
\end{prop}
This is proved in~\bref{app:radford-proof}.

\subsubsection{} We thus have a basis in the center $\ZZ$ of $\XXX$
associated with the Radford map:
\begin{equation}\label{eq:radford-basis}
  \mbox{}\kern-13pt\begin{gathered}
    \Radford^{+}(p_+,p_-),\qquad
    \underset{\strut\mbox{\footnotesize$1\leq r\leq p_+\!-\!1$}}{
      \RadfordNESW(r,p_-),}\qquad
    \underset{\strut\mbox{\footnotesize$(r,r')\in\setR$}}{
      \RadfordUpUp(r,r'),}\qquad
    \underset{\strut\mbox{\footnotesize$1\leq r'\leq p_-\!-\!1$}}{
      \RadfordNWSE(p_+,r'),}\qquad
    \Radford^{-}(p_+,p_-),\\[4pt]
    \underset{\strut\mbox{\footnotesize$1\leq r\leq p_+\!-\!1$}}{
      \Radford^{+}(r,p_-),~~\Radford^{-}(p_+\!-\!r,p_-),}\quad
    \underset{\strut\mbox{\footnotesize$
        \begin{array}[t]{c}
          1\leq r\leq p_+\!-\!1,\\
          1\leq r'\leq p_-\!-\!1
        \end{array}$}}{
      \RadfordNESW(r,r'),~~\RadfordNWSE(r,r'),}\quad
    \underset{\strut\mbox{\footnotesize$1\leq r'\leq p_-\!-\!1$}}{
      \Radford^{+}(p_+, r'),~\Radford^{-}(p_+, p_-\!-\!r'),}\\[4pt]
    \underset{\strut\mbox{\footnotesize$(r,r')\in\setR$}}{
      \Radford^{+}(r,r'),~\Radford^{-}(p_+\!-\!r,r'),~
      \Radford^{-}(r,p_-\!-\!r'),~\Radford^{+}(p_+\!-\!r,p_-\!-\!r')}.
  \end{gathered}\kern-7pt
\end{equation}
The first line comprises $1\,{+}\,(p_+\,{-}\,1)\,{+}\,
\half(p_+\,{-}\,1)(p_-\,{-}\,1)\,{+}\,(p_-\,{-}\,1)\,{+}\,1 =
\half(p_+\,{+}\,1)\cdot(p_-\,{+}\,1)$ elements that are expressed as the
corresponding primitive idempotents plus possible nilpotent elements.
The second line comprises
$2(p_+\,{-}\,1)\,{+}\,2(p_+\,{-}\,1)(p_-\,{-}\,1)+2(p_-\,{-}\,1)$
elements in the radical, and the bottom line
$4\cdot\half(p_+\,{-}\,1)(p_-\,{-}\,1)$ elements in the radical with
zero pairwise products.

The elements of this basis are called the \textit{``Radford'' central
  elements} in what follows. They are a ``mixing'' of the canonical
central elements in~\bref{prop-center} within each subalgebra
in~\eqref{center-decomposition}, and are used in decompositions in
what follows instead of the canonical central elements
in~\bref{prop-center}.

\subsubsection{More Radford-map notation}\label{Radford-notation} 
Certain linear combinations of the Radford central elements are
extensively used in what follows.  For the convenience of the reader,
we summarize the notation for these combinations here.  We define
$\RadfordProj(r,r')$ with $(r,r')\in\setI$ as
\begin{alignat}{2}\label{Proj-notation}
  \RadfordProj(r,r')
  &= \Radford^{+}(r,r')
  + \Radford^{-}(p_+\!-\!r,r')\\*
  &\quad{}+ \Radford^{-}(r,p_-\!-\!r')
  + \Radford^{+}(p_+\!-\!r,p_-\!-\!r'),&\quad &(r,r')\in\setR,\notag\\
  \label{Proj-notation2}
  \RadfordProj(r,p_-)&={}\,\Radford^{+}(r,p_-) + \Radford^{-}(p_+ -
  r,p_-),&\quad &1\leq r\leq p_+\!-\!1,\\
  \label{Proj-notation3}
  \RadfordProj(p_+,r')&={}\,\Radford^{+}(p_+,r') +
  \Radford^{-}(p_+,p_-\!-\!r'),&\quad &1\leq r'\leq p_-\!-\!1,\\
  \RadfordProj(p_+,p_-)&={}\,\Radford^{+}(p_+,p_-),\\
  \RadfordProj(0,p_-)&={}\,\Radford^{-}(p_+,p_-)
  \label{Proj-notation5}
\end{alignat}
(as we see in what follows, these form a basis in the linear span of
the Drinfeld-map image of projective-module characters.  We also set
\begin{multline}\label{vvarphi-def}
  \vvarphi(r,r')=
  (p_+\!-\!r)(p_-\!-\!r')\,\Radford^+(r,  r')
  - r(p_-\!-\!r')\,\Radford^{-}(p_+\!-\!r, r')\\*
  {}- (p_+\!-\!r) r' \Radford^{-}(r, p_-\!-\!r')
  + r r' \Radford^{+}(p_+\!-\!r, p_-\!-\!r'),
  \quad (r,r')\in\setR.
\end{multline}
Next, we define $\hrhoS(r,r')$ with $(r,r')\in\setILeft$ as
\begin{alignat}{2}
  \label{hrhoright}
  \hrhoS(r,r')&=
  (p_+\!-\!r)
  \bigl(\Radford^{+}(r,r')+\Radford^{-}(r,p_-\!-\!r')
  \bigr)&&\\*
  &\quad{}-r\bigl(\Radford^{-}(p_+\!-\!r,r')
  + \Radford^{+}(p_+\!-\!r, p_-\!-\!r')\bigr),
  &\quad &(r,r')\in\setR,\notag\\
  \label{hrhoright-bdr}
  \hrhoS(r,p_-)
  &=(p_+\!-\!r)\,\Radford^{+}(r,p_-)
  - r \Radford^{-}(p_+\!-\!r,p_-),&\quad &1\leq r\leq p_+\!-\!1,\\
  \intertext{and $\hrhoB(r,r')$ with $(r,r')\in\setIRight$ as}
  \label{hrholeft}
  \hrhoB(r,r')&=
  (p_-\!-\!r')\bigl(\Radford^+(r, r') +
  \Radford^-(p_+\!-\!r,r')\bigr)\\*
  &\quad{}- r'\bigl(\Radford^-(r, p_-\!-\!r')
  + \Radford^+(p_+\!-\!r, p_-\!-\!r')\bigr),
  &\quad &(r,r')\in\setR,\notag\\
  \label{hrholeft-bdr}
  \hrhoB(p_+, r')
  &=(p_-\!-\!r')\,\Radford^{+}(p_+,r')
  - r' \Radford^{-}(p_+, p_-\!-\!r'),&\quad &1\leq r'\leq p_-\!-\!1.
\end{alignat}
We also define $\RadfordS(r,r')$ for $(r,r')\in\setILeft$ as
\begin{alignat}{2}
  \label{hrhoRight}
  \NRadfordS(r,r')
  &=
  (-1)^{r'}\NRadfordNESW(r,r')
  -(-1)^{p_- + r'}
  \NRadfordNESW(p_+\!-\!r, p_-\!-\!r'),&
  \ &(r,r')\in\setR,
  \\
  \label{hrhoRight2}
  \NRadfordS(r,p_-)
  &= (-1)^{p_-}\NRadfordNESW(r,p_-),&\ &1\leq r\leq p_+\!-\!1,\\
  \intertext{and $\RadfordB(r,r')$ for $(r,r')\in\setIRight$ as}
  \label{hrhoLeft}
  \NRadfordB(r,r')
  &=
  (-1)^{r}\NRadfordNWSE(r,r')
  -(-1)^{p_+ + r}
  \NRadfordNWSE(p_+\!-\!r, p_-\!-\!r'),&
  \ &(r,r')\in\setR,
  \\
  \label{hrhoLeft2}
  \NRadfordB(p_+,r')
  &=  
  (-1)^{p_+}\NRadfordNWSE(p_+, r'),&\ &1\leq r'\leq p_-\!-\!1.
\end{alignat}

\section{The factorizable and ribbon structures
  for~$\smash{\XXX}$}\label{sec:drinfeld} In this section, we study
the Drinfeld map and a ribbon structure for $\XXX$.  We find the
necessary ingredient for the Drinfeld map (the $M$-matrix)
in~\bref{sec:M}.  As a preparation for the study of the
$\SLiiZ$-action, we then evaluate the Drinfeld images of the $\gamma$
basis in~$\Ch$.  Together with the ``Radford'' basis constructed in
Sec.~\ref{sec:radford}, this gives two bases in the center that play
a crucial role in the study of the modular group action in
Sec.~\ref{sec:SL2Z}.

We consider the Drinfeld map in~\bref{sec:Drinfeld} and find a ribbon
element in~\bref{XXX-ribbon}.  As a short detour
in~\bref{sec:Grring-U}, we characterise the Grothendieck ring in terms
of Chebyshev polynomials.

\subsection{$M$-matrix}\label{sec:M} For a quasitriangular Hopf
algebra $A$ with the universal $R$-matrix $R$, the $M$-matrix is
defined as
\begin{equation*}
  M=R_{21}R_{12}\in A\tensor A.
\end{equation*}\pagebreak[3]
It satisfies the relations $(\Delta\tensor\id)(M)=R_{32}M_{13}R_{23} $
and $M\Delta(x)=\Delta(x) M$ $\forall x\in A$.  If in addition~$M$ can
be represented as\pagebreak[3]
\begin{equation*}
  M=\sum_I \pbwd_I\tensor \pbwdd_I,
\end{equation*}
where~$\pbwd_I$ and~$\pbwdd_I$ are two \textit{bases} in~$A$, the
Hopf algebra~$A$ is called \textit{factorizable}.

For a Hopf algebra $A$ with an $M$-matrix, the Drinfeld map
$\Drinfeld:A^*\to A$ is defined by
\begin{equation}\label{drinfeld-def}
  \Drinfeld(\beta)=(\beta\tensor\id)M.
\end{equation}
In a factorizable Hopf algebra $A$, the restriction of the Drinfeld
map to the space $\Ch$ of $q$-characters gives an isomorphism
$\Ch(A)\xrightarrow{\simeq}\cZ(A)$ of associative
algebras~\cite{[Drinfeld]}.

\subsubsection{The $M$-matrix for $\smash{\XXX}$}\label{sec:the-M}
As any Drinfeld double, the double $D(\mathscr{H})$ whose construction
was outlined in~\bref{sec:from} (see~\cite{[FGST3]} for the detailed
calculation) is endowed with a universal $R$-matrix.  To evaluate it
explicitly, we construct the elements of the basis dual to the PBW
basis in $\mathscr{H}$, Eqs.~\eqref{ejmn}, in terms of the generators
introduced in~\eqref{basis} and~thus~find
\begin{multline}\label{q-double-R}
  R=\ffrac{1}{4p_+p_-}
  \sum_{m=0}^{p_+-1}\sum_{n=0}^{p_--1}
  \sum_{j=0}^{4p_+p_--1}
  \sum_{c=0}^{4p_+p_--1}
  \ffrac{(-1)^{m+n}(\qp^{p_-}\!- \qp^{-p_-})^m
    (\qm^{p_+}\!- \qm^{-p_+})^n}{
    [m]_+![n]_-!}\\[-3pt]
  {}\times
  \qp^{-mj-p_-\frac{m(m-1)}{2}} \qm^{nj-p_+\frac{n(n-1)}{2}}\q^{-cj}
  \,\kk^j\ep^m\fm^n\tensor\fp^m\emi^n\dkk^c.
\end{multline}

The $M$-matrix for $D(\mathscr{H})$ is $M=R_{21}R_{12}$ with this $R$.
By simple calculation, it then follows that the $M$-matrix $\bar M$
for $\bar D(\mathscr{H})$ (see~\eqref{D-bar}) is given by
\begin{multline}\label{bar-M}
  \bar M=
  \ffrac{1}{2p_+p_-}\!\!\sum_{j'=0}^{2p_+p_--1}\sum_{j=0}^{2p_+p_--1}
  \sum_{m=0}^{p_+-1}\sum_{n=0}^{p_+-1}
  \sum_{m'=0}^{p_--1}\sum_{n'=0}^{p_--1}
  \ffrac{(\qp^{-p_-}\!- \qp^{p_-})^{m+n}
    (\qm^{-p_+}\!- \qm^{p_+})^{m'+n'}}{
    [m]_+![m']_-![n]_+![n']_-!}\\*[-3pt]
  {}\times
  \qp^{(j - j') m + p_-\frac{m(m+1)}{2} - p_-\frac{n(n-1)}{2}}
  \qm^{(j' - j) m' + p_+\frac{m'(m'+1)}{2} - p_+\frac{n'(n'-1)}{2}}
  \q^{2 j j'}\\*[-3pt]
  {}\times
  \fp^{n}\ep^{m}\emi^{n'}\fm^{m'}\kk^{2j}\tensor
  \ep^{n}\fp^m\fm^{n'}\emi^{m'}\kk^{2j'},
\end{multline}
which involves only even powers of $k$, and therefore $\bar M\in
\XXX\tensor \XXX$.  To repeat this in words, the attempted calculation
of the $M$-matrix for~$\bar D(\mathscr{H})$ has yielded the $M$-matrix
for its $\XXX$ subalgerba (and hence $\bar D(\mathscr{H})$ is not a
factorizable quantum group; as noted in the introduction, the relevant
object is the pair $\XXX\subset\bar D(\mathscr{H})$ of quantum groups
one of which has an $M$-matrix and the other an
$R$-matrix).\footnote{This is rather natural from the conformal field
  theory standpoint: the $R$-matrix acting on two vertex operators
  describes their braiding, while the $M$-matrix describes its
  ``square,'' the operation of taking one of the operators
  \textit{around} the other.  The braiding requires choosing one of
  the two possible directions, which corresponds to choosing how
  $\kk=\sqrt{K}$ acts on the representations.}  We next use $\bar M$
to construct the Drinfeld map for~$\XXX$.

\subsection{The Drinfeld map for $\smash{\XXX}$}\label{sec:Drinfeld}
\subsubsection{Drinfeld map of the Grothendieck ring} 
In considering the Drinfeld map,\pagebreak[3] we first restrict it to
the image of the Grothendieck ring in the space $\Ch$ of
$q$-characters under the map $\XX\mapsto\qTr_{\XX}$ (see~\eqref{qCh}),
that is, consider the map
\begin{equation}\label{eq:Dr-hom}
  \Drinfeld\circ\qTr:
  \Grring(A)
  \to
  \cZ(A),
\end{equation}
\textit{which is a homomorphism of associative commutative algebras}
for any factorizable ribbon Hopf algebra~$A$.

With the balancing element for $A=\XXX$ given by~\eqref{balancing},
Eq.~\eqref{eq:Dr-hom} becomes
\begin{equation}\label{cchi-def}
  \begin{split}
    \Grring_{2p_+ p_-}&\to\ZZ\\    
    \XX^{\alpha}_{r,r'}&\mapsto\Drinfeld^{\alpha}(r,r')
    \equiv \Drinfeld(\gamma^{\alpha}(r,r'))
    = (\Tr_{\XX^{\alpha}_{r,r'}}\tensor\id)
    \bigl((K^{p_- - p_+}\tensor\one)\,
    \bar M
    \bigr),
  \end{split}
\end{equation}
where $1\leq r\leq p_+$, $ 1\leq r'\leq p_-$, and $\alpha=\pm$.
Clearly, $\Drinfeld^{+}(1,1)=\one$.  The other~$\Drinfeld^\pm(r,r')$
are evaluated as follows.

\begin{prop}\label{prop-eval}
  For $1\leq r\leq p_+$ and $ 1\leq r'\leq p_-$,
  \begin{equation}\label{the-cchi}
    \begin{split}
      \Drinfeld^{+}(r,r')&=\Drinfeld_{(+)}(r)\Drinfeld_{(-)}(r'),\\
      \Drinfeld^{-}(r,r')&=(-1)^{p_- + p_+}\Drinfeld^{+}(r,r')
      K^{p_+ p_-},
    \end{split}
  \end{equation}
  where
  \begin{multline}\label{cchi-pm}
    \Drinfeld_{(\pm)}(r)=
    \smash[t]{(-1)^{r-1}\sum_{a=0}^{r-1}\sum_{m=0}^{a}}
    (\q_{\pm}^{p_\mp}\!- \q_{\pm}^{-p_\mp})^{2m}
    \q_{\pm}^{p_\mp m (m + r  - 2a) + p_\mp(r-1-2a)}\\*
    {}\times{\qbin{r\!-\!a\!+\!m\!-\!1}{m}}_\pm {\qbin{a}{m}}_\pm
    e_{\pm}^{m} f_{\pm}^m K_\pm^{-m-r+1+2a}.
  \end{multline}
  Moreover, for $1\leq r\leq p_+$, $ 1\leq r'\leq p_-$, and
  $\alpha=\pm$,\ \ $\Drinfeld^{\alpha}(r,r')$ decomposes in terms of
  the Radford central elements as \textup{(}recall the notation
  in~\eqref{Proj-notation}--\eqref{hrhoLeft2}
  \textup{)}
  \begin{multline}\label{chi-decompose}
    \Drinfeld^{\alpha}(r,r')
    = \ffrac{r r'}{\iippA}\,\RadfordProj(p_+, p_-)
    + (-1)^{r p_- + r' p_+ + \beta p_+p_-}
    \ffrac{r r'}{\iippA}\,\RadfordProj(0, p_-)
    \\
    \shoveleft{\quad{}- 
      r'\sum_{s=1}^{p_+ - 1}
      (-1)^{(r - 1)p_- + (\beta p_- + r')(s + p_+)}
      \ffrac{\qp^{p_- r s}\!- \qp^{-p_- r s}}{\iippA}\,
      \NRadfordNESW(s, p_-)}
    \\
    \shoveright{{} +
      r r'\sum_{s=1}^{p_+ - 1} (-1)^{r p_- + (\beta p_- - r')(p_+ - s)}
      \ffrac{\qp^{p_- r s}\!+ \qp^{-p_- r s}}{\iippA}\,
      \RadfordProj(s, p_-)}\\
    \shoveleft{\quad
      {}-  r\sum_{s'=1}^{p_- - 1}
      (-1)^{(\beta p_+ + r)(s' + p_-) + (r' - 1)p_+}
      \ffrac{\qm^{p_+ r' s'}\!- \qm^{-p_+ r' s'}}{\iippA}\,
      \NRadfordNWSE(p_+, s')}\\
    \shoveright{{}+ 
      r r'\sum_{s'=1}^{p_- - 1} (-1)^{r' p_+ + (\beta p_+ - r)(p_- - s')}
      \ffrac{\qm^{p_+ r' s'}\!+\qm^{-p_+ r' s'}}{\iippA}\,
      \RadfordProj(p_+, s')}\\
    \shoveleft{{}
      +\!\!\sum_{(s,s')\in\setR}\!\!\!
      (-1)^{(\beta p_+ + r - 1)s' + (\beta p_- + r' - 1)s}
      \ffrac{(\qp^{p_- r s}\!-\!\qp^{-p_- r s})
        (\qm^{p_+ r' s'}\!-\!\qm^{-p_+ r' s'})}{
        \iippA}\,
      \NRadfordUpUp(s, s')}
    \\
    {} -
    \!\!\sum_{(s,s')\in\setR}\!\!\!
    (-1)^{(\beta p_+ - r) s' + (\beta p_- - r') s}
    \Bigl(    
    r'
    \ffrac{(\qm^{p_+ r' s'}\!+\!\qm^{-p_+ r' s'})
      (\qp^{p_- r s}\!-\!\qp^{-p_- r s})}{\iippA}\,
    \NRadfordS(s,s')
    \\*
    \shoveright{{}+
      r
      \ffrac{(\qp^{p_- r s}\!+\!\qp^{-p_- r s})
        (\qm^{p_+ r' s'}\!-\!\qm^{-p_+ r' s'})}{\iippA}\,
      \NRadfordB(s,s')
      \Bigr)}
    \\
    {} + 
    r r'\sum_{(s,s')\in\setR}\!\!
    (-1)^{(\beta p_+ - r) s' + (\beta p_- - r') s}
    \ffrac{(\qp^{p_- r s}\!+ \qp^{-p_- r s})
      (\qm^{p_+ r' s'}\!+ \qm^{-p_- r' s'})}{\iippA}\,
    \RadfordProj(s,s'),
  \end{multline}
  where $\alpha = (-1)^\beta$ \textup{(}that is, $\beta=0$ for
  $\Drinfeld^{+}(r,r')$ and $\beta=1$ for
  $\Drinfeld^{-}(r,r')$\textup{)}.
\end{prop}
  
Because~\eqref{eq:Dr-hom} is an algebra homomorphism, the product
$\Drinfeld^{\alpha}(r,r')\Drinfeld^{\alpha'}(s,s')$ is expressed
through the $\Drinfeld^{\pm}(t,t')$ just by the formulas
in~\bref{thm:Gr-ring} with $\XX^\pm\to\Drinfeld^\pm$ and
$\tXX^\pm\to\tilde\Drinfeld^\pm$ (cf.~\cite{[FGST]}).

The proof of~\eqref{the-cchi}--\eqref{cchi-pm} is by straightforward
calculation based on~\eqref{eq:EmFm-prod2} and~\eqref{eq:EmFm-prod}.
Using~\eqref{on-vectors}, we then find that
\begin{multline*}
  \Tr_{\XX^{\alpha}(r,r')}\bigl(
  \fp^{m}\ep^{m}\emi^{m'}\fm^{m'} K^{j + p_- - p_+}\bigr)={}\\
  \shoveleft{\quad{}
    =\sum_{a=m}^{r-1}\sum_{a'=0}^{r'-1-m'}
    \alpha^{m p_- + m' p_+ + j + p_- - p_+}(-1)^{m(r'-1)+m'(r-1)+r+r'}
    ([m]_+![m']_-!)^2}\\*
  {}\times{\qbin{r\!-\!a\!+\!m\!-\!1}{m}}_+ {\qbin{a}{m}}_+
  {\qbin{r'\!-\!a'\!-\!1}{m'}}_- {\qbin{a'\!+\!m'}{m'}}_-
  \qp^{(r-1-2a)(j + p_-)}
  \qm^{(r'-1-2a')(j - p_+)},
\end{multline*}
which after some rearrangements
gives~\eqref{the-cchi}--\eqref{cchi-pm} in accordance
with~\eqref{cchi-def} for the $M$-matrix in~\eqref{bar-M}.  The
decomposition then follows by calculating the action on projective
modules, as explained in~\bref{sec:coeffs}.

\begin{rem}\label{Drinfeld-factors}
  Unlike the $M$-matrix, the Drinfeld map in~\eqref{the-cchi} factors
  into a product of the Drinfeld maps~\eqref{cchi-pm}
  for the two $\overline{\mathscr{U}}_{Q_{\pm}}s\ell(2)$ quantum
  groups with $Q_{\pm}=\q_{\pm}^{p_{\mp}}$.  The corresponding
  multiplication formulas are known from~\cite{[FHST],[FGST]}:
  \begin{equation}\label{the-fusion-12}
    \Drinfeld_{(\pm)}(r)\Drinfeld_{(\pm)}(r')
    =\smash[b]{\sum_{\substack{r''=|r - r'| + 1\\
        \mathrm{step}=2}}^{r + r' - 1}}
    \tilde\Drinfeld_{(\pm)}(r''),
  \end{equation}
  where
  \begin{equation*}
    \tilde\Drinfeld_{(\pm)}(r)
    =
    \begin{cases}
      \Drinfeld_{(\pm)}(r),&1\leq r\leq p_\pm,\\
      \Drinfeld_{(\pm)}(2p_\pm\!-\!r)
      + 2 (-1)^{p_- + p_+}\Drinfeld_{(\pm)}(r\!-\!p_\pm)K_\pm^{p_\pm},
      & p_\pm\!+\!1 \leq r
      \leq 2p_\pm\!-\!1.
    \end{cases}
  \end{equation*}
  We recall that $K_+^{p_+}=K_-^{p_-}=K^{p_+ p_-}$ in~$\XXX$.
  
  Naturally, relations~\eqref{the-fusion-12} become the corresponding
  formulas in~\cite{[FGST]} in the cases where either $p_+ = 1$ or
  $p_- = 1$.  
\end{rem}

\subsubsection{Drinfeld images of pseudotraces} We next evaluate the
Drinfeld map on the $\gamma^{\NESW}$, $\gamma^{\NWSE}$, and
$\gamma^{\UpUp}$ in~\eqref{eq:gamma-basis}. This gives
\begin{multline}\label{chiNWSE-prepare}
  \NDrinfeldNWSE(r,r')=\Drinfeld(\Ngamma^{\NESW}(r,r'))={}\\
  {}= (-1)^{r'}
  (\Tr_
  {\mathbb{P}_{r,r'}}
  \tensor\id)
  \bigl((K^{p_- - p_+}\tensor\one)\,
  \bar M (\NsigmaNESW(r,r')\tensor\id)\bigr)
\end{multline}
\vspace{-\baselineskip}
\begin{multline}\label{chiNESW-prepare}
  \NDrinfeldNESW(r,r') =\Drinfeld(\Ngamma^{\NWSE}(r,r'))\\*
  {}=(-1)^{r}
  (\Tr_
  {\mathbb{P}_{r,r'}}
  \tensor\id)
  \bigl((K^{p_- - p_+}\tensor\one)\,
  \bar M (\NsigmaNWSE(r,r')\tensor\id)\bigr)
\end{multline}
for all $1\leq r\leq p_+\!-\!1$, $1\leq r'\leq p_-\!-\!1$,
and
\begin{multline}
  \NDrinfeldDnDn(r,r')=\Drinfeld(\Ngamma^{\UpUp}(r,r'))\\
  {}=(-1)^{r+r'}
  (\Tr_{\mathbb{P}_{r,r'}}
  \tensor\id)
  \bigl((K^{p_- - p_+}\tensor\one)\,\bar M
  (\NsigmaUpUp(r,r')\tensor\id) \bigr)
\end{multline}
for $(r,r')\in\setR$.

\begin{prop}\label{prop:skew-drinfeld-decompose}
  The central elements $\DrinfeldNWSE(r,r')$, $\DrinfeldNESW(r,r')$,
  and $\DrinfeldDnDn(r,r')$ decompose in terms of the Radford central
  elements as
  \begin{multline}\label{chiNWSE-decompose}
    (-1)^{r'}\NDrinfeldNWSE(r,r')
    ={} 
     r'\sum_{s=1}^{p_+ - 1} (-1)^{r'(s + p_+) + p_- r}\,
    \ffrac{\qp^{p_- r s}\!-\!\qp^{-p_- r s}}{p_- \iiA}\,
    \hrhoS(s,p_-)
    \\
    \qquad{}+
    r'\!\!\!\!\sum_{(s,s')\in\setR}\!\!  (-1)^{r s' + r' s}
    \ffrac{(\qp^{p_- r s}\!-\!\qp^{-p_- r s})
      (\qm^{p_+ r' s'}\!+ \qm^{-p_+ r' s'})}{p_- \iiA}\,
    \hrhoS(s,s')
    \\
    {}-
    \!\!\sum_{(s,s')\in\setR}\!\! (-1)^{r s' + r' s}
    \ffrac{(\qp^{p_- r s}\!-\!\qp^{-p_- r s})
      (\qm^{p_+ r' s'}\!- \qm^{-p_+ r' s'})}{p_- \iiA}\,
    \NvvarphiNWSE(s,s'),
  \end{multline}
  \vspace{-.9\baselineskip}
  \begin{multline}\label{chNWSE(r,p_-)-decompose}
    (-1)^{p_-}\NDrinfeldNWSE(r,p_-)
    {}=
    \ffrac{\qp^{p_-}\!-\!\qp^{-p_-}}{\iiA}\Biggl(
    \sum_{s=1}^{p_+-1}
    (-1)^{p_-(s + p_+ +  r)}
    {[r s]}_+
    \hrhoS(s,p_-)
    \\*
    {}+2\!\!\!\sum_{(s,s')\in\setR}\!\!\!
    (-1)^{(r+p_+) s' + p_- s}
    [r s]_+
    \hrhoS(s,s')
    \!\Biggr),
  \end{multline}
  \vspace{-.9\baselineskip}
  \begin{multline}
    \label{chiNESW-decompose}
    (-1)^{r}
    \NDrinfeldNESW(r,r')
    ={} 
    r\sum_{s'=1}^{p_- - 1} (-1)^{r(s' + p_-) + p_+ r'}\,
    \ffrac{\qm^{p_+ r' s'}\!-\!\qm^{-p_+ r' s'}}{p_+ \iiA}\,
    \hrhoB(p_+, s')
    \\
    \qquad{} +
    r\!\!\!\sum_{(s,s')\in\setR}\!\!  (-1)^{r' s + r s'}
    \ffrac{(\qp^{p_- r s}\!+ \qp^{-p_- r s})
      (\qm^{p_+ r' s'}\!-\!\qm^{-p_+ r' s'})}{p_+ \iiA}\,
    \hrhoB(s,s')
    \\
    {}-
    \!\!\sum_{(s,s')\in\setR}\!\!(-1)^{r' s + r s'}
    \ffrac{(\qp^{p_- r s}\!- \qp^{-p_- r s})
      (\qm^{p_+ r' s'}\!-\!\qm^{-p_+ r' s'})}{p_+ \iiA}\,
    \NvvarphiNESW(s,s'),
  \end{multline}
  \vspace{-.9\baselineskip}
  \begin{multline}\label{chNESW(p_+,r')-decompose}
    (-1)^{p_+}
    \NDrinfeldNESW(p_+,r')
    =
    \ffrac{\qm^{p_+}\!-\!\qm^{-p_+}}{\iiA}\Biggl(
    \sum_{s'=1}^{p_--1}
    (-1)^{p_+(s' + p_- +  r')}\,
    {[r' s']}_- 
    \hrhoB(p_+, s')
    \\*
    {}+2\!\!\!\sum_{(s,s')\in\setR}\!\!\!
    (-1)^{(r'+p_-) s + p_+ s'}
    [r' s']_-\,
    \hrhoB(s,s')
    \!\Biggr),
  \end{multline}
  \vspace{-.9\baselineskip}
  \begin{multline}
    \label{chiDnDn-decompose}
    (-1)^{r+r'}\NDrinfeldDnDn(r,r')
    ={}
    \!\!\!\sum_{(s,s')\in\setR}\!\!(-1)^{r s' + r' s}
    \ffrac{(\qp^{p_- r  s}\!- \qp^{-p_- r  s})
      (\qm^{p_+ r'  s'}\!- \qm^{-p_+ r'  s'})}{
      \iiA}\,\vvarphi(s,s'),
  \end{multline}
  where we use the notation in~\bref{Radford-notation}
  and set
  \begin{equation}\label{vvarphiNWSE+}
    \begin{aligned}
      \NvvarphiNWSE(r,r')&=
      (-1)^r(p_+\!-\!r)\,
      \NRadfordNWSE(r, r')
      +(-1)^{p_+ + r} r\,
      \NRadfordNWSE(p_+\!-\!r,p_-\!-\!r'),
      \\    
      \NvvarphiNESW(r,r')&=
      (-1)^{r'}(p_-\!-\!r')\,    
      \NRadfordNESW(r, r')
      +(-1)^{p_- + r'} r'\,
      \NRadfordNESW(p_+\!-\!r, p_-\!-\!r'),
    \end{aligned}    
  \end{equation}
  for $ (r,r')\in\setR$.
\end{prop}

\begin{proof}
  A calculation similar to the one in the proof of~\bref{prop-eval}
  shows that the $\DrinfeldNWSE(r,r')$, $\DrinfeldNESW(r,r')$, and
  $\DrinfeldDnDn(r,r')$ are explicitly given by
  \begin{align}\label{chiNWSE-eval}
    \NDrinfeldNWSE(r,r')
    &=
    (-1)^{r'}\bigl(\Nttheta_{(+)}(r)
    -(-1)^{p_+ + p_-}\Nttheta_{(+)}(p_+\!-\!r)K^{p_+ p_-}\bigr)
    \Drinfeld_{(-)}(r')\\
    \intertext{for $1\leq r\leq p_+-1$ and $1\leq r'\leq p_-$,}
    \NDrinfeldNESW(r,r')&=
    (-1)^{r}\Drinfeld_{(+)}(r)\bigl(\Nttheta_{(-)}(r')
    -(-1)^{p_+ + p_-}\Nttheta_{(-)}(p_-\!-\!r')K^{p_+ p_-}\bigr),
    \label{chiNESW-eval}\\
    \intertext{for $1\leq r\leq p_+$ and $1\leq r'\leq p_--1$, and}
    \label{chiDnDn-eval}
    \NDrinfeldDnDn(r,r')&=
    (-1)^{r+r'}\bigl(\Nttheta_{(+)}(r)
    -(-1)^{p_+ + p_-}\Nttheta_{(+)}(p_+\!-\!r)K^{p_+ p_-}\bigr)\\*
    &\qquad\quad{}\times\bigl(\Nttheta_{(-)}(r')
    -(-1)^{p_+ + p_-}\Nttheta_{(-)}(p_-\!-\!r')K^{p_+ p_-}\bigr),
    \notag
  \end{align}
  for $(r,r')\in\setR$, where $\Drinfeld_{(\pm)}$ are defined
  in~\eqref{cchi-pm} and
  \begin{multline*}
    \Nttheta_{(\pm)}(r)= (-1)^{r}[r]_{\pm}
    \sum_{a=0}^{r-1}\sum_{m=0}^{p_\pm-1}
    \ffrac{(\q_{\pm}^{p_\mp}\!- \q_{\pm}^{-p_\mp})^{2m-1}}{
      ([m]_\pm!)^2}\,
    \q_{\pm}^{p_\mp m(m+r-2a) + p_\mp(r-1-2a)}
    \\*
    {}\times [x^1]\bigl(\CCpm{r,a}^{m}(x)\bigr)
    e_{\pm}^{m} f_{\pm}^m K^{-p_\mp (m+r-1-2a)},
  \end{multline*}
  where the polynomials $\CCpm{r,a}^{m}(x)$ are defined
  in~\eqref{CCpolynom}.
  
  The rest is shown by evaluation of the action of $\DrinfeldNWSE$,
  $\DrinfeldNESW$, and $\DrinfeldDnDn$ (that is, of
  $\Nttheta_{(\pm)}(r)$ and $\Drinfeld_{(+)}(r)$) on projective
  modules, see~\bref{sec:coeffs}.  Because of the Drinfeld map
  factorization (cf.~\bref{Drinfeld-factors}), the demonstration
  largely reduces to separate calculations in each of the two
  $\overline{\mathscr{U}}_{Q_{\pm}}s\ell(2)$ quantum groups with
  $Q_{\pm}=\q_{\pm}^{p_{\mp}}$, cf.~\cite{[FGST]}.
\end{proof}

\subsubsection{Grothendieck ring as a polynomial quotient
  ring}\label{sec:Grring-U} We next characterize the Groth\-endieck
ring as follows.  Let 
\begin{gather}\label{eq:chebyshev-sin}
  \cheb_s(2\cos t)=\mfrac{\sin s t}{\sin t},\quad s\geq1,
\end{gather}
be the Chebyshev polynomials of the second kind.  They satisfy (and
are determined by) the recursive relation
\begin{equation}\label{cheb-rec}
  x\cheb_s(x)=\cheb_{s-1}(x)+\cheb_{s+1}(x),
  \quad s \geq 2,
\end{equation}
with the initial data $\cheb_1(x)=1$, $\cheb_2(x)=x$.

\begin{prop}\label{prop:Cheb}
  The Grothendieck ring $\Grring_{2 p_+ p_-}$ of $\XXX$ is the
  quotient of \ $\oC[x,y]$ over the ideal generated by the polynomials
  \begin{align*}
    &\cheb_{2p_+ + 1}(x) - \cheb_{2p_+ - 1}(x) - 2,\\
    &\cheb_{2p_- + 1}(y) - \cheb_{2p_- - 1}(y) - 2,\\
    &\cheb_{p_+ + 1}(x) - \cheb_{p_+ - 1}(x)
    -\cheb_{p_- + 1}(y) + \cheb_{p_- - 1}(y).
  \end{align*}
  Moreover, the basis of irreducible representations
  $\XX^{\pm}_{r,r'}$ is given by \textup{(}the image under the
  quotient mapping of\textup{)} the respective polynomials
  $P^{\pm}_{r,r'}$, $1\leq r\leq p_+$, $1\leq r'\leq p_-$, where
  \begin{align*}
    P^+_{r,r'}(x,y)&=\cheb_{r}(x)\cheb_{r'}(y),\\
    P^-_{r,r'}(x,y)&
    =\Bigl(\half\,\cheb_{p_+ + r}(x)-\half\,\cheb_{p_+ - r}(x)\Bigr)
    \cheb_{r'}(y).
  \end{align*}
\end{prop}

To show~\bref{prop:Cheb}, we first of all recall, once again, that the
Drinfeld map $\Drinfeld$ is an algebra isomorphism.  Next, it follows
from~\eqref{the-Groth} that the Grothendieck ring is
\textit{generated} by the elements
\begin{equation*}
  \Drinfeld^{+}(2,1)
  =\cas_+,\qquad
  \Drinfeld^{+}(1,2)
  = \cas_-
\end{equation*}
(see~\eqref{eq:casimir-plus} and~\eqref{eq:casimir-minus}), which,
moreover, satisfy Eqs.~\eqref{psi-pm}, where we can rewrite
$\psi_{\pm}$ as
\begin{align*}
  \psi_+(x)
  &=\cheb_{2p_+ + 1}(x) - \cheb_{2p_+ - 1}(x) - 2,\\
  \psi_-(y)
  &=\cheb_{2p_- + 1}(y) - \cheb_{2p_- - 1}(y) - 2,
\end{align*}
and the equation
\begin{equation}\label{U=U}
  \cheb_{p_+ + 1}(\cas_+) - \cheb_{p_+ - 1}(\cas_+)
  =\cheb_{p_- + 1}(\cas_-) - \cheb_{p_- - 1}(\cas_-)
\end{equation}
(where both sides are actually equal to $(-1)^{p_+ + p_-} 2 K^{p_+
  p_-}$).  The statement in~\bref{prop:Cheb} is now immediate.

\subsection{The $\smash{\XXX}$ ribbon element}\label{XXX-ribbon} 
We recall that a \textit{ribbon Hopf algebra}~\cite{[RSts]} is a
quasitriangular Hopf algebra equipped with an invertible central
element $\ribbon$, called the \textit{ribbon element}, such that
\begin{equation*}
  S(\ribbon)=\ribbon,\quad\epsilon(\ribbon)=1,
  \quad
  \Delta(\ribbon)=M^{-1}(\ribbon\tensor\ribbon).
\end{equation*}
In a ribbon Hopf algebra, $\ribbon^2= \sqs S(\sqs)$ and, in fact,
\begin{equation}\label{balance-ribbon}
  \ribbon=\sqs\balance^{-1},
\end{equation}
where $\balance$ is the balancing element and 
\begin{gather}\label{canon-sqs}
  \sqs = \cdot\bigl((S\tensor\id)R_{21}\bigr)
\end{gather}
(where $\cdot(a\tensor b)=ab$) is the \textit{canonical
  element}~\cite{[Drinfeld]}.

\begin{prop}\label{prop:ribbon}
  The $\XXX$-ribbon element $\ribbon$ is given by
  \begin{equation}\label{ribbon-factors}
    \ribbon = \bar\ribbon\,\ribbon^*,
  \end{equation}
  where
  \begin{equation}\label{bar-ribbon}
    \bar\ribbon=
    \sum_{(r,r')\in\setI}\!\!\!
    e^{2i\pi\Delta_{r,r'}}\Idem(r,r'),
  \end{equation}
  with
  \begin{equation}\label{Delta-rs}
    \Delta_{r,r'}=\ffrac{(p_+r'-p_-r)^2-(p_+-p_-)^2}{4p_+p_-}
  \end{equation}
  \textup{(}and $\Idem(r,r')$ being the primitive idempotents in the
  center\textup{)}, and
  \begin{equation*}
    \ribbon^*
    = \Bigl(\one
    + \ffrac{1}{p_+}\,
    \NDrinfeldNWSE(1, 1)
    \Bigr)
    \Bigl(\one
    + \ffrac{1}{p_-}\,    
    \NDrinfeldNESW(1, 1)
    \Bigr).
  \end{equation*}
\end{prop}

\begin{proof}
  First, using the $R$-matrix in~\eqref{q-double-R}, we explicitly
  calculate the canonical element~$\sqs$ in $\bar D(\mathscr{H})$ as
  \begin{multline*}
    \sqs\! = \! \ffrac{1 + i}{2\sqrt{\!p_+ p_-\mathstrut}}\!
    \sum_{m=0}^{p_+-1}\sum_{r=0}^{p_+ - 1}
    \sum_{n=0}^{p_--1}\sum_{s=0}^{p_- - 1}\!
    \ffrac{(\qp^{p_-}\!-\!\qp^{-p_-})^m
      (\qm^{p_+}\!-\!\qm^{-p_+})^n\!\!}{{[m]_+! [n]_-!}}\,
    \qp^{p_-\left(\!\frac{m(m+3)}{2} - \frac{r^2}{2}\right)}
    \qm^{p_+\left(\!\frac{n(n+3)}{2} - \frac{s^2}{2}\right)}\\*
    {}\times(-1)^{r s} \bigl(\fp^{m} K_+^{r - m} \ep^{m}\bigr)
    \bigl(\one + (-1)^{p_+ s + p_- r}(-i)^{p_+ p_-} K^{p_+ p_-}\bigr)
    \bigl(\emi^{n} K_-^{s + n} \fm^{n}\bigr).
  \end{multline*}
  The calculation involved the simple Gaussian sum
  \begin{equation*}
    \sum_{j=0}^{2p_+p_--1}
    e^{\frac{i\pi j^2}{2 p_+ p_-}}=(1+ i)\sqrt{p_+ p_-\mathstrut}.
  \end{equation*}
  The expression obtained shows, in particular, that $\sqs$ is an
  element of $\XXX$, not of $\bar D(\mathscr{H})$.  We then find the
  ribbon element in accordance with~\eqref{balance-ribbon}
  and~\eqref{balancing}, as
  \begin{equation*}
    \ribbon=\sqs K^{p_- - p_+}.
  \end{equation*}
  
  We next decompose the ribbon element $\ribbon$ with respect to the
  basis of the Radford central elements.  First, with the explicit
  expression for $\ribbon$, we find its action on the irreducible
  representations $\XX^+_{r,r'}$.  An elementary calculation shows
  that the action of the ribbon element on irreducible modules is
  given by
  \begin{align*}
    \ribbon\bigr|_{\XX_{r,r'}^{+}} =
    (-1)^{r r' + 1}\qp^{p_-\frac{r^2 - 1}{2}}
    \qm^{p_+ \frac{r'^2 - 1}{2}}
    &= e^{2i\pi\Delta_{r,r'}},
    \\
    \ribbon\bigr|_{\XX_{p_+,p_-}^{-}}
    &=e^{2 i \pi \Delta_{0, p_-}},
  \end{align*}
  thus associating $\Delta_{r,r'}$ with each pair $(r,r')\in\setI$.
  
  With some more work (as for other decompositions obtained above; also
  cf.~\cite{[FGST]}), we calculate the nondiagonal part of the ribbon
  element action on projective modules and thus find the coefficients in
  the decomposition of the ribbon element with respect to the canonical
  central elements, as explained in~\bref{sec:coeffs}.  With the
  notation in~\bref{Radford-notation}, 
  we have
  \begin{multline}\label{ribbon-decompose}
    \ribbon =\sum_{(r,r')\in\setI}\!\!\!
    e^{2i\pi\Delta_{r,r'}}\Idem(r,r')+{}\\
    +\!\sum_{(r,r')\in\setR}\!\!\!  (-1)^r e^{2i\pi\Delta_{r,r'}}
    \ffrac{\qm^{p_+ r'}\!- \qm^{-p_+ r'}}{ 4p_-^2}\bigl(
    r'\vSW(r,r')\!-\!(p_-\!-\!r')\vNE(r,r')\bigr)\\
    \shoveright{{}+ \!\sum_{(r,r')\in\setR}\!\!\!  (-1)^{r'}
      e^{2i\pi\Delta_{r,r'}} \ffrac{\qp^{p_- r}\!- \qp^{-p_-
          r}}{4p_+^2}
      \bigl(r\vSE(r,r')\!-\!(p_+\!-\!r)\vNW(r,r')\bigr)}\\
    \shoveright{{}+\!\sum_{(r,r')\in\setR}\!\!\!  (-1)^{r+r'}
      e^{2i\pi\Delta_{r,r'}} \ffrac{(\qp^{p_- r}\!-\!\qp^{-p_- r})
        (\qm^{p_+ r'}\!-\!\qm^{-p_+ r'})}{ \iippA}\, \vvarphi(r,r')}
    \\
    \shoveright{{}-\sum_{r'=1}^{p_- - 1} e^{2i\pi\Delta_{p_+,r'}}
      \ffrac{(-1)^{p_+ + p_- + r'} (\qm^{p_+ r'}\!- \qm^{-p_+ r'})}{
        \iippA}\, \hrhoB(p_+,r')\quad}
    \\*
    {}-\sum_{r=1}^{p_+ - 1} e^{2i\pi\Delta_{r,p_-}} \ffrac{(-1)^{p_-
        + p_+ + r} (\qp^{p_- r}\!- \qp^{-p_- r})}{\iippA}\,
    \hrhoS(r,p_-).
  \end{multline}
  
  The formulas in~\bref{prop:radford-decompose} of course allow
  expressing this through the Radford-map images, similarly to other
  decompositions of central elments, but the above form involving
  primitive idempotents shows that $\ribbon$ factors into semisimple
  and unipotent (one-plus-nilpotent) parts as $\ribbon =
  \bar\ribbon\,\ribbon^*$, with $\bar\ribbon$ just as in the
  proposition and $\ribbon^*$ given by $\one$ plus all the terms
  in~\eqref{ribbon-decompose} except the first sum, with all the
  factors $e^{2i\pi\Delta_{r,r'}}$ replaced by unity.  The formulas
  in~\bref{sec:Drinfeld} then make it easy to verify that
  \begin{equation*}
    \ribbon^*
    = \one
    + \ffrac{1}{p_+ p_-}\,
    \NDrinfeldDnDn(1, 1)
    + \ffrac{1}{p_-}\,
    \NDrinfeldNESW(1, 1)
    + \ffrac{1}{p_+}\,
    \NDrinfeldNWSE(1, 1),
  \end{equation*}
  whence the formula in the proposition follows immediately.
\end{proof}

\section{$\SLiiZ$-representations on the $\XXX$
  center}\label{sec:SL2Z} The aim of this section is to prove
Theorem~\bref{thm:q-structure}.  We define the $\SLiiZ$-action on the
$\XXX$ center $\cZ$ as
\begin{equation}\label{theST}
  \begin{split}
    \modS:{}& a\mapsto \Radford\bigl(\Drinfeld^{-1}(a)\bigr),\\
    \modT:{}& a\mapsto\AT\,\modS(\ribbon\modS^{-1}(a)),
  \end{split}
\end{equation}
which follows~\cite{[Lyu],[LM],[Kerler]} with insignificant
variations.  Here,
\begin{equation*}
  c=13-6\ffrac{p_+}{p_-} - 6\ffrac{p_-}{p_+},
\end{equation*}
which is the central charge of the $(p_+,p_-)$ conformal field theory
model.

\subsection{$\SLiiZ$-adapted bases in the
  center}\label{sec:SL2Z-acts} We have two ``distinguished'' bases
in the center $\cZ$ of $\XXX$, associated with the Radford and the
Drinfeld maps.  The ``Radford'' basis is given
in~\eqref{eq:radford-basis}.  The ``Drinfeld'' basis consists of the
$2 p_+
p_-\,{+}\,(p_+\,{-}\,1)p_-\,{+}$\linebreak[0]$\,p_+(p_-\,{-}\,1)\,{+}\,
\half(p_+\,{-}\,1)(p_-\,{-}\,1)$ elements
\begin{align*}
  &\Drinfeld^{\pm}(r,r'),
  \;\quad 1\leq r\leq p_+,\quad 1\leq r'\leq p_-,\\
  &\DrinfeldNWSE(r,r'),\
  \begin{array}[t]{l}
    1\leq r\leq p_+\!-\!1,\\
    1\leq r'\leq p_-,
  \end{array}\qquad
  \DrinfeldNESW(r,r'),\
  \begin{array}[t]{l}
    1\leq r\leq p_+,\\
    1\leq r'\leq p_-\!-\!1,
  \end{array}
  \\
  &\DrinfeldDnDn(r,r'),\quad (r,r')\in\setR.
\end{align*}

{}From the definitions, the $\modS$ mapping relates the two bases as
\begin{align}
  \modS^{-1}\bigl(\Radford^\pm(r,r')\bigl)
  &=\Drinfeld^\pm(r,r'),\label{S-pm}\\
  \modS^{-1}\bigl(\NRadfordNESW(r,r')\bigr)
  &=\NDrinfeldNWSE(r,r'),
  \label{S-RadfordNESW}\\
  \modS^{-1}\bigl(\NRadfordNWSE(r,r')\bigr)  
  &=\NDrinfeldNESW(r,r'),
  \\
  \modS^{-1}\bigl(\NRadfordUpUp(r, r')\bigr)
  &= \NDrinfeldDnDn(r, r').\label{S-DnDn}
\end{align}
In accordance with \cite{[Lyu]},
\begin{equation*}
  \modS^2\bigm|_{\cZ}=\id,
\end{equation*}
which, clearly, requires a proper choice of normalization in our
formulas, most easily checked by calculating
$\modS^{-1}(\one)=\modS^{-1}\bigl(\Drinfeld^{+}(1,1)\bigr)=
\Radford^+(1,1) =\coint$.
(As we have noted, this fixes the value of $\zeta$
in~\eqref{Rad-normalization}.)

\subsection{Modular group action on the
  center}\label{sec:modular-eval} We now prove the most tedious part
of Theorem~\bref{thm:q-structure}, the structure of the $\SLiiZ$
representation on $\cZ$ in Claim~\ref{claim:decomp}.  The proof
amounts to a series of lemmas.

\begin{lemma}\label{lemma:R-min}
  For $(r,r')\in\setR$, the elements
  \begin{equation*}
    \NvvarphiJoin(r,r') = \NvvarphiNESW(r, r') - \NvvarphiNWSE(r, r')
  \end{equation*}
  \textup{(}see the notation in~\eqref{vvarphiNWSE+}\textup{)} span a
  $\half(p_+\!-\!1)(p_-\!-\!1)$-dimen\-sion\-al
  $\SLiiZ$-sub\-rep\-resen\-tation in~$\cZ$.  Moreover, this
  representation is isomorphic to $\Rmin$
  in~\eqref{W-structure}.\footnote{We remind the reader that this is
    the $\SLiiZ$ representation in the minimal $(p_+,p_-)$ conformal
    field theory model.}
\end{lemma}
\begin{proof}
  First, the formulas
  \begin{align*}
    \ribbon\,\NvvarphiNESW(r, r')    
    &= e^{2i\pi\Delta_{r,r'}}\NvvarphiNESW(r, r')    
    + e^{2i\pi\Delta_{r,r'}}\vvarphi(r,r'),\\
    \ribbon\,\NvvarphiNWSE(r, r')    
    &= e^{2i\pi\Delta_{r,r'}}\NvvarphiNWSE(r, r')    
    + e^{2i\pi\Delta_{r,r'}}\vvarphi(r,r')
  \end{align*}
  (which are a direct consequence of~\bref{prop:radford-decompose}
  and~\eqref{ribbon-decompose}) show that
  \begin{equation*}
    \ribbon\,\NvvarphiJoin(r, r')
    = e^{2i\pi\Delta_{r,r'}}\NvvarphiJoin(r, r').
  \end{equation*}
  Second, we find $\modS(\vvarphiNESW(r, r'))$ from
  Eqs.~\eqref{S-RadfordNESW} and~\eqref{chiNWSE-decompose}, which
  imply
  \begin{multline}\label{relevant}    
    \modS(\NvvarphiNESW(r, r'))=
    (-1)^{r'}\bigl((p_-\!-\!r') 
    \NDrinfeldNWSE(r, r')
    +(-1)^{p_-} r'
    \NDrinfeldNWSE(p_+\!-\!r, p_-\!-\!r')\bigr)\\
    = -\ffrac{1}{\iiA}
    \!\sum_{(s,s')\in\setR} \!\!\!
    (-1)^{r s' + r' s}
    (\qp^{p_- r s}\!-\!\qp^{-p_- r s})
    (\qm^{p_+ r' s'}\!- \qm^{-p_+ r' s'})
    \NvvarphiNWSE(s,s'),
  \end{multline}
  where we also used~\eqref{radfordNWSE-decompose}.  With a similar
  expression for $\modS(\vvarphiNWSE(r, r'))$, we have
  \begin{equation*}
    \modS\bigl(\NvvarphiJoin(r, r')\bigr)
    =\ffrac{1}{\iiA}
    \!\sum_{(s,s')\in\setR}\!\!\!
    (-1)^{r s' + r' s} (\qp^{p_- r s}\!- \qp^{-p_- r s})
    (\qm^{p_- r' s'}\!- \qm^{-p_- r' s'})
    \NvvarphiJoin(s, s').
  \end{equation*}
  The lemma is proved.
\end{proof}

We once again recall the notation in~\eqref{vvarphiNWSE+} and set
\begin{equation*}
  \Nppsi(r,r')= \NvvarphiNESW(r,r') + \NvvarphiNWSE(r,r').
\end{equation*}
\begin{lemma}\label{lemma:triplet}
  The elements $\RadfordUpUp(r,r')$, $\ppsi(r,r')$, and
  $\vvarphi(r,r')$ with $(r,r')\in\setR$ span a  
  $\frac{3}{2}(p_+\!-\!1)(p_-\!-\!1)$-dimensional
  $\SLiiZ$-subrepresentation $\mathscr{R}$ in~$\cZ$ and,
  moreover,
  \begin{equation*}
    \mathscr{R}=\oC^3\tensor\Rmin,
  \end{equation*}
  where $\oC^3$ is the symmetrized square of the standard
  two-dimensional representation and $\Rmin$ is the representation
  in~\eqref{W-structure}.
\end{lemma}

\begin{proof}
  Clearly,
  \begin{equation*}
    \ribbon\,\vvarphi(r,r')=e^{2i\pi\Delta_{r,r'}}\vvarphi(r,r')
  \end{equation*}
  from~\eqref{ribbon-decompose} and~\bref{prop:radford-decompose}.  We
  next evaluate the action of the ribbon element
  on~$\RadfordUpUp(r,r')$ using~\bref{prop:radford-decompose}, which
  gives
  \begin{equation*}
    \ribbon\NRadfordUpUp(r,r')=
    e^{2i\pi\Delta_{r,r'}}\bigl(\NRadfordUpUp(r,r')
    + (-1)^{r+r'}
    \vvarphi(r,r')
    {}+ (-1)^{r + r'}
    \Nppsi(r,r')\bigr).
  \end{equation*}
  Next, a similar calculation yields
  \begin{equation*}
    \ribbon\,\Nppsi(r,r')
    = e^{2i\pi\Delta_{r,r'}}\bigl(\Nppsi(r,r')
    +2 \vvarphi(r,r')\bigr).
  \end{equation*} 
  
  As regards the $\modS$ action, we recall~\eqref{S-DnDn}
  and~\eqref{chiDnDn-decompose}, which we rewrite as
  \begin{equation*}
    \modS^{-1}(\NRadfordUpUp(r,r'))=
    \ffrac{(-1)^{r+r'}}{\iiA}
    \sum_{(s,s')\in\setR}\!\!\!
    (-1)^{r s' + r' s}
    (\qp^{p_- r s}\!- \qp^{-p_- r s})
    (\qm^{p_+ r' s'}\!- \qm^{-p_+ r' s'})
    \vvarphi(s,s').
  \end{equation*}
  Next, repeating the calculation in the proof of~\bref{lemma:R-min},
  we obtain
  \begin{equation*}    
    \modS^{-1}(\Nppsi(r,r'))
    = -\ffrac{1}{\iiA}
    \!\sum_{(s,s')\in\setR}\!\!\!(-1)^{r s' + r' s}
    (\qp^{p_- r s}\!- \qp^{-p_- r s})
    (\qm^{p_+ r' s'}\!- \qm^{-p_+ r' s'})
    \Nppsi(s,s').
  \end{equation*}
  
  This completes the proof.
\end{proof}

\begin{lemma}\label{lemma:Lsh}
  For $(r,r')\in\setILeft$, the elements $\hrhoS(r,r')$
  \textup{(}see~\eqref{hrhoright} and~\eqref{hrhoright-bdr}\textup{)}
  and $\RadfordS(r,r')$ \textup{(}see~\eqref{hrhoRight} and
  \eqref{hrhoRight2}\textup{)} span a\pagebreak[3]
  $(p_+\,{-}\,1)(p_-\,{+}\,1)$-dimensional
  $\SLiiZ$-subrep\-resen\-tation $\RRRight$ in~$\cZ$.  Moreover, its
  structure is
  \begin{equation*}
    \RRRight=\oC^2\tensor\RRight,
  \end{equation*}
  where $\oC^2$ is the standard two-dimensional $\SLiiZ$
  representation and $\RRight$ is the representation
  in~\eqref{W-structure}.
\end{lemma}

\begin{proof}
  We first calculate the ribbon action.
  {}From~\eqref{ribbon-decompose} and~\bref{prop:radford-decompose},
  we have
  \begin{align*}    
    \ribbon\NRadfordS(r,r')
    &=e^{2i\pi\Delta_{r,r'}}    
    \NRadfordS(r,r')
    + e^{2i\pi\Delta_{r,r'}}    
    \hrhoS(r,r')\\[-6pt]
  \intertext{for $(r,r')\in\setR$, and}
    \ribbon\NRadfordS(r,p_-)
    &=e^{2i\pi\Delta_{r,p_-}}    
    \NRadfordS(r,p_-)
    + e^{2i\pi\Delta_{r,p_-}}\hrhoS(r,p_-).
  \end{align*}
  {}From~\eqref{hrhoright} and~\eqref{hrhoright-bdr}, evidently,
  \begin{equation*}
    \ribbon\hrhoS(r,r')
    =e^{2i\pi\Delta_{r,r'}}\,\hrhoS(r,r')
  \end{equation*}
  for all the $\hrhoS$ elements in the lemma.
  
  Next, the $\modS$-action is calculated as
  \begin{multline*}
    \modS^{-1}(\NRadfordS(r,r'))
    ={}
    (-1)^{r'}\bigl(    
    \NDrinfeldNWSE(r,r')
    -(-1)^{p_-}
    \NDrinfeldNWSE(p_+\!-\!r,p_-\!-\!r')\bigr)\\
    ={}\ffrac{1}{\iiA}
    \sum_{s=1}^{p_+ - 1} (-1)^{r'(s + p_+) + p_- r}
    (\qp^{p_- r s}\!-\!\qp^{-p_- r s})
    \hrhoS(s,p_-)\\*
    {}+\ffrac{1}{\iiA}
    \!\!\sum_{(s,s')\in\setR}\!\!\!  (-1)^{r s' + r' s}
    (\qp^{p_- r s}\!-\!\qp^{-p_- r s})
    (\qm^{p_+ r' s'}\!+\!\qm^{-p_+ r' s'})
    \hrhoS(s,s'),
  \end{multline*}
  where $(r,r')\in\setR$ and we used~\eqref{chiNWSE-decompose}, and,
  similarly,
  \begin{multline*}
    \modS^{-1}(\NRadfordS(r,p_-))
    ={}
    (-1)^{p_-}
    \NDrinfeldNWSE(r,p_-)\\
    ={}\ffrac{1}{\iiA}
    \sum_{s=1}^{p_+ - 1} (-1)^{p_-(s + p_+ + r)}\,
    (\qp^{p_- r s}\!-\!\qp^{-p_- r s})
    \hrhoS(s,p_-)\\*
    {}+\ffrac{1}{\iiA}
    \!\sum_{(s,s')\in\setR}\!\!\!
    2(-1)^{(p_+ + r) s' + p_- s}
    (\qp^{p_- r s}\!-\!\qp^{-p_- r s})
    \hrhoS(s,s').
  \end{multline*}
  Conversely, from~\eqref{chi-decompose} we have
  \begin{multline*}
    \modS^{-1}(\hrhoS(r,r'))=\\
    {}=(p_+\!-\!r)\bigl(\Drinfeld^{+}(r,r')
    +\Drinfeld^{-}(r,p_-\!-\!r')
    \bigr)-r\bigl(\Drinfeld^{-}(p_+\!-\!r,r')+
    \Drinfeld^{+}(p_+\!-\!r, p_-\!-\!r')\bigr)\\
    {}=- \ffrac{1}{\iiA} \biggl(\sum_{s=1}^{p_+ - 1}
    (-1)^{r p_- + r'(s + p_+)}
    (\qp^{p_- r s}\!-\!\qp^{-p_- r s})
    \NRadfordS(s,p_-)\\*
    {}+\sum_{(s,s')\in\setR}\!\!
    (-1)^{rs' + r's}
    (\qp^{p_- r s}\!-\!\qp^{-p_- r s})
    (\qm^{p_+ r' s'}\!+\qm^{-p_+ r' s'})
    \NRadfordS(s,s')\!\biggr),
  \end{multline*}
  which completes the proof.
\end{proof}

Totally similarly, we establish the following lemma.
\begin{lemma}\label{lemma:Rsh}
  For $(r,r')\in\setIRight$, the elements $\hrhoB(r,r')$
  \textup{(}see~\eqref{hrholeft} and~\eqref{hrholeft-bdr}\textup{)}
  and $\RadfordB(r,r')$ \textup{(}see~\eqref{hrhoLeft}
  and~\eqref{hrhoLeft2}\textup{)} span a
  $(p_+\,{+}\,1)(p_-\,{-}\,1)$-dimensional
  $\SLiiZ$-subrep\-resen\-tation $\RRLeft$ in~$\cZ$.  Moreover, its
  structure is
  \begin{equation*}
    \RRLeft=\oC^2\tensor\RLeft,
  \end{equation*}
  where $\oC^2$ is the standard two-dimensional $\SLiiZ$
  representation and $\RLeft$ is the representation
  in~\eqref{W-structure}.
\end{lemma}

To formulate our next (and final) lemma, we recall the notation
in~\eqref{Proj-notation}--\eqref{Proj-notation5}.
\begin{lemma}\label{lemma-proj}
  The elements $\RadfordProj(r,r')$, $(r,r')\in\setI$, span a
  $\half(p_+\,{+}\,1)(p_-\,{+}\,1)$-dimen\-sion\-al
  $\SLiiZ$-subrepresentation $\Rproj$ in~$\cZ$.
\end{lemma}
\begin{proof}
  Clearly, for $(r,r')\in\setI$, the ribbon element acts as
  \begin{equation*}
    \ribbon\RadfordProj(r,r')
    =e^{2i\pi\Delta_{r,r'}}\,\RadfordProj(r,r').
  \end{equation*}
  As regards the $\modS$ action, we recall~\eqref{S-pm}
  and~\eqref{chi-decompose}, which we rewrite, for $(r,r')\in\setR$
  and using notation~\eqref{Proj-notation}--\eqref{Proj-notation5}, as
  \begin{align*}
    \modS^{-1}(\RadfordProj(r,r'))={}& \Drinfeld^{+}(r,r') +
    \Drinfeld^{-}(p_+\!-\!r,r') + \Drinfeld^{-}(r,p_-\!-\!r') +
    \Drinfeld^{+}(p_+\!-\!r,p_-\!-\!r')\\
    ={}&\ffrac{1}{\iiA}
    \Bigl(\RadfordProj(p_+,p_-) +
    (-1)^{rp_- + r'p_+}\,\RadfordProj(0,p_-)\\*
    &{}+\sum_{s=1}^{p_+ - 1} (-1)^{r p_- + r'(p_+ - s)}
    (\qp^{p_{-} r s}\!+ \qp^{-p_{-} r s})
    \RadfordProj(s,p_-)\\*
    &{}+\sum_{s'=1}^{p_- - 1} (-1)^{r' p_+ + r(p_- - s')}
    (\qm^{p_{+} r' s'}\!+ \qp^{-p_{+} r' s'})
    \RadfordProj(p_+,s')\\*
    &{}+\!\!\sum_{(s,s')\in\setR}\!\!
    (-1)^{r s' + r' s}(\qp^{p_{-} r s}\!+ \qp^{-p_{-} r s})
    (\qm^{p_{+} r' s'}\!+ \qp^{-p_{+} r' s'})
    \RadfordProj(s,s')\!\Bigr),
  \end{align*}
  and similarly for $\modS^{-1}(\RadfordProj(r,p_-))$,
  $\modS^{-1}(\RadfordProj(p_+,r'))$, and
  $\modS^{-1}(\Radford^{\pm}(p_+,p_-))$.  The lemma is proved.
\end{proof}

This completes the proof of Claim~\ref{claim:decomp}
in~\bref{thm:q-structure}.

\subsubsection{}Claim~\ref{claim:grring} is now obvious if we recall
from~\eqref{chi-decompose} that the Drinfeld-map image of the
Grothendieck ring in $\cZ$ is spanned by
\begin{alignat*}{4}
  &\RadfordUpUp(r,r'),&\quad (r,r')&\in\setR,&
  &\RadfordProj(r,r'),& (r,r')&\in\setI,\\*
  &\RadfordS(r,r'),& (r,r')&\in\setILeft,&\qquad
  &\RadfordB(r,r'),&\quad (r,r')&\in\setIRight.
\end{alignat*}

\subsection{Comparison with the $\SLiiZ$-representation on generalized
  characters of the $W$-algebra}Claim~\ref{claim:equiv}
in~\bref{thm:q-structure} is established by direct comparison with
formulas in~\cite[Lemmas~4.3 and 4.4]{[FGST3]}.  First,
with~\bref{lemma:R-min}, we define
\begin{equation*}
  \chi_{r,r'}=\NDrinfeldNWSE(r, r') - \NDrinfeldNESW(r, r')
  =\modS\,\NvvarphiJoin(r,r'),
  \quad (r,r')\in\setR.
\end{equation*}
Then \textit{the $\chi_{r,r'}$ transform under $\modS$ and $\modT$
  in~\eqref{theST} as the functions $\chi_{r,r'}(\tau)$
  in~\cite{[FGST3]} transform under $\tau\mapsto-\frac{1}{\tau}$ and
  $\tau\mapsto\tau+1$.}  (In particular, $\modT\chi_{r,r'} =
\AT\,\modS(\ribbon\,\NvvarphiJoin(r,r')) =\AT\,
e^{2i\pi\Delta_{r,r'}}\chi_{r,r'}$.)

Next, with~\bref{lemma:triplet}, we define
\begin{align*}
  \rho_{r,r'}&=\modS\,\vvarphi(r,r')
  =\ffrac{-2\sqrt{2}}{\sqrt{p_+ p_-\mathstrut}}
  \!\sum_{(s,s')\in\setR}\!\!\!
  (-1)^{(r+1)s' + (r'+1)s}
  \sin\!\ffrac{\pi p_- r s}{p_+}\,\sin\!\ffrac{\pi p_+ r' s'}{p_-}\,
  \NRadfordUpUp(s,s'),\\
  \psi_{r,r'}&=\modS\,\Nppsi(r,r')
  =\ffrac{2\sqrt{2}}{\sqrt{p_+ p_-\mathstrut}}
    \!\sum_{(s,s')\in\setR}\!\!\!(-1)^{r s' + r' s}
    \sin\!\ffrac{\pi p_- r s}{p_+}\,
    \sin\!\ffrac{\pi p_+ r' s'}{p_-}\,\Nppsi(s,s'),\\
  \varphi_{r,r'}&=(-1)^{r+r'}\NDrinfeldDnDn(r,r'),
\end{align*}
where $(r,r')\in\setR$.  Then \textit{$\rho_{r,r'}$, $\psi_{r,r'}$,
  and $\varphi_{r,r'}$ transform under $\modS$ and $\modT$
  in~\eqref{theST} as the respective functions $\rho_{r,r'}(\tau)$,
  $\psi_{r,r'}(\tau)$, and $\varphi_{r,r'}(\tau)$ in~\cite{[FGST3]}
  transform under $\tau\mapsto-\frac{1}{\tau}$ and
  $\tau\mapsto\tau+1$.}  (In particular, $\modT\psi_{r,r'}=\AT\,
e^{2i\pi\Delta_{r,r'}}(\psi_{r,r'}+2\rho_{r,r'})$ and
$\modT\varphi_{r,r'}=\AT\,e^{2i\pi\Delta_{r,r'}}(\varphi_{r,r'} +
\psi_{r,r'}+ \rho_{r,r'})$.)

Similarly, with~\bref{lemma:Lsh}, we set
\begin{alignat*}{2}
  \varphi^{\boxslash}_{r,r'}&=
  -(-1)^{r'}\NDrinfeldNWSE(r,r')
  +(-1)^{p_- + r'}\NDrinfeldNWSE(p_+\!-\!r, p_-\!-\!r'),&
  \quad &(r,r')\in\setR,\\
  \varphi^{\boxslash}_{r,0}&=
  (-1)^{p_-}\NDrinfeldNWSE(p_+\!-\!r,p_-),&\ &1\leq r\leq p_+\!-\!1,
  \\
  \intertext{and}
  \rho^{\boxslash}_{r,r'}&=-\modS\,\hrhoS(r,r')
  =r\bigl(\Drinfeld^{-}(p_+\!-\!r,r')
  + \Drinfeld^{+}(p_+\!-\!r, p_-\!-\!r')\bigr)\\*
  &\qquad\qquad\qquad\qquad{}-(p_+\!-\!r)
  \bigl(\Drinfeld^{+}(r,r')+\Drinfeld^{-}(r,p_-\!-\!r')\bigr),
  &\quad &(r,r')\in\setR,
  \\
  \rho^{\boxslash}_{r,0}&=\modS\,\hrhoS(p_+\!-\!r,p_-)
  =r\,\Drinfeld^{+}(p_+\!-\!r,p_-) - (p_+\!-\!r)\Drinfeld^{-}(r,p_-).
  &\ &1\leq r\leq p_+\!-\!1.
\end{alignat*}
Then \textit{the $\rho_{r,r'}^{\boxslash}$ and
  $\varphi_{r,r'}^{\boxslash}$ thus defined transform under $\modS$
  and $\modT$ in~\eqref{theST} as the respective functions
  $\rho_{r,r'}^{\boxslash}(\tau)$ and
  $\varphi_{r,r'}^{\boxslash}(\tau)$ in~\cite{[FGST3]} transform under
  $\tau\mapsto-\frac{1}{\tau}$ and $\tau\mapsto\tau+1$.}

We give some more details here, mainly because of a minor difference
between the conventions in this paper and in~\cite{[FGST]}.  First,
from the formulas in the proof of~\bref{lemma:Lsh}, it follows that
\begin{multline*}
  \modS\varphi^{\boxslash}_{r,r'}=
  \ffrac{\sqrt{2}\,i}{\sqrt{\mathstrut p_+p_-}}\Bigl( \sum_{s=1}^{p_+
    - 1}(-1)^{sr'}\sin\!\ffrac{\pi p_-rs}{p_+}\,
  \rho^{\boxslash}_{s,0}
  \\
  {}+\sum_{(s,s')\in\setii}\!\!  2 (-1)^{rs'+sr'}\sin\!\ffrac{\pi
    p_-rs}{p_+}\, \cos\ffrac{\pi
    p_+r's'}{p_-}\,\rho^{\boxslash}_{s,s'} \Bigr),
\end{multline*}
for all the $\varphi^{\boxslash}$ defined above, which coincides with the
$\modS$-transformation formula for $\varphi^{\boxslash}_{r,r'}(\tau)$
found in~\cite{[FGST3]}.  Next, from the last formula in the proof
of~\bref{lemma:Lsh} and formulas
\eqref{hrhoRight}--\eqref{hrhoRight2}, it follows that
\begin{multline*}
  \modS\rho^{\boxslash}_{r,r'}= \ffrac{-\sqrt{2}\,i}{\sqrt{\mathstrut
      p_+p_-}}\Bigl( \sum_{s=1}^{p_+ - 1}(-1)^{sr'}\sin\!\ffrac{\pi
    p_-rs}{p_+}\, \varphi^{\boxslash}_{s,0}
  \\
  {}+\sum_{(s,s')\in\setii}\!\!  2 (-1)^{rs'+sr'}\sin\!\ffrac{\pi
    p_-rs}{p_+}\, \cos\ffrac{\pi
    p_+r's'}{p_-}\,\varphi^{\boxslash}_{s,s'} \Bigr)
\end{multline*}
for $(r,r')\in\setR$, which is also the formula for
$\rho^{\boxslash}_{r,r'}(-\frac{1}{\tau})$ in~\cite{[FGST3]}.
Finally, it follows from the definitions (see~\eqref{hrhoright-bdr})
that $\modS\rho^{\boxslash}_{r,0}=\hrhoS(p_+\,{-}\,r,p_-)
=r\modS\Drinfeld^+(p_+\,{-}\,r,p_-)
-(p_+\,{-}\,r)\modS\Drinfeld^-(r,p_-)$.  Using~\eqref{chi-decompose}
now gives
\begin{multline*}
  r\Drinfeld^+(p_+\!-\!r,p_-)-(p_+\!-\!r)\Drinfeld^-(r,p_-)
  = \ffrac{\sqrt{2}\,i}{\sqrt{\mathstrut p_+ p_-\!}}
  \Bigl(\sum_{s=1}^{p_+-1}(-1)^{(r-1)p_-} \sin\!\ffrac{\pi p_- r
    s}{p_+}\,
  \NRadfordNESW(s,p_-)\\*
  +\!\!\sum_{(s,s')\in\setR}\!\!2(-1)^{r s'} \sin\!\ffrac{\pi p_- r
    s}{p_+}\, \NRadfordS(s,s') \Bigr),
\end{multline*}
which readily shows that $\modS\rho^{\boxslash}_{r,0}$ is given by the
above formula for $\modS\rho^{\boxslash}_{r,r'}$ at $r'=0$ (as
in~\cite{[FGST3]}).

It next follows from the ribbon element action in~\bref{lemma:Lsh}
that
\begin{equation*}
  \modT\varphi^{\boxslash}_{r,r'}
  =\AT\,e^{2i\pi\Delta_{r,r'}}(\varphi^{\boxslash}_{r,r'}
  +\rho^{\boxslash}_{r,r'})
\end{equation*}
for all $\varphi^{\boxslash}_{r,r'}$ defined above (for
$\varphi^{\boxslash}_{r,0}$, the simple derivation involves the
identity $\Delta_{r,0}=\Delta_{p_+-r,p_-}$,
see~\eqref{Delta-rs}
).  Equally easily, \bref{lemma:Lsh} also implies that
\begin{equation*}
  \modT\rho^{\boxslash}_{r,r'}
  =\AT\,e^{2i\pi\Delta_{r,r'}}\rho^{\boxslash}_{r,r'}
\end{equation*}
for all $\rho^{\boxslash}_{r,r'}$ defined above.

The correspondence involving the central elements in~\bref{lemma:Rsh}
is totally similar.  Finally, for the central elements
in~\bref{lemma-proj}, the correspondence is also established
straightforwardly, and we omit the details.

\subsection{Factorization of the $\SLiiZ$ action} We next consider
factorization Claim~\ref{claim:factorization}
in~\bref{thm:q-structure}.  As in a simpler case studied
in~\cite{[FGST]}, this factorization is related to the factorization
of the ribbon element in~\eqref{ribbon-factors}.  We first note that
with the $M$-matrix expressed as
$M=(\ribbon\tensor\ribbon)\Delta(\ribbon)^{-1}$, it readily
follows~\cite{[Lyu]} that
\begin{equation}
  \modS(\ribbon)=\ribbon^{-1}.
\end{equation}

We now define a map $\xxi:\Ch\to\cZ$, which is an isomorphism of
vector spaces and intertwines the coadjoint and adjoint actions
of~$\XXX$,~as
\begin{equation}\label{xxi-def}
  \xxi(\beta)=(\beta\tensor\id)M^*,
\end{equation}
where
\begin{equation}\label{M-star-def}
  M^* =(\ribbon^*\tensor\ribbon^*)
  \Delta(\modS(\ribbon^*)).
\end{equation}
{}From the $\modS$-transformation formulas in~\bref{sec:modular-eval},
we immediately obtain
\begin{equation*}
  \modS\bigl(\ribbon^*\bigr)
  =  \coint
  + \ffrac{1}{p_+ p_-}\,\NRadfordUpUp(1, 1)
  + \ffrac{1}{p_+}\,\NRadfordNESW(1, 1)
  + \ffrac{1}{p_-}\,\NRadfordNWSE(1, 1).
\end{equation*}

The maps $\cZ\to\cZ$ given by
\begin{equation}
  \modS^*=\Radford\circ\xxi^{-1}
\end{equation}
and
\begin{equation}
  \bar\modS=\xxi\circ\Drinfeld^{-1}
\end{equation}
then provide a factorization of $\modS:\cZ\to\cZ$,
\begin{equation}
  \modS=\modS^*\bar\modS.
\end{equation}
The $\SLiiZ$-representation on the center is thus factored into two
representations, generated by $(\modS^*,\modT^*)$ and
$(\bar\modS,\bar\modT)$.

The representation generated by $(\modS^*,\modT^*)$ can be further
factored in accordance with the factorization of $\ribbon^*$ as
\begin{equation}\label{ribbon-factors3}
  \ribbon^* = \ribbonplus\,\ribbonminus,
\end{equation}
where 
\begin{align*}
  \ribbonpm
  ={}& 
  \begin{cases}
    \one
    + \ffrac{1}{p_+}\,
    \NDrinfeldNWSE(1, 1)\\[6pt]
    \one
    + \ffrac{1}{p_-}\,
    \NDrinfeldNESW(1, 1)
  \end{cases}\\
  ={}& \one
  - \sum_{m=1}^{p_\pm - 1}\sum_{a=m-1}^{p_\pm - 1}\!
  (-1)^{m}
  \ffrac{(\q_{\pm}^{p_\mp}\!- \q_{\pm}^{-p_\mp})^{2m-1}}{p_\pm [m]_\pm}\,
  \q_{\pm}^{p_\mp m(m - 1 - 2a) - p_\mp (2 + 2a)}
  \\*
  &\qquad\qquad\qquad{}\times
  {\qbin{a}{m\!-\!1}}_\pm^2
  e_{\pm}^{m} f_{\pm}^m K^{-p_\mp (m - 2 - 2a)}
\end{align*}
and we calculated $\DrinfeldNWSE(1,1)$ and $\DrinfeldNESW(1,1)$
from~\eqref{chiNWSE-eval} and~\eqref{chiNESW-eval} using~\eqref{C-low}
and \eqref{C-top}.  This leads to a further factorization of the
$\SLiiZ$-action on the center, constructed similarly
to~\eqref{xxi-def}--\eqref{M-star-def}.

\section{Conclusions}
A remarkable correspondence between modular group actions in
nonsemisimple (logarithmic) conformal field theory models and on
quantum groups deserves further study.  Another quite interesting
problem is to extend the category of finite-dimensional
$\XXX$-rerpesentations such that it becomes \textit{equivalent} to the
category of the $\WWW$ algebra representations realized in logarithmic
$(p_+,p_-)$ conformal field theory models.

\subsubsection*{Acknowledgments} We are grateful to A.~Belavin and
S.~Parkhomenko for the useful discussions.  We thank G.~Mutafian for
help at the early stages of this work.  This paper was supported in
part by the RFBR Grant 04-01-00303 and the RFBR--JSPS Grant
05-01-02934YaF\_a.  AMS was also supported in part by the RFBR Grant
05-01-00996. Part of the paper was written when the authors were
visiting Kyoto University, and we are grateful to T.~Miwa for
hospitality.

\appendix

\section{Some indecomposable
  $\smash{\XXX}$-modules}\label{app:priojective}

\subsection{$\smash{\modPr^{\pm,\pm}_{r,r'}}$ modules}
\label{app:(1,p)-proj-mod-base} As a preparation for the description
of projective modules, we introduce modules $\modPr^{\alpha,+}_{r,r'}$
and $\modPr^{\alpha,-}_{r,r'}$, where $\alpha=\pm$, with the
respective structure of subquotients described as
\begin{equation}\label{schem-(1,p)-proj}
  \includegraphics[bb=1.7in 9in 7.6in 10.36in, clip]{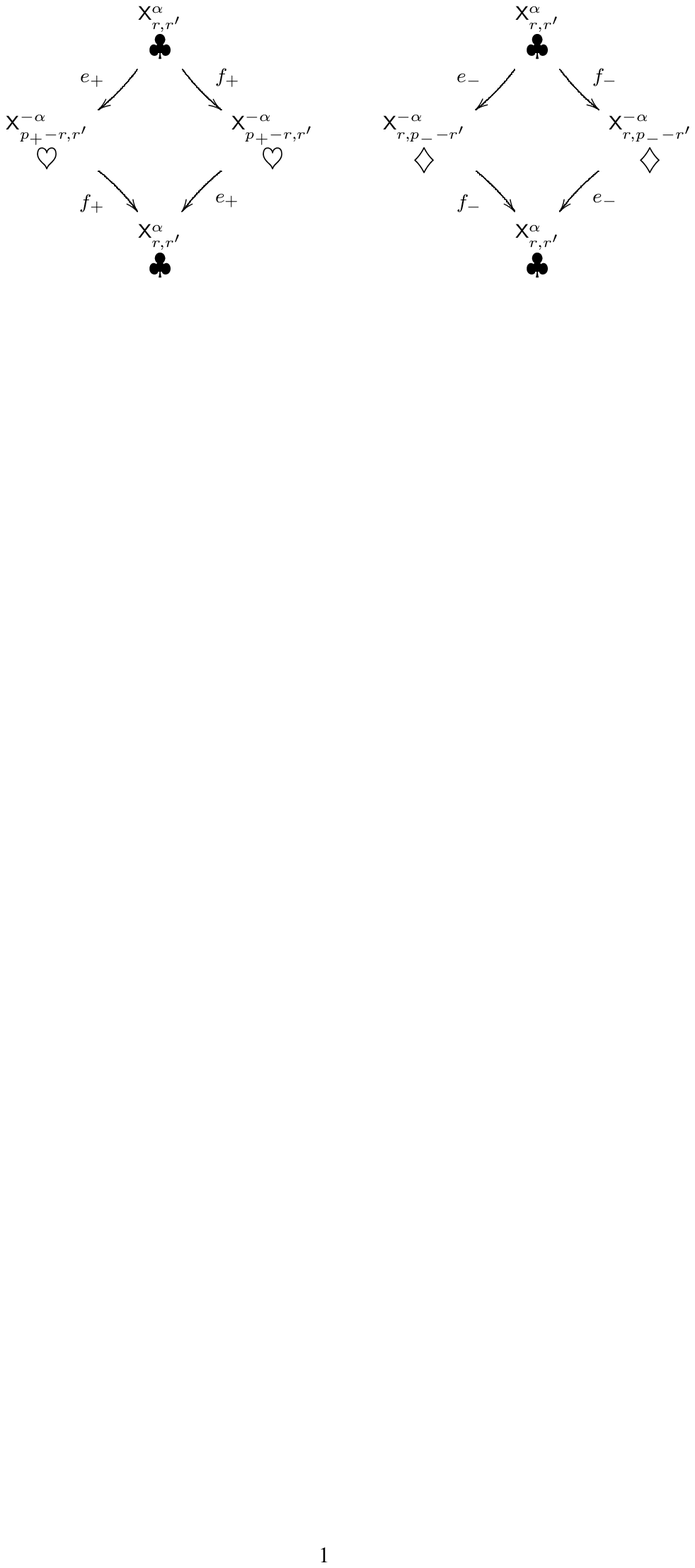}
\end{equation}
For $1\leq r\leq p_+-1$ and $1\leq r'\leq
p_-$,  the $\modPr^{\pm,+}_{r,r'}$ module has the basis
\begin{equation}\label{app:(1,p)-projective-basis}
  \{\prUp{t}_{n,n'},\prDown{t}_{n,n'}\}_{\substack{0\le n\le r-1\\0\le
      n'\le r'-1}}
  \cup\{\prLeft{t}_{k,k'},\prRight{t}_{k,k'}\}_{\substack{0\le k\le
      p_+-r-1\\0\le k'\le r'-1}},
\end{equation}
where $\prUp{t}_{n,n'}$ is the basis corresponding to the top
$\XX^{\pm}_{r,r'}$ module in~\eqref{schem-(1,p)-proj},
$\prDown{t}_{n,n'}$ to the bottom $\XX^{\pm}_{r,r'}$,
$\prRight{t}_{k,k'}$ to the right $\XX^{\mp}_{p_+-r,r'}$, and
$\prLeft{t}_{k,k'}$ to the left $\XX^{\mp}_{p_+-r,r'}$ module.  We fix
the arbitrariness involved in choosing $\prRight{t}$, $\prLeft{t}$,
and $\prUp{t}$ by specifying the $\XXX$-action
in~$\modPr^{\alpha,+}_{r,r'}$ as
\begin{multline*}
  K\prUp{t}_{n,n'}
  =\alpha \qp^{r-1-2n}\,\qm^{r'-1-2n'}\prUp{t}_{n,n'},\quad
  K\prDown{t}_{n,n'}
  =\alpha \qp^{r-1-2n}\,\qm^{r'-1-2n'}\prDown{t}_{n,n'},\\
  0\le n\le r\!-\!1,\quad
  0\le n'\le r'\!-\!1,
\end{multline*}
\begin{multline*}
  K\prRight{t}_{k,k'}=-\alpha \qp^{p_+-r-1-2k}\,\qm^{r'-1-2k'}
  \prRight{t}_{k,k'},\quad
  K\prLeft{t}_{k,k'}
  =-\alpha \qp^{p_+-r-1-2k}\,\qm^{r'-1-2k'}\prLeft{t}_{k,k'},\\
  0\le k\le p_+\!-\!r\!-\!1,\quad
  0\le k'\le r'\!-\!1,
\end{multline*}
\begin{align}\label{app:(1,p)-proj-XXX-action-plus}
  \ep\prUp{t}_{n,n'}&=
  \begin{cases}
    (\alpha)^{p_-}(-1)^{r'-1}
    \qint{n}_+\qint{r\!-\!n}_+\prUp{t}_{n-1,n'}+\prDown{t}_{n-1,n'},
    &1\le n\le r-1,\\
    \prLeft{t}_{p_+\!-\!r\!-\!1,n'}, &  n=0,\\
  \end{cases}
  \notag\\
  \ep\prRight{t}_{k,k'}&=
  \begin{cases}
    (-\alpha)^{p_-}(-1)^{r'-1}
    \qint{k}_+\qint{p_+\!-\!r\!-\!k}_+\prRight{t}_{k-1,k'},
    &1\le k\le p_+\!-\!r\!-\!1,\\
    \prDown{t}_{r-1,k'}, & k=0,\\
  \end{cases}
  \notag\\
  \ep\prLeft{t}_{k,k'}
  &=(-\alpha)^{p_-}(-1)^{r'-1}\qint{k}_+\qint{p_+\!-\!r\!-\!k}_+
  \prLeft{t}_{k-1,k'},
  \;
  \begin{array}[t]{ll}
    0\le k\le p_+\!-\!r\!-\!1\\
    (\text{with}\;\prLeft{t}_{-1,k'}\equiv0),
  \end{array}
  \\
  \ep\prDown{t}_{n,n'}
  &=(\alpha)^{p_-}(-1)^{r'-1}\qint{n}_+\qint{r\!-\!n}_+
  \prDown{t}_{n-1,n'},
  \quad 0\le n\le r\!-\!1
  \quad(\text{with}\;\prDown{t}_{-1,n'}\equiv0),\notag\\
  \fp\prUp{t}_{n,n'}&=
  \begin{cases}
    \prUp{t}_{n+1,n'}, & 0\le n\le r\!-\!2,\\
    \prRight{t}_{0,n'}, & n=r\!-\!1,
  \end{cases}\quad
  \fp\prLeft{t}_{k,k'}=
  \begin{cases}
    \prLeft{t}_{k+1,k'}, & 0\le k\le p_+\!-\!r\!-\!2,\\
    \prDown{t}_{0,k'}, & k=p_+\!-\!r\!-\!1,\\
  \end{cases}
  \notag\\
  \fp\prRight{t}_{k,k'}&=\prRight{t}_{k+1,k'}, \quad
  0 \le k\le p_+\!-\!r\!-\!1
  \quad(\text{with}\;\prRight{t}_{p_+-r,k'}\equiv0),\notag\\
  \fp\prDown{t}_{n,n'}&=\prDown{t}_{n+1,n'}, \quad 0\le
    n\le r\!-\!1\quad(\text{with}\;\prDown{t}_{r,n'}\equiv0),\notag
\end{align}
where $0\le n',k'\le r'\!-\!1$, and
\begin{align*}
  \emi\prUp{t}_{n,n'}
  &=(\alpha)^{p_+}(-1)^{r-1}\qint{n'}_-
  \qint{r'\!-\!n'}_-\prUp{t}_{n,n'-1}
  \quad(\text{with}\;\prUp{t}_{n,-1}\equiv0),&\\
  \emi\prRight{t}_{k,k'}
  &=(-\alpha)^{p_+}(-1)^{p_+-r-1}
  \qint{k'}_-\qint{r'\!-\!k'}_-\prRight{t}_{k,k'-1}
  \quad(\text{with}\;\prRight{t}_{k,-1}\equiv0),&\\
  \emi\prLeft{t}_{k,k'}
  &=(-\alpha)^{p_+}(-1)^{p_+-r-1}\qint{k'}_-
  \qint{r'\!-\!k'}_-\prLeft{t}_{k,k'-1}
  \quad(\text{with}\;\prLeft{t}_{k,-1}\equiv0),&\\
  \emi\prDown{t}_{n,n'}
  &=(\alpha)^{p_+}(-1)^{r-1}\qint{n'}_-
  \qint{r'\!-\!n'}_-\prDown{t}_{n,n'-1}
  \quad(\text{with}\;\prDown{t}_{n,-1}\equiv0),&\\
  \fm\prUp{t}_{n,n'}
  &=\prUp{t}_{n,n'+1}\quad(\text{with}\;\prUp{t}_{n,r'}\equiv0),\quad
  \fm\prRight{t}_{k,k'}
  =\prRight{t}_{k,k'+1}\quad(\text{with}\;\prRight{t}_{k,r'}\equiv0),\\
  \fm\prLeft{t}_{k,k'}
  &=\prLeft{t}_{k,k'+1}\quad(\text{with}\;\prLeft{t}_{k,r'}\equiv0),
  \quad
  \fm\prDown{t}_{n,n'}
  =\prDown{t}_{n,n'+1}\quad(\text{with}\;\prDown{t}_{n,r'}\equiv0),
\end{align*}
where $0\le n\le r\!-\!1$, $0\le k\le p_+\!-\!r\!-\!1$, and $0\le
n',k'\le r'\!-\!1$.

For $1\leq r\leq p_+$ and $1\leq r'\leq p_-\!-\!1$, we similarly define
the $\modPr^{\pm,-}_{r,r'}$ module.

\subsection{Projective modules}\label{app:proj-mod} We describe the
projective modules $\PP^{\pm}_{r,r'}$ for $1\leq r\leq p_+$ and $1\leq
r'\leq p_-$.
First, for $1\leq r\leq p_+\!-\!1$ and $1\leq r'\leq p_-\!-\!1$, the
projective module $\PP^{\pm}_{r,r'}$ is a second extension of the
Verma module $\VV^{\pm}_{r,r'}$,
\begin{equation}\label{schem-proj-Verma}
  \xymatrix@=12pt{
    &&\VV^{\pm}_{r,r'}
    \ar@/^/[dl]_{\emi} \ar@/_/[dr]^{\ep}&\\
    &\VV^{\mp}_{r,p_--r'}
    \ar@/^/[dr]_{\ep}&
    &\VV^{\mp}_{p_+-r,r'}
    \ar@/_/[dl]^{\emi}\\
    &&\VV^{\pm}_{p_+-r,p_--r'}&
  }
\end{equation}
and is therefore a ``deck'' of sixteen subquotients shown in
Fig.~\ref{app:proj-pic-plus}.  We note that $\dim\PP^{\pm}_{r,r'}=4
p_+ p_-$ for $1\leq r\leq p_+\!-\!1$ and $1\leq r'\leq p_-\!-\!1$.
Second, $\PP^{\pm}_{r,p_-}=\modPr^{\pm,+}_{r,p_-}$ for $1\leq r\leq
p_+\!-\!1$ and $\PP^{\pm}_{p_+,r'}=\modPr^{\pm,-}_{p_+,r'}$ for $1\leq
r'\leq p_-\!-\!1$, where the modules $\modPr^{\pm,+}_{r,p_-}$ and
$\modPr^{\pm,-}_{p_+,r'}$ are defined above
in~\bref{app:(1,p)-proj-mod-base}. Third, $\PP^{\pm}_{p_+,
  p_-}=\XX^{\pm}_{p_+, p_-}$ are irreducible.

The module $\PP^{\pm}_{r,r'}$ has the basis
\begin{equation}\label{app:projective-basis}
  \prBullet{t}_{n,n'},\prBullet{b}_{n,n'},
  \prBullet{r}_{k,k'},\prBullet{l}_{k,k'},
\end{equation}
where $\sbullet\in\{\suparrow,\sdownarrow,\srightarrow,
\sleftarrow\}$, $0\le n'\le r'-1$, $0\le k'\le p_--r'-1$, $0\le n,k\le
r-1$ if $\sbullet\in\{\suparrow, \sdownarrow\}$,\\ $0\le n,k\le
p_+-r-1$ if $\sbullet\in\{\srightarrow, \sleftarrow\}$, and
$\prBullet{t}_{n,n'}$ is the basis corresponding to the top module
$\modPr^{\pm,+}_{r,r'}$, $\prBullet{b}_{n,n'}$ to the bottom
$\modPr^{\pm,+}_{r,r'}$, $\prBullet{r}_{k,k'}$ to the right
$\modPr^{\mp,+}_{r,p_--r'}$, and $\prBullet{l}_{k,k'}$ to the left
module $\modPr^{\mp,+}_{r,p_--r'}$ in the diagram
\begin{equation}\label{schem-proj-(1,p)-proj-plus}
  \xymatrix@=12pt{
    &&\modPr^{\pm,+}_{r,r'}
    \ar@/^/[dl]_{\emi} \ar@/_/[dr]^{\fm}&\\
    &{\modPr^{\mp,+}_{r,p_--r'}} \ar@/^/[dr]_{\fm}&
    &{\modPr^{\mp,+}_{r,p_--r'}}
    \ar@/_/[dl]^{\emi}\\
    &&\modPr^{\pm,+}_{r,r'}&
  }
\end{equation}

The action of the generators $K$, $\ep$, and $\fp$ on
$\prBullet{t}_{n,n'}$ is defined in
\eqref{app:(1,p)-proj-XXX-action-plus} and the action of $\emi$ and
$\fm$ is given by
\begin{align*}
  \emi\prBullet{t}_{n,n'}&=
  \begin{cases}
    (\pm1)^{p_+}(-1)^{r-1}\qint{n'}_-
    \qint{r'\!-\!n'}_-\prBullet{t}_{n,n'-1}+\prBullet{b}_{n,n'-1},
    &1 \le n'\le r'\!-\!1,\\
    \prBullet{l}_{n,p_--r'-1}, & n'=0,\\
  \end{cases}
  \\
  \fm\prBullet{t}_{n,n'}&=
  \begin{cases}
    \prBullet{t}_{n,n'+1}, &0\le n'\le r'\!-\!2,\\
    \prBullet{r}_{n,0}, & n'=r'\!-\!1.
  \end{cases}
  \\
  \intertext{where $\sbullet\in\{\suparrow,\sdownarrow\}$ and $0\le
    n\le r-1$,}
  \emi\prBullet{t}_{n,n'}&=
  \begin{cases}
    (\mp1)^{p_+}(-1)^{p_+-r-1}\qint{n'}_-
    \qint{r'\!-\!n'}_-\prBullet{t}_{n,n'-1}+\prBullet{b}_{n,n'-1},
    &1\le n'\le r'\!-\!1,\\
    \prBullet{l}_{n,p_--r'-1}, & n'=0,\\
  \end{cases}
  \\
  \fm\prBullet{t}_{n,n'}&=
  \begin{cases}
    \prBullet{t}_{n,n'+1}, &0\le n'\le r'{-}2,\\
    \prBullet{r}_{n,0}, &  n'=r'{-}1.
  \end{cases}
\end{align*}
where $\sbullet\in\{\srightarrow,\sleftarrow\}$ and $0\le n\le
p_+-r-1$.

The action of the generators $K$, $\ep$, $\fp$, and $\fm$ on
$\prBullet{r}_{n,n'}$ coincides with the action on
$\modPr^{\mp,+}_{r,p_--r'}$ (see \bref{app:(1,p)-proj-mod-base}). The
action of $\emi$ is given by
\begin{align*}
  \emi\prBullet{r}_{k,k'}&=
  \begin{cases}
    (\mp1)^{p_+}(-1)^{r-1}
    \qint{k'}_-\qint{p_-\!-\!r'\!-\!k'}_-\prBullet{r}_{k,k'-1},
    &1\le k'\le p_-\!-\!r'\!-\!1,\\
    \prBullet{b}_{k,r'-1}, & k'=0,\\
  \end{cases}
  \\ \intertext{where $\sbullet\in\{\suparrow,\sdownarrow\}$ and $0\le
    k\le r-1$,} \emi\prBullet{r}_{k,k'}&=
  \begin{cases}
    (\pm1)^{p_+}(-1)^{p_+-r-1}
    \qint{k'}_-\qint{p_-\!-\!r'\!-\!k'}_-\prBullet{r}_{k,k'-1},
    &1\le k'\le p_-\!-\!r'\!-\!1,\\
    \prBullet{b}_{k,r'-1}, & k'=0,\\
  \end{cases}
\end{align*}
where $\sbullet\in\{\srightarrow,\sleftarrow\}$ and $0\le k\le
p_+-r-1$.

The action of the generators $K$, $\ep$, $\fp$, and $\emi$ on
$\prBullet{l}_{n,n'}$ coincides with the action on
$\modPr^{\mp,+}_{r,p_--r'}$ (see \bref{app:(1,p)-proj-mod-base}). The
action of $\fm$ is given by
\begin{align*}
  \fm\prBullet{l}_{k,k'}&=
  \begin{cases}
    \prBullet{l}_{k,k'+1}, & 0\le k'\le p_-\!-\!r'\!-\!2,\\
    \prBullet{b}_{k,0}, &  k'=p_-\!-\!r'\!-\!1,\\
  \end{cases}
  \\ \intertext{where $\sbullet\in\{\suparrow,\sdownarrow\}$ and $0\le
    k\le r-1$,} \fm\prBullet{l}_{k,k'}&=
  \begin{cases}
    \prBullet{l}_{k,k'+1}, & 0\le k'\le p_-\!-\!r'\!-\!2,\\
    \prBullet{b}_{k,0}, & k'=p_-\!-\!r'\!-\!1,\\
  \end{cases}
\end{align*}
where $\sbullet\in\{\srightarrow,\sleftarrow\}$ and $0\le k\le
p_+\!-\!r\!-\!1$.

The $\XXX$-action on $\prBullet{b}_{n,n'}$ is defined in
\bref{app:(1,p)-proj-mod-base}.

\section{The center of $\smash{\XXX}$}\label{app:the-center}
Here, we find the center of $\XXX$ using the isomorphism between the
center and the algebra of \textit{bimodule} endomorphisms of the
regular representation. In~\bref{sec-RegDecompose}, we study
decomposition of the regular representation as a $\XXX$-bimodule.
In~\bref{sec:bimod-center}, we describe the algebra of bimodule
endomorphisms of the regular representation.
In~\bref{sec:basis-explicit}, the canonical central elements are also
constructed explicitly in terms of the algebra generators.
In~\bref{app:calculations}, we collect several calculations needed in
the proofs in this paper.

\subsection{Decomposition of the regular representation}
\label{sec-RegDecompose}
We study the decomposition of the regular representation as a
$\XXX$-bimodule. For this, we closely follow the strategy proposed
in~\cite{[FGST]} for $(p,1)$ models. The starting point is to
introduce indecomposable $\XXX$-bimodules $\modQ^{\pm,\pm}_{r,r'}$
defined as follows.  For $1\leq r\leq p_+-1$ and $1\leq r'\leq p_-$,
$\modQ^{\pm,+}_{r,r'}$ is composed of eight subquotients with the left
$\XXX$-action 
represented as
\begin{equation*}
  \includegraphics[bb=1.5in 9in 7in 10.4in, clip]{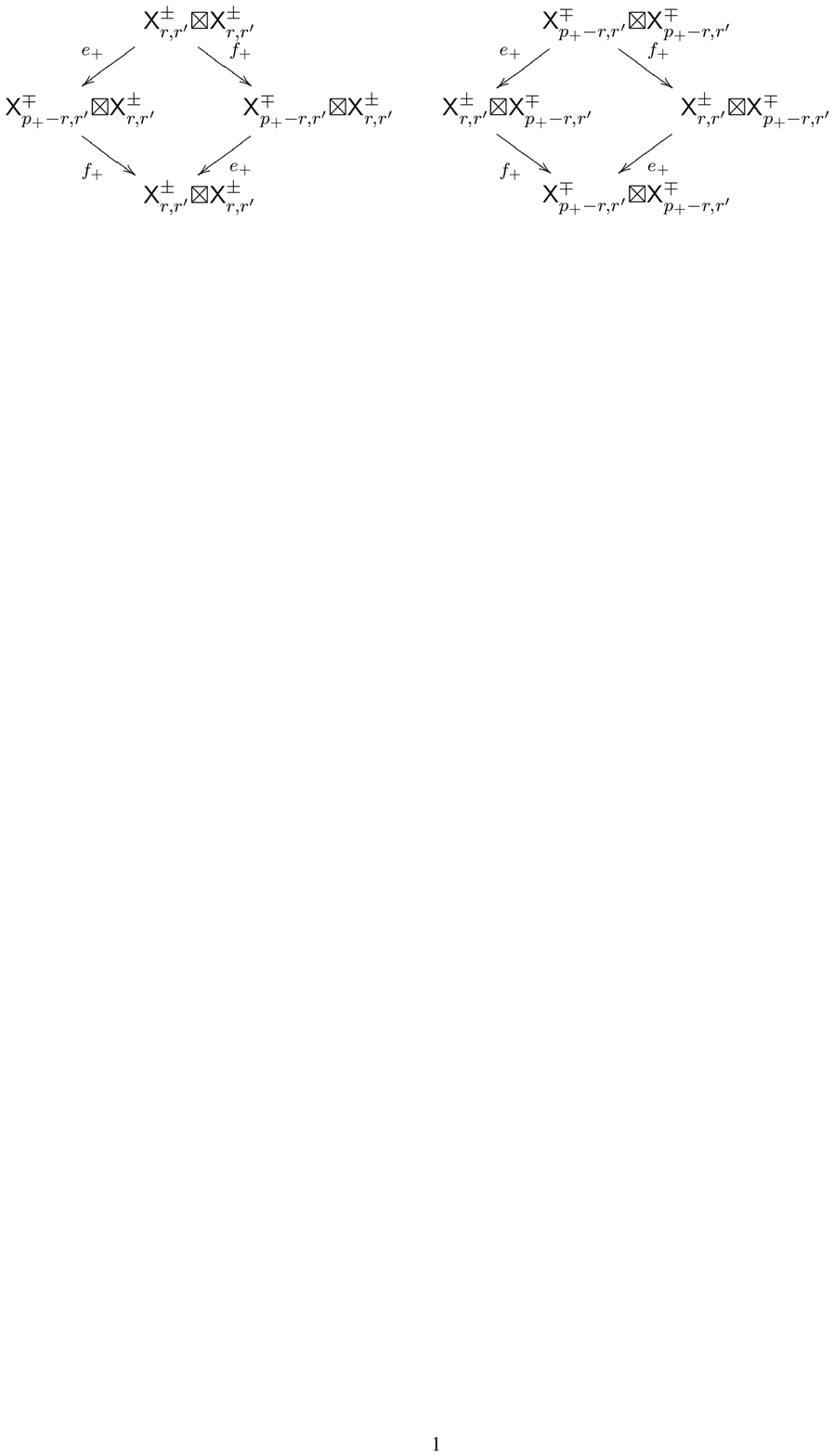}
\end{equation*}
and the right action as
\begin{equation*}
  \includegraphics[bb=1.5in 9in 7in 10.4in, clip]{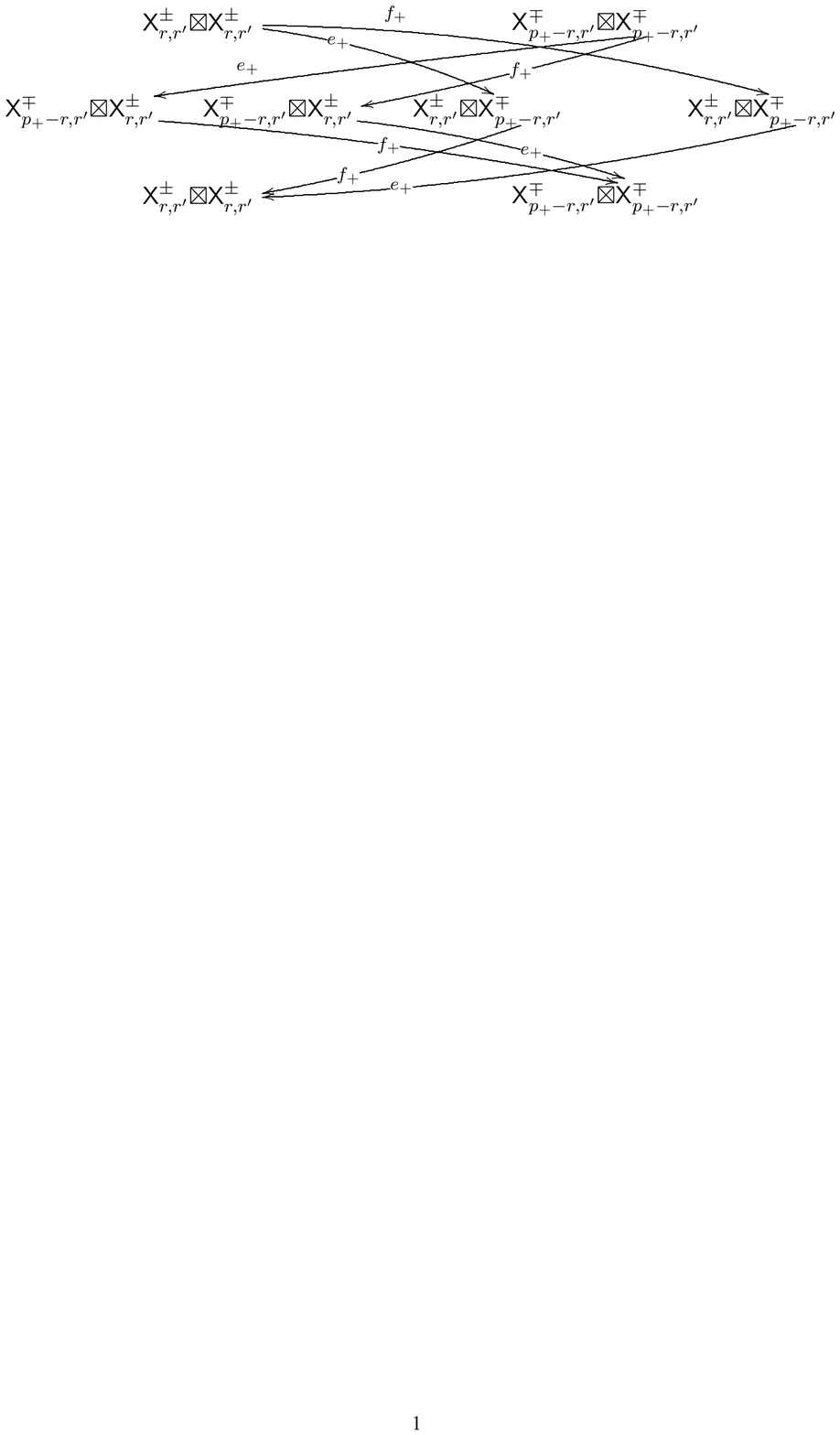}
  \end{equation*}
  In these diagrams, as we have noted, it is understood that the
  $\XXX$-action on each subquotient is changed compared\pagebreak[3]
  with the action on the corresponding irreducible representation in
  agreement with the arrows connecting a given subquotient with
  others.  For example, $\ep$ maps the highest-weight vector of
  $\XX^{+}_{r,r'}$ into the lowest-weight vector
  in~$\XX^{-}_{p_+-r,r'}$.

For $1\leq r\leq p_+$ and $1\leq r'\leq p_--1$, the bimodules
$\modQ^{\pm,-}_{r,r'}$ are defined in the same way with the
substitution $\ep\mapsto\emi$, $\fp\mapsto\fm$, and
$\XX^{\mp}_{p_+-r,r'}\mapsto\XX^{\mp}_{r,p_--r'}$ in the above
diagrams.

\begin{prop}
  As a $\XXX$-bi\-module, the regular representation decomposes as
  \begin{equation*}
    \mathsf{Reg}=\smash[b]{\bigoplus_{(r,r')\in\setR}\modQ(r,r')\oplus
    \bigoplus_{r=1}^{p_+-1}\modQ(r,p_-)\oplus
    \bigoplus_{r'=1}^{p_--1}\modQ(p_+,r')
    \oplus\modQ(p_+,p_-)\oplus\modQ(0,p_-)},
  \end{equation*}
  where
  \begin{enumerate}
  \item $\modQ(p_+,p_-)=\XX^{+}_{p_+,p_-}\boxtimes\XX^{+}_{p_+,p_-}$
    and $\modQ(0,p_-)=\XX^{-}_{p_+,p_-}\boxtimes\XX^{-}_{p_+,p_-}$ are
    simple bimodules;
    
  \item $\modQ(r,p_-)=\modQ^{+,+}_{r,p_-}$ and
    $\modQ(p_+,r')=\modQ^{+,-}_{p_+,r'}$ for $1\leq r\leq p_+-1$ and
    $1\leq r'\leq p_--1$;
    
  \item the bimodules $\modQ(r,r')$, $1\leq r\leq p_+-1$ and $1\leq
    r'\leq p_--1$, are indecomposable and admit two length-$3$
    filtrations:
    \begin{equation}\label{filtr-Q+(r,r')}
      0\subset\repR^+_2(r,r')\subset\repR^+(r,r')\subset\modQ(r,r'),
    \end{equation}
    where the structure of subquotients is given by
    \begin{align*}
      \modQ(r,r')/\repR^+(r,r')&= \modQ^{+,+}_{r,r'}\oplus
      \modQ^{-,+}_{r,p_--r'},\\
      \repR^+(r,r')/\repR^+_2(r,r') &= 2\modQ^{-,+}_{r,p_--r'} \oplus
      2\modQ^{+,+}_{r,r'},\\
      \repR^+_2(r,r') &= \modQ^{+,+}_{r,r'}\oplus
      \modQ^{-,+}_{r,p_--r'},
    \end{align*}
    and
    \begin{equation}\label{filtr-Q-(r,r')}
      0\subset\repR^-_2(r,r')\subset\repR^-(r,r')\subset\modQ(r,r'),
    \end{equation}
    where the structure of subquotients is given by
    \begin{align*}
      \modQ(r,r')/\repR^-(r,r')&= \modQ^{+,-}_{r,r'}\oplus
      \modQ^{-,-}_{p_+-r,r'},\\
      \repR^-(r,r')/\repR^-_2(r,r') &= 2\modQ^{-,-}_{p_+-r,r'} \oplus
      2\modQ^{+,-}_{r,r'},\\
      \repR^-_2(r,r') &= \modQ^{+,-}_{r,r'}\oplus
      \modQ^{-,-}_{p_+-r,r'}.
    \end{align*}
  \end{enumerate}
\end{prop}
We do not repeat the standard steps leading to this statement, the
derivation is totally similar to the one in~\cite{[FGST]}.

\subsection{Bimodule homomorphisms and the
  center}\label{sec:bimod-center} Here, we study the algebra of
bimodule endomorphisms of the regular representation $\mathsf{Reg}$.
Simultaneously, we derive the structure of the $\XXX$ center because
bimodule endomorphisms of the regular representation are in a
one-to-one correspondence with elements in the center.

\subsubsection{}\label{bimodule--center} Clearly, the algebra of
bimodule endomorphisms of $\mathsf{Reg}$ decomposes, in agreement
with~\eqref{center-decomposition}, as
\begin{multline*}
  \End\bigl(\mathsf{Reg}\bigr)
  =\bigoplus_{(r,r')\in\setR}
  \End\bigl(\modQ(r,r')\bigr)\oplus
  \bigoplus_{r=1}^{p_+-1}\End\bigl(\modQ(r,p_-)\bigr)\oplus\\[-4pt]
  \bigoplus_{r'=1}^{p_--1}\End\bigl(\modQ(p_+,r')\bigr)\oplus
  \End\bigl(\modQ(p_+,p_-)\bigr)\oplus\End\bigl(\modQ(0,p_-)\bigr).
\end{multline*}

For each $\modQ(r,r')$, there is a bimodule endomorphism
$\idem(r,r'):\mathsf{Reg}\to\mathsf{Reg}$ that acts as identity on
$\modQ(r,r')$ and is zero on~$\modQ(s,s')$ unless $s= r$ and $s'= r'$.
These endomorphisms give rise to $\half(p_+\,{+}\,1)(p_-\,{+}\,1)$
primitive idempotents (also denoted by $\idem(r,r')$) in the center
of~$\XXX$.

Next, for each $\modQ(r,r')$ with $(r,r')\in\setR$,
\begin{enumerate}

\item there is a homomorphism $\vNE(r,r'):\modQ(r,r')\to\modQ(r,r')$
  (defined up to a nonzero factor) whose image is
  $\modQ^{+,+}_{r,r'}$; in other words, $\vNE(r,r')$ sends the
  quotient $\modQ^{+,+}_{r,r'}$ (see~\eqref{filtr-Q+(r,r')}) into the
  subbimodule $\modQ^{+,+}_{r,r'}$ at the bottom of $\modQ(r,r')$ and
  is zero on $\modQ(s,s')$ unless $s= r$ and $s'= r'$.

\item Similarly, there is a central element associated with the
  homomorphism $\vNW(r,r'):\modQ(r,r')\to\modQ(r,r')$ with the image
  $\modQ^{+,-}_{r,r'}$, i.e., sending the quotient
  $\modQ^{+,-}_{r,r'}$ (see~\eqref{filtr-Q-(r,r')}) into the
  subbimodule $\modQ^{+,-}_{r,r'}$ (and zero on~$\modQ(s,s')$ unless
  $s= r$ and $s'= r'$).

\item A central element is associated with the homomorphism
  $\vSW(r,r'):\modQ(r,r')\to\modQ(r,r')$ with the image
  $\modQ^{-,+}_{r,p_--r'}$, i.e., sending the quotient
  $\modQ^{-,+}_{r,p_--r'}$ (see~\eqref{filtr-Q+(r,r')}) into the
  subbimodule $\modQ^{-,+}_{r,p_--r'}$ (and zero on~$\modQ(s,s')$
  unless $s= r$ and $s'= r'$).
  
\item Similarly, there is a central element associated with the
  homomorphism $\vSE(r,r'):\modQ(r,r')\to\modQ(r,r')$ with the image
  $\modQ^{-,-}_{p_+-r,r'}$, i.e., sending the quotient
  $\modQ^{-,-}_{p_+-r,r'}$ (see~\eqref{filtr-Q-(r,r')}) into the
  subbimodule $\modQ^{-,-}_{p_+-r,r'}$ (and zero on~$\modQ(s,s')$
  unless $s= r$ and $s'= r'$).
  
\item There is also a homomorphism
  $\wUp(r,r'):\modQ(r,r')\to\modQ(r,r')$ (defined up to a nonzero
  factor) whose image is $\XX^{+}_{r,r'}\boxtimes\XX^{+}_{r,r'}$; in
  other words, $\wUp(r,r')$ sends the quotient
  $\XX^{+}_{r,r'}\boxtimes\XX^{+}_{r,r'}$ into the subbimodule
  $\XX^{+}_{r,r'}\boxtimes\XX^{+}_{r,r'}$ at the bottom of
  $\modQ(r,r')$ and is zero on~$\modQ(s,s')$ unless $s= r$ and $s'=
  r'$.
  
\item Similarly, there is a homomorphism
  $\wRight(r,r'):\modQ(r,r')\to\modQ(r,r')$ whose image is
  $\XX^{-}_{p_+-r,r'}\boxtimes\XX^{-}_{p_+-r,r'}$; in other words,
  $\wRight(r,r')$ sends the quotient
  $\XX^{-}_{p_+-r,r'}\boxtimes\XX^{-}_{p_+-r,r'}$ into the subbimodule
  $\XX^{-}_{p_+-r,r'}\boxtimes\XX^{-}_{p_+-r,r'}$ at the bottom of
  $\modQ(r,r')$ and is zero on~$\modQ(s,s')$ unless $s= r$ and $s'=
  r'$.
  
\item There is also a homomorphism
  $\wLeft(r,r'):\modQ(r,r')\to\modQ(r,r')$ whose image is
  $\XX^{-}_{r,p_--r'}\boxtimes\XX^{-}_{r,p_--r'}$; in other words,
  $\wLeft(r,r')$ sends the quotient
  $\XX^{-}_{r,p_--r'}\boxtimes\XX^{-}_{r,p_--r'}$ into the subbimodule
  $\XX^{-}_{r,p_--r'}\boxtimes\XX^{-}_{r,p_--r'}$ at the bottom of
  $\modQ(r,r')$ and is zero on~$\modQ(s,s')$ unless $s= r$ and $s'=
  r'$.
  
\item Similarly, there is also a homomorphism
  $\wDown(r,r'):\modQ(r,r')\to\modQ(r,r')$ whose image is
  $\XX^{+}_{p_+-r,p_--r'}\boxtimes\XX^{+}_{p_+-r,p_--r'}$; in other
  words, $\wDown(r,r')$ sends the quotient
  $\XX^{+}_{p_+-r,p_--r'}\boxtimes\XX^{+}_{p_+-r,p_--r'}$ into the
  subbimodule $\XX^{+}_{p_+-r,p_--r'}\boxtimes\XX^{+}_{p_+-r,p_--r'}$
  at the bottom of $\modQ(r,r')$ and is zero on~$\modQ(s,s')$ unless
  $s= r$ and $s'= r'$.
\end{enumerate}

These maps give $4(p_+\!-\!1)(p_-\!-\!1)$ nilpotent elements
$\wLeft(r,r')$, $\wUp(r,r')$, $\wRight(r,r')$, $\wDown(r,r')$,
$\vNW(r,r')$, $\vNE(r,r')$, $\vSE(r,r')$, $\vSW(r,r')$. {}From the
structure of the maps, it is obvious that relations~\eqref{rad-times}
are satisfied and that all other compositions are zero.

Next, for each $\modQ(r,p_-)$ with $1\leq r\leq p_+-1$, there is a
homomorphism $\vUp(r,p_-):\modQ(r,p_-)\to\modQ(r,p_-)$ with the image
$\XX^{+}_{r,p_-}\boxtimes\XX^{+}_{r,p_-}$; in other words,
$\vUp(r,p_-)$ sends the quotient
$\XX^{+}_{r,p_-}\boxtimes\XX^{+}_{r,p_-}$ into the subbimodule
$\XX^{+}_{r,p_-}\boxtimes\XX^{+}_{r,p_-}$ at the bottom of
$\modQ(r,p_-)$ and is zero on~$\modQ(r',p_-)$ with~\hbox{$r'\neq r$}.
Similarly, for each $r=1,\dots,p_+-1$, there is a central element
associated with the homomorphism
$\vRight(r,p_-):\modQ(r,p_-)\to\modQ(r,p_-)$ with the image
$\XX^{-}_{p_+-r,p_-}\boxtimes\XX^{-}_{p_+-r,p_-}$, i.e., sending the
quotient $\XX^{-}_{p_+-r,p_-}\boxtimes\XX^{-}_{p_+-r,p_-}$ into the
subbimodule $\XX^{-}_{p_+-r,p_-}\boxtimes\XX^{-}_{p_+-r,p_-}$ (and
zero on~$\modQ(r',p_-)$ with $r'\neq r$).  Next, for each
$\modQ(p_+,s)$ with $1\leq s\leq p_--1$, there are similar
homomorphisms $\vUp(p_+,s)$ and
$\vLeft(p_+,s):\modQ(p_+,s)\to\modQ(p_+,s)$.

This gives $2(p_+\!-\!1) + 2(p_-\!-\!1)$ more elements $\vUp(r,p_-)$,
$\vRight(r,p_-)$, $\vUp(p_+,s)$, $\vLeft(p_+,s)$, which are
obviously in the radical of the center.

The above is just restated in~\bref{prop-center}.

\subsubsection{Canonical decomposition of central
  elements}\label{sec:coeffs} Any central element $A$ can be
decomposed in the above central elements as
\begin{multline}\label{decomp-general}
  A=\sum_{(r,r')\in\setI} a_{r,r'} \idem(r,r')
  + \sum_{\substack{(r,r')\in\setR \\ \sbullet \in
      \{\skewarrows\}}}\cRadbullet{r,r'} \Radbullet(r,r')\\
  {}+\sum_{\substack{(r,r')\in\setR \\ \sbullet \in
      \{\rightarrows\}}}\cNilpbullet{r,r'}\Nilpbullet(r,r')
  + \sum_{\substack{r=1\\\sbullet \in
      \{\suparrow,\srightarrow\}}}^{p_+-1} \cRadbullet{r}
  \Radbullet(r,p_-) + \sum_{\substack{r'=1\\\sbullet \in
      \{\suparrow,\sleftarrow\}}}^{p_--1} \cNilpbullet{r'}
  \Radbullet(p_+,r').
\end{multline}
It immediately follows from~\bref{bimodule--center} that the
coefficient $a_{r,r'}$ is the eigenvalue of $A$ in the irreducible
representation $\XX^{+}_{r,r'}$.  To determine the
$\cRadbullet{r,r'}$, $\cNilpbullet{r,r'}$, $\cRadbullet{r}$, and
$\cNilpbullet{r'}$ coefficients similarly, we fix the normalizations
such that in terms of the respective bases described
in~\bref{app:proj-mod}, $\Radbullet(r,r')$, $\Nilpbullet(r,r')$,
$\Radbullet(r,p_-)$, and $\Radbullet(p_+,r')$ act as
\begin{align*}
  \vNE(r,r')\,\prBullet{t}_{n,n'}&=\prBullet{b}_{n,n'},
  \ \ \sbullet \in \{\rightarrows\},
  \quad\text{in $\PP^{+}_{r,r'}$ and $\PP^{-}_{p_+-r,r'}$},\\
  \intertext{with $\vNE(r,r')$ being identically zero on any module
  other than $\PP^{+}_{r,r'}$ or $\PP^{-}_{p_+-r,r'}$; similarly,}
  \vNW(r,r')\,\prUp{\mathsf{x}}_{n,n'}&=\prDown{\mathsf{x}}_{n,n'},
  \ \ \mathsf{x}\in\{\mathsf{t},\mathsf{r},\mathsf{l},\mathsf{b}\},
  \quad\text{in $\PP^{+}_{r,r'}$ and $\PP^{-}_{r,p_--r'}$},\\
  \vSW(r,r')\,\prBullet{t}_{n,n'}&=\prBullet{b}_{n,n'},
  \ \ \sbullet \in \{\rightarrows\},
  \quad\text{in $\PP^{-}_{r,p_--r'}$ and $\PP^{+}_{p_+-r,p_--r'}$},\\
  \vSE(r,r')\,\prUp{\mathsf{x}}_{n,n'}&=\prDown{\mathsf{x}}_{n,n'},
  \ \ \mathsf{x}\in\{\mathsf{t},\mathsf{r},\mathsf{l},\mathsf{b}\},
  \quad\text{in $\PP^{-}_{p_+-r,r'}$ and $\PP^{+}_{p_+-r,p_--r'}$},\\
  \wUp(r,r')\,\prUp{t}_{n,n'}&=\prDown{b}_{n,n'}
  \quad\text{in $\PP^{+}_{r,r'}$},\\
  \wRight(r,r')\,\prUp{t}_{n,n'}&=\prDown{b}_{n,n'}
  \quad\text{in $\PP^{-}_{p_+-r,r'}$},\\
  \wLeft(r,r')\,\prUp{t}_{n,n'}&=\prDown{b}_{n,n'}
  \quad\text{in $\PP^{-}_{r,p_--r'}$},\\
  \wDown(r,r')\,\prUp{t}_{n,n'}&=\prDown{b}_{n,n'}
  \quad\text{in $\PP^{+}_{p_+-r,p_--r'}$},\\
  \vUp(r,p_-)\,\prUp{t}_{n,n'}&=\prDown{t}_{n,n'}
  \quad\text{in $\PP^{+}_{r,p_-}$},\\
  \vRight(r,p_-)\,\prUp{t}_{n,n'}&=\prDown{t}_{n,n'}
  \quad\text{in $\PP^{-}_{p_+-r,p_-}$},\\
  \vUp(p_+,r')\,\prUp{t}_{n,n'}&=\prDown{t}_{n,n'}
  \quad\text{in $\PP^{+}_{p_+,r'}$},\\
  \vLeft(p_+,r')\,\prUp{t}_{n,n'}&=\prDown{t}_{n,n'}
  \quad\text{in $\PP^{-}_{p_+,p_--r'}$}.
\end{align*}
Then the coefficients in~\eqref{decomp-general} are determined from
the relation
\begin{alignat*}{3}
  A\prBullet{t}_{n,n'}&=\cRadNE{r,r'}\prBullet{b}_{n,n'}&
  \quad&\text{in}\quad
  \PP^{+}_{r,r'},&
  \quad&\sbullet\in\{\rightarrows\},\\
  A\prUp{\mathsf{x}}_{n,n'}&=\cRadNW{r,r'}\prDown{\mathsf{x}}_{n,n'}&
  \quad&\text{in}\quad \PP^{-}_{r,p_--r'},&
  \quad&\mathsf{x}\in\{\mathsf{t},\mathsf{r},\mathsf{l},\mathsf{b}\},\\
  A\prBullet{t}_{n,n'}&=\cRadSW{r,r'}\prBullet{b}_{n,n'}&
  \quad&\text{in}\quad\PP^{+}_{p_+-r,p_--r'},&
  \quad&\sbullet\in\{\rightarrows\},\\
  A\prUp{\mathsf{x}}_{n,n'}&=\cRadSE{r,r'}\prDown{\mathsf{x}}_{n,n'}&
 \quad&\text{in}\quad \PP^{-}_{p_+-r,r'},&
  \quad&\mathsf{x}\in\{\mathsf{t},\mathsf{r},\mathsf{l},\mathsf{b}\},\\
  A\prUp{t}_{n,n'}&=\cNilpUp{r,r'}\;\prDown{b}_{n,n'}&
  \quad&\text{in}\quad
  \PP^{+}_{r,r'},\\
  A\prUp{t}_{n,n'}&=\cNilpRight{r,r'}\;\prDown{b}_{n,n'}&
  \quad&\text{in}\quad
  \PP^{-}_{p_+-r,r'},\\
  A\prUp{t}_{n,n'}&=\cNilpLeft{r,r'}\;\prDown{b}_{n,n'}&
  \quad&\text{in}\quad
  \PP^{-}_{r,p_--r'},\\
  A\prUp{t}_{n,n'}&=\cNilpDown{r,r'}\;\prDown{b}_{n,n'}&
  \quad&\text{in}\quad
  \PP^{+}_{p_+-r,p_--r'},\\
  A\prUp{t}_{n,n'}&=\cRadUp{r}\; \prDown{t}_{n,n'}&
  \quad&\text{in}\quad
  \PP^{+}_{r,p_-},\\
  A\prUp{t}_{n,n'}&=\cRadRight{r}\; \prDown{t}_{n,n'}&
  \quad&\text{in}\quad
  \PP^{-}_{p_+-r,p_-},\\
  A\prUp{t}_{n,n'}&=\cNilpUp{r'}\; \prDown{t}_{n,n'}&
  \quad&\text{in}\quad
  \PP^{-}_{p_+,r'},\\
  A\prUp{t}_{n,n'}&=\cNilpLeft{r'}\; \prDown{t}_{n,n'}&
  \quad&\text{in}\quad
  \PP^{-}_{p_+,p_--r'}.
\end{alignat*}
These formulas are used many times in explicit calculations in the
text.

\subsection{Explicit construction of the canonical basis in the
  center}\label{sec:basis-explicit} To explicitly construct the
canonical central elements in~\bref{prop-center} in terms of the
$\XXX$ generators, we proceed similarly to~\cite{[FGST],[Kerler]}.

Let $(r,r')\in\setI$ (see~\eqref{eq:setI}).  We set
\begin{align*}
  \beta_+(r, r') &= (-1)^{r'}(\qp^{p_- r}\!+ \qp^{-p_- r}),\\
  \beta_-(r, r') &= (-1)^{r}(\qm^{p_+ r'}\!+ \qm^{-p_+ r'}).
\end{align*}
These are roots of the respective polynomials $\psi_+$ and $\psi_-$
(see~\bref{sec:casimirs}).  Moreover, \textit{as $(r,r')$ range over
  the set $\setI$, all roots of $\psi_+$ and $\psi_-$ occur among the
  values taken by $\beta_+(r,r')$ and $\beta_-(r,r')$ respectively}
(not necessarily once).  We then define
\begin{equation*}
  \psi_{\pm,r, r'}(x) = 
  \begin{cases}
    \ffrac{\psi_\pm(x)}{x - \beta_\pm(r, r')},
    & \beta_\pm(r, r') = 2\quad\text{or}\quad\beta_\pm(r, r') =
    -2,\\[6pt]
    \ffrac{\psi_\pm(x)}{(x - \beta_\pm(r, r'))^2}&\text{otherwise}
  \end{cases}
\end{equation*}
and set
\begin{equation*}
  \nilp_{\pm}(r,r')=
  \begin{cases}
    0, &\beta_\pm(r, r') = 2\quad\text{or}\quad\beta_\pm(r, r') = -2,\\
    (\cas_\pm - \beta_\pm(r, r'))\psi_{\pm,r, r'}(\cas_\pm)
    &\text{otherwise}
  \end{cases}
\end{equation*}
and
\begin{gather*}
  \idem_\pm(r, r') =
    \ffrac{1}{\psi_{\pm,r, r'}(\beta_\pm(r,r'))}
    \Bigl(\psi_{\pm,r, r'}(\cas_\pm)-                  
    \ffrac{\psi'_{\pm,r, r'}(\beta_\pm(r, r'))}{
      \psi_{\pm,r, r'}(\beta_\pm(r,r'))}\,
    \nilp_\pm(r, r')\!\Bigr),
\end{gather*}
where $\psi'_{\pm,r,r'}(x)=\dd\psi_{\pm,r,r'}(x)/\dd x$. 

Then
\begin{equation}\label{the-idems}
  \Idem(r,r')=\idem_+(r,r')\idem_-(r,r'),\quad
  (r,r')\in\setI,
\end{equation}
are the $\half(p_+\,{+}\,1)(p_-\,{+}\,1)$ primitive idempotents in the
center.  

The constants $\psi_{\pm,r, r'}(\beta_\pm(r,r'))$ involved in the
normalization factors are easily found explicitly.  For this, we
recall that the Chebyshev polynomials are eigenfunctions of a
second-order differential operator,
\begin{gather}\label{cheb-eigen}
  \Bigl(\!(x^2-4)\ffrac{\partial^2}{\partial x^2}
  +3x\ffrac{\partial}{\partial x}+1\Bigr)
  \cheb_s(x)=s^2\cheb_s(x).
\end{gather}
Whenever $\beta_\pm(r,r')\neq2$ and $\beta_\pm(r,r')\neq-2$, i.e.,
$a=\beta_\pm(r,r')$ is a multiplicity-$2$ root of $\psi_\pm(x)$
in~\bref{sec:Grring}, we use~\eqref{cheb-eigen} and~\eqref{cheb-rec}
to find
\begin{equation*}
  (a^2-4)\psi_{\pm,r, r'}(a)
  = 4p_\pm^2 + 2p_\pm a \cheb_{2p_\pm}(a)
\end{equation*}
and hence
\begin{equation}
  \psi_{+,r, r'}(\beta_+(r,r'))=
  \ffrac{4p_+^2}{(\qp^{p_- r}\!- \qp^{-p_- r})^2},
  \quad
  \psi_{-,r, r'}(\beta_-(r,r'))=
  \ffrac{4p_-^2}{(\qm^{p_+ r'}\!- \qm^{-p_+ r'})^2}.
\end{equation}
Also, if $\beta_\pm(r,r')=2$, then
\begin{equation*}
  \psi_{\pm,r,r'}(\beta_\pm(r,r'))=
  \ffrac{4p_\pm^2}{3} + \ffrac{4p_\pm}{3} \cheb_{2p_\pm}(2)
  =\ffrac{4p_\pm^2}{3} + \ffrac{8p_\pm^2}{3}
  =4p_\pm^2,
\end{equation*}
and if $\beta_\pm(r,r')=-2$, then
\begin{equation*}
  \psi_{\pm,r,r'}(\beta_\pm(r,r'))=
  -\ffrac{4p_\pm^2}{3} + \ffrac{4p_\pm}{3} \cheb_{2p_\pm}(-2)
  =-4p_\pm^2.
\end{equation*}

To construct the canonical nilpotent elements, we next introduce the
operators\footnote{The notation is correlated with the top, left,
  right, and bottom modules in~\eqref{schem-Verma}; these
  $\boldsymbol{\pi}$ operators project on the weights (i.e., the
  eigenvalues of $K$) occurring in the respective irreducible
  modules.}
\begin{align*}
  \piUp(r, r') &= 
  \ffrac{1}{2 p_+ p_-}
  \sum_{j=0}^{2 p_+ p_- - 1}
  \sum_{\substack{a=-r + 1\\\text{step}=2}}^{r - 1}
  \sum_{\substack{a'=-r' + 1\\\text{step}=2}}^{r' - 1}
  \qp^{-a j}\qm^{-a' j}K^j,
  \\
  \piLeft(r, r') &= \ffrac{1}{2 p_+ p_-}
  \sum_{j=0}^{2 p_+ p_- - 1}
  \sum_{\substack{a=-r + 1\\\text{step}=2}}^{r - 1}
  \sum_{\substack{a'=-p_- + r' + 1\\\text{step}=2}}^{p_- - r' - 1}
  (-1)^j \qp^{-a j}\qm^{-a' j}K^j,
  \\
  \piRight(r, r') &= 
  \ffrac{1}{2 p_+ p_-} \sum_{j=0}^{2 p_+ p_- - 1}
  \sum_{\substack{a=-p_+ + r + 1\\\text{step}=2}}^{p_+ - r - 1}
  \sum_{\substack{a'=-r' + 1\\\text{step}=2}}^{r' - 1}
  (-1)^j \qp^{-a j}\qm^{-a' j} K^j,
  \\
  \piDown(r, r') &= \ffrac{1}{2 p_+ p_-}
  \sum_{j=0}^{2 p_+ p_- - 1}
  \sum_{\substack{a=-p_+ + r + 1\\\text{step}=2}}^{p_+ - r - 1}
  \sum_{\substack{a'=-p_- + r' + 1\\\text{step}=2}}^{p_- - r' - 1}
  \qp^{-a j}\qm^{-a' j}K^j.
\end{align*}
These are pairwise orthogonal projection operators.  We also note that
\begin{equation*}
  \piUp(r, r')+\piLeft(r, r')+\piRight(r, r')+\piDown(r, r')
  =\half(\one + (-1)^{p_- (r - 1) + p_+ (r' - 1)}K^{p_+ p_-})
\end{equation*}
and
\begin{equation*}
  \piLeft(r, p_-)=\piRight(p_+, r')=
  \piDown(r,p_-)=\piDown(p_+,r')=0.
\end{equation*}

For $(r,r')\in\setR$ (see~\eqref{set-R}), we use the above projectors
to define the central elements
\begin{equation}\label{the-v}
  \begin{aligned}
    \vNE(r,r')&=
    \idem_+(r,r')\nilp_-(r,r')\bigl(\piUp(r,r')+\piRight(r,r')\bigr),
    \\
    \vSW(r,r')&=
    \idem_+(r,r')\nilp_-(r,r')\bigl(\piLeft(r,r')+\piDown(r,r')\bigr),\\
    \vNW(r,r')&=
    \nilp_+(r,r')\idem_-(r,r')\bigl(\piUp(r,r')+\piLeft(r,r')\bigr),
    \\
    \vSE(r,r')&=
    \nilp_+(r,r')\idem_-(r,r')\bigl(\piRight(r,r')+\piDown(r,r')\bigr).
  \end{aligned}
\end{equation}
Clearly, their nonzero products are those in~\eqref{rad-times}, with
\begin{align*}
  \wUp(r,r')&=\nilp_+(r,r')\nilp_-(r,r')\piUp(r,r'),\\
  \wRight(r,r')&=\nilp_+(r,r')\nilp_-(r,r')\piRight(r,r'),\\
  \wLeft(r,r')&=\nilp_+(r,r')\nilp_-(r,r')\piLeft(r,r'),\\
  \wDown(r,r')&=\nilp_+(r,r')\nilp_-(r,r')\piDown(r,r').
\end{align*}
The idempotent $\Idem(r,r')$ acts as identity on each of these $8$
elements.  This gives $(1+8)\cdot\half(p_+\!-\!1)(p_-\!-\!1)$ central
elements, labeled by $(r,r')\in\setR$.  In addition, for $1\leq r\leq
p_+\!-\!1$, there are the central radical elements
\begin{equation}\label{bdr1}
  \begin{aligned}
    \vUp(r,p_-)&=\nilp_+(r,p_-)\idem_-(r,p_-)\piUp(r,p_-),\\
    \vRight(r,p_-)&=\nilp_+(r,p_-)\idem_-(r,p_-)\piRight(r,p_-),
  \end{aligned}
\end{equation}
on which $\Idem(r,p_-)$ acts as identity (the total of $(1+2)\cdot(p_+
- 1)$ central elements), and for $1\leq r'\leq p_-\!-\!1$, there are the
central radical elements
\begin{equation}\label{bdr2}
  \begin{aligned}
    \vUp(p_+,r')&=
    \idem_+(p_+,r')\nilp_-(p_+,r')\piUp(p_+,r'),
    \\
    \vLeft(p_+,r')&=\idem_+(p_+,r')\nilp_-(p_+,r')\piLeft(p_+,r'),
  \end{aligned}
\end{equation}
on which $\Idem(p_+,r')$ acts as identity (the total of $(1+2)\cdot(p_-
- 1)$ central elements).

\subsection{Some $\smash{\XXX}$ calculations}\label{app:calculations}
\subsubsection{Proof of~\bref{prop:gammaUpUp}}
\label{app:proof-prop:gammaUpUp}
Here, we prove about a quarter of the relations
in~\bref{prop:gammaUpUp} by establishing a similar statement for
$\gamma(r,r')$ with $\alphaUp{\sdownarrow} = \alphaDown{\bullet} =
\betaDown{\bullet} =\betaUp{\bullet}=0$, in which case the trace
actually restricts to the submodule
$\modPr^{+,+}_{r,r'}\oplus\modPr^{-,+}_{p_+-r,r'}
\subset\mathbb{P}_{r,r'}$ (see the definition of the modules
$\modPr^{\pm,\pm}_{r,r'}$ in~\bref{app:(1,p)-proj-mod-base}), and
which eventually (after we ensure that $\alphaUp{\suparrow} =
\alphaUp{\srightarrow}$) gives $\gamma^{\NESW}(r,r')$
in~\eqref{gammaNESW-def}.

We thus temporarily write $\gamma^{\NESW}(r,r')$ for $\gamma(r,r')$ as
defined in~\eqref{last-azat-label
} in the case where
$\alphaUp{\sdownarrow} = \alphaDown{\bullet} = \betaDown{\bullet}
=\betaUp{\bullet}=0$.  It is a $q$-character if and only if
(see~\eqref{Ch-def})
\begin{equation*}
  0=\gamma^{\NESW}(x y)- \gamma^{\NESW}(S^2(y)x) \equiv
  \Tr_{\mathbb{P}_{r,r'}}(\balance^{-1}x[y,\sigma^{\NESW}]),
\end{equation*}
where, as we have just noted,
\begin{equation*}
  \Tr_{\mathbb{P}_{r,r'}}(\balance^{-1}x[y,\sigma^{\NESW}])
  =
  \Tr_{\modPr^{+,+}_{r,r'}\oplus\modPr^{-,+}_{p_+-r,r'}}
  (\balance^{-1}x[y,\sigma^{\NESW}])
\end{equation*}
It now follows from the formulas in~\bref{app:(1,p)-proj-mod-base}
that $[K,\sigma^{\NESW}(r,r')]=0$ and the other commutation relations
in $\modPr^{+,+}_{r,r'}$ are given by
\begin{align*}
  [\ep,\sigma^{\NESW}(r,r')]\nprDown{b}{\suparrow}_{n,n'} &=
  \begin{cases}
    \alphaUp{\suparrow}\nprDown{b}{\suparrow}_{n-1,n'},&
    1\leq n\leq r\!-\!1,\\
    \alphaUp{\suparrow}\nprLeft{b}{\suparrow}_{p_+-r-1,n'},& n=0,
  \end{cases}
  \\
  [\ep,\sigma^{\NESW}(r,r')]\nprUp{b}{\suparrow}_{n,n'} &=
  \begin{cases}
    - \alphaUp{\suparrow}\nprUp{b}{\suparrow}_{n-1,n'},& 1\leq n\leq
    r\!-\!1,\\ 0,&n=0,
  \end{cases}
  \\
  [\ep,\sigma^{\NESW}(r,r')]\nprRight{b}{\suparrow}_{k,k'} &=
  \begin{cases}
    0,& 1\leq k\leq p_+\!-\!r\!-\!1,\\
    - \alphaUp{\suparrow}\nprUp{b}{\suparrow}_{r-1,k'},& k=0,
  \end{cases}
  \\
  [\ep,\sigma^{\NESW}(r,r')]\nprLeft{b}{\suparrow}_{k,k'} &=0, \quad
  1\leq k\leq p_+\!-\!r\!-\!1.
\end{align*}
Writing the trace operation as
\begin{multline*}
  \Tr_{\modPr^{+,+}_{r,r'}}(x) =\sum_{n=1}^{r-1}\sum_{n'=0}^{r'-1}
  \bra{\nprUp{b}{\suparrow}_{n,n'}}x\ket{\nprUp{b}{\suparrow}_{n,n'}}
  +\sum_{n=1}^{p_+-r-1}\sum_{n'=0}^{r'-1}
  \bra{\nprLeft{b}{\suparrow}_{n,n'}}x
  \ket{\nprLeft{b}{\suparrow}_{n,n'}}\\*
  +\sum_{n=1}^{p_+-r-1}\sum_{n'=0}^{r'-1}
  \bra{\nprRight{b}{\suparrow}_{n,n'}}x
  \ket{\nprRight{b}{\suparrow}_{n,n'}}
  +\sum_{n=1}^{r-1}\sum_{n'=0}^{r'-1}
  \bra{\nprDown{b}{\suparrow}_{n,n'}}x
  \ket{\nprDown{b}{\suparrow}_{n,n'}},
\end{multline*}
we then have
\begin{equation*}
  \Tr_{\modPr^{+,+}_{r,r'}}(\balance^{-1}x[\ep,\sigma^{\NESW}])={}
  \alphaUp{\suparrow}\sum_{n'=0}^{r'-1}
  \bra{\nprDown{b}{\suparrow}_{0,n'}}\balance^{-1}x
  \ket{\nprLeft{b}{\suparrow}_{p_+-r-1,n'}}
  -\alphaUp{\suparrow}\sum_{n'=0}^{r'-1}
  \bra{\nprRight{b}{\suparrow}_{0,n'}}\balance^{-1}x
  \ket{\nprUp{b}{\suparrow}_{r-1,n'}}
\end{equation*}
and, similarly,
\begin{multline*}
  \Tr_{\modPr^{-,+}_{p_+-r,r'}}(\balance^{-1}x[\ep,\sigma^{\NESW}])
  ={}\\
  {}=\alphaUp{\srightarrow}\sum_{n'=0}^{r'-1}
  \bra{\nprDown{b}{\srightarrow}_{0,n'}}\balance^{-1}x
  \ket{\nprLeft{b}{\srightarrow}_{r-1,n'}}
  -\alphaUp{\srightarrow}\sum_{n'=0}^{r'-1}
  \bra{\nprRight{b}{\srightarrow}_{0,n'}}\balance^{-1}x
  \ket{\nprUp{b}{\srightarrow}_{p_+-r-1,n'}}.
\end{multline*}
Summing these two contributions of the modules $\modPr^{+,+}_{r,r'}$
and $\modPr^{-,+}_{p_+-r,r'}$ and using the obvious correspondence
between the bases in the respective composition factors in
$\modPr^{+,+}_{r,r'}$ and $\modPr^{-,+}_{p_+-r,r'}$, we immediately
obtain that the trace in~\eqref{q-char-cond} with $y=e_+$ vanishes if
and only if $\alphaUp{\suparrow} = \alphaUp{\srightarrow}$.

Similar calculations for the other generators ($\fp$, $\emi$, and
$\fm$) give no more conditions for $\gamma^{\NESW}(r,r')$ to be a
$q$-character.  This completes the proof in the simpler case of
$\gamma^{\NESW}$.  The case of $\gamma^{\UpUp}$ is equally uneventful
but somewhat lengthier.

\subsubsection{Casimir operators}\label{sec:casimirs}
In the calculations, we often use the well-known quantum-$s\ell(2)$
identity (see, e.g.,~\cite{[ChP]})
\begin{equation}\label{eq:EmFm-prod2}
  \fp^m \ep^m
  =(-1)^m \prod_{s=0}^{m-1}
  \ffrac{\cas_+ + \qp^{p_-(2s+1)}K_+ +\qp^{-p_-(2s+1)}K_+^{-1}}{
    (\qp^{p_-}\!-\!\qp^{-p_-})^2},\quad m<p_+,
\end{equation}
where the element
\begin{equation}\label{eq:casimir-plus}
  \cas_+=-\qp^{p_-}K_+^{-1} - \qp^{-p_-}K_+
  -(q_+^{p_-} - q_+^{-p_-})^2 \ep \fp
\end{equation}
is central in $\XXX$ and therefore acts by the same eigenvalue on all
vectors in each irreducible representation: using~\eqref{on-vectors},
we readily calculate
\begin{equation}\label{cas-value}
  \cas_+\bigr|_{\XX_{r,r'}^{\alpha}} =
  \alpha^{p_-}(-1)^{r'}(\qp^{p_- r}\!+ \qp^{-p_- r}).
\end{equation}

Similarly,
\begin{equation}\label{eq:EmFm-prod}
  \emi^{m'} \fm^{m'}=(-1)^{m'}\prod_{s'=0}^{m'-1}
  \ffrac{\cas_- + \qm^{-p_+(2s'+1)}K_-
    +\qm^{p_+(2s'+1)}K_-^{-1}}{
    (\qm^{p_+}\!-\!q_-^{-p_+})^2},
  \quad m'<p_-,
\end{equation}
where
\begin{equation}
  \label{eq:casimir-minus}
  \cas_-=-\qm^{p_+}K_-^{-1} - \qm^{-p_+}K_-
  -(q_-^{p_+} - q_-^{-p_+})^2 \emi \fm
\end{equation}
is central and
\begin{equation}\label{cas-value2}
  \cas_-\bigr|_{\XX_{r,r'}^{\alpha}} =
  \alpha^{p_+}(-1)^{r}(\qm^{p_+ r'}\!+ \qm^{-p_+ r'}).
\end{equation}

The elements $\cas_+$ satisfy the equations
\begin{equation}\label{psi-pm}
  \psi_+(\cas_+)=0,\qquad \psi_-(\cas_-)=0,
\end{equation}
where
\begin{equation*}
  \psi_{\pm}(x)
  =
  \prod_{r=0}^{p_\pm - 1} (x + \q_\pm^{p_\mp r}\!+ \q_\pm^{-p_\mp r})
  (x - \q_\pm^{p_\mp r} - \q_\pm^{-p_\mp r})
\end{equation*}
(another relation between $\cas_+$ and $\cas_-$, characteristic of the
$\XXX$ quantum group rather than of the quantum $s\ell(2)$, is given
in~\eqref{U=U}).  We note that $\psi_+(x)$ (and similarly $\psi_-(x)$)
has two roots of multiplicity~$1$, which are $\pm2$, and all the other
roots are of multiplicity~$2$.

In finding decompositions of central elements as explained
in~\bref{sec:coeffs}, we also need the action of $\cas_\pm$ on some
indecomposable modules.  For example, the action of $\cas_+$ on
$\modPr^{\alpha,+}_{r,r'}$ (see~\bref{app:(1,p)-proj-mod-base}) can be
written as
\begin{equation}\label{cas-plus-P}
  \cas_+\bigr|_{\modPr^{\alpha,+}_{r,r'}} =
  \alpha^{p_-}(-1)^{r'}(\qp^{p_- r}\!+ \qp^{-p_- r})
  - (\qp^{p_-}\!- \qp^{-p_-})^2
  \sum_{n,n'}\prDown{t}_{n,n'}\ffrac{\dd}{\dd\prUp{t}_{n,n'}}
\end{equation}
in terms of the bases introduced in~\bref{app:(1,p)-proj-mod-base}.
It now follows from~\eqref{cas-plus-P} and~\eqref{eq:EmFm-prod2} that
\begin{equation*}
  \fp^m \ep^m\,\prUp{t}_{a,a'}
  =\prod_{s=0}^{m-1}
  \Bigl(
  \prDown{t}_{a,a'}\ffrac{\dd}{\dd\prUp{t}_{a,a'}}
  {}+{} \alpha^{p_-}(-1)^{r'-1}[s-a+r]_+\,[a-s]_+\Bigr)\,
  \prUp{t}_{a,a'},
  \quad m<p_+.
\end{equation*}
Writing $\fp^m \ep^m\,\prUp{t}_{a,a'}=
A_{\scriptstyle\uparrow}\prUp{t}_{a,a'}
+A_{\scriptstyle\downarrow}\prDown{t}_{a,a'}$, we see that
\begin{equation*}
  A_{\scriptstyle\downarrow}=
  \alpha^{p_- (m-1)}(-1)^{(m-1)(r'-1)}\,
  [x^1]\bigl(\CCplus{r,a}^{m}(x)\bigr),
\end{equation*}
where we define the polynomials
\begin{equation}\label{CCpolynom}
  \CCpm{r,a}^m(x)=
  \prod_{s=0}^{m - 1} (x + [s - a + r]_\pm [a - s]_\pm)
\end{equation}
and use $[x^n](f(x))$ to denote the coefficient at $x^n$.  

Similar formulas can be easily written for the action of~$\cas_-$ and
for the action of~$\cas_\pm$ on other indecomposable modules.

We note that
\begin{equation*}
  [x^0]\bigl(\CCpm{r,a}^m(x)\bigr)
  =([m]_\pm!)^2{\qbin{a}{m}}_\pm{\qbin{r\!-\!a\!+\!m\!-\!1}{m}}_\pm
\end{equation*}
and
\begin{align}\label{C-low}
  [x^1]\bigl(\CCplus{1,a}^{m}(x)\bigr) &=(-1)^{m + 1}[m]_+!
  [m\!-\!1]_+!{\qbin{a\!-\!1}{a\!-\!m}}_+^2,
  \\
  \intertext{with $[x^1]\bigl(\CCplus{1,0}^{m}(x)\bigr) =(-1)^{m +
      1}[m]_+! [m\!-\!1]_+!$ for $m\geq 1$, and}
  \label{C-top}
  [x^1]\bigl(\CCplus{p_+ - 1,a}^{m}(x)\bigr)
  &=
  (-1)^{(p_- + 1)(m + 1)} [m\!-\!1]_+! [m]_+!
  {\qbin{a}{m\!-\!1}}_+^2.
\end{align}

\subsubsection{Proof of
  \bref{prop:radford-decompose}}\label{app:radford-proof} We begin
with calculating the comultiplication of the cointegral~as
\begin{multline}\label{Delta(c)}
  \Delta(\coint)= \zeta \sum_{r=0}^{p_+ - 1}\!\sum_{m=0}^{p_+ -
    1}\!\sum_{s=0}^{p_- - 1}\!  \sum_{n=0}^{p_- - 1}
  \sum_{\ell=0}^{2 p_+ p_- - 1} \!\!
  \qp^{-p_-(m + r + 1) (m + r + 2)}
  \qm^{-p_+ (n + s + 1) (n + s + 2)}\\*
  {}\times (-1)^{r + m + n + s}
  \fp^{p_+ - r - 1} \ep^{m} \fm^{n}
  \emi^{p_- - 1 - s} k^{2\ell - 2p_-(m + 1) + 2p_+(n + 1)}\tensor{}\\*
  {}\tensor \fp^{r} \ep^{p_+ - 1 - m} \fm^{p_- - 1 - n} \emi^{s}
  k^{2\ell + 2p_-(r + 1) - 2p_+(s + 1)}
\end{multline}
and then evaluate the traces in the definition of
$\Radford^{\pm}(r,r')$
using~\eqref{eq:EmFm-prod2}--\eqref{cas-value2}.\ \ This gives
\begin{equation}\label{phi}
  \Radford^{\,\alpha}(r,r')=
  \zeta
  \Radford^{\,\alpha}_{(+)}(r,r')
  \bigl(\one + \alpha^{p_+ p_-} (-1)^{p_+(r'-1)+p_-(r-1)}K^{p_+ p_-}
  \bigr)
  \Radford^{\,\alpha}_{(-)}(r,r'),
\end{equation}  
where $1\leq r\leq p_+$ and $ 1\leq r'\leq p_-$, and
\begin{multline*}
  \Radford^{\,\alpha}_{(+)}(r,r')=
  \sum_{m=0}^{p_+ - 1}\sum_{a=m}^{r - 1}\sum_{b=0}^{p_+ - 1}
  (-1)^{b(r'-1)} \alpha^{b p_-}([m]_+!)^2  
  \\
  {}\times
  \qp^{p_- (r - 1 - 2a)(b - m)}
  {\qbin{r\!-\!a\!+\!m\!-\!1}{m}}_+
  {\qbin{a}{m}}_+
  \fp^{p_+ - m - 1} \ep^{p_+ - m - 1}K_+^{b - m}
\end{multline*}
and
\begin{multline*}
  \Radford^{\,\alpha}_{(-)}(r,r')=
  \sum_{m'=0}^{p_- - 1}\sum_{a'=0}^{r' - m' - 1}\sum_{b'=0}^{p_- - 1}
  (-1)^{b'(r-1)}\alpha^{b' p_+}([m']_-!)^2\\
  {}\times
  \qm^{p_+(r' - 1 - 2a' - 2m')(b' + m')}
  {\qbin{r'\!-\!a'\!-\!1}{m'}}_-
  {\qbin{a'\!+\!m'}{m'}}_-
  \fm^{p_- - m' - 1} \emi^{p_- - m' - 1}K_-^{b' + m'}.
\end{multline*}

As regards the central elements $\RadfordNESW(r,r')$ and
$\RadfordNWSE(r,r')$, we recall the nondiagonal action of $\cas_\pm$
on indecomposable modules, see~\bref{sec:casimirs}.  Straightforward
calculations then give the following formulas:
\begin{align}\label{radfordNESW-eval}
  \NRadfordNESW(r,r') &=
  \Nllambda^+(r, r')
  -(-1)^{p_-}
  \Nllambda^{-}(p_+\!-\!r, r')\\
  \intertext{for $1\leq r\leq p_+-1$ and $1\leq r'\leq p_-$,}
  \label{radfordNWSE-eval}
  \NRadfordNWSE(r,r')&
  =\Nbarllambda^{+}(r, r')
  -(-1)^{p_+}
  \Nbarllambda^{-}(r, p_-\!-\!r')
\end{align}
for $1\leq r\leq p_+$ and $1\leq r'\leq p_--1$, and
\begin{multline}
  \label{radfordUpUp-eval}
  \NRadfordUpUp(r,r')
  =\NLLambda^{+}(r, r')
  -(-1)^{p_-}\NLLambda^{-}(p_+\!-\!r, r')\\
  {}-(-1)^{p_+}\NLLambda^{-}(r, p_-\!-\!r')
  +(-1)^{p_- + p_+}\NLLambda^{+}(p_+\!-\!r, p_-\!-\!r')
\end{multline}
for $(r,r')\in\setR$,  where
\begin{multline*}
  \Nllambda^{\alpha}(r, r')
  = \zeta\ffrac{[r]_+}{\qp^{p_-}\!- \qp^{-p_-}}
  \sum_{m=0}^{p_+ - 1}\sum_{m'=0}^{p_- - 1}
  \sum_{a=0}^{r - 1}\sum_{a'=m'}^{r' - 1} \sum_{\ell=0}^{2 p_+ p_- -
    1}\!\!\!  \alpha^{\ell - p_-} (-1)^{r' - 1}\\
  \qquad{}\times\qp^{(\ell - p_- m)(r - 1 - 2a)}
  \qm^{(\ell + p_+ m')(r' - 1 - 2a')}
  [x^1]\bigl(\CCplus{r,a}^{m}(x)\bigr)\,
  [y^0]\bigl(\CCminus{r',a'}^{m'}(y)\bigr)\\
  {}\times\fp^{p_+ - 1 - m} \ep^{p_+ - 1 - m} \fm^{p_- - 1 - m'}
  \emi^{p_- - 1 - m'} K^{\ell - p_- m + p_+ m'},
\end{multline*}
\vspace{-1.1\baselineskip}
\begin{multline*}
  \Nbarllambda^{\alpha}(r, r')
  = \zeta\ffrac{[r']_-}{\qm^{p_+}\!- \qm^{-p_+}}\sum_ {m'=0}^{p_- - 1}\sum_{m=0}^{p_+ - 1}
  \sum_{\ell=0}^{2 p_+ p_- - 1} \sum_{a'=0}^{r' - 1}
  \sum_{a=m}^{r - 1} \alpha^{\ell - p_+} (-1)^{r - 1}\\*
  \qquad{}\times\qm^{(\ell - p_+ m')(r' - 1 - 2a')}
  \qp^{(\ell + p_- m)(r - 1 - 2a)}
  [x^0]\bigl(\CCplus{r,a}^{m}(x)\bigr)\,
  [y^1]\bigl(\CCminus{r',a'}^{m'}(y)\bigr)\\*
  {}\times\fm^{p_- - 1 - m'} \emi^{p_- - 1 - m'} \fp^{p_+ - 1 - m}
  \ep^{p_+ - 1 - m} K^{\ell - p_+ m' + p_- m},
\end{multline*}
and
\begin{multline*}  
  \NLLambda^{\alpha}(r, r')
  = \zeta\ffrac{[r]_+\,[r']_-}{(\qp^{p_-}\!-\!\qp^{-p_-})
    (\qm^{p_+}\!-\!\qm^{-p_+})}
  \sum_{m=0}^{p_+ - 1}\sum_{m'=0}^{p_- - 1}
  \sum_{a=0}^{r - 1}\sum_{a'=0}^{r' - 1} \sum_{\ell=0}^{2 p_+ p_- -
    1}\!\!\!  \alpha^{\ell + p_- + p_+} (-1)^{r' + r}\\*
  {}\times\qp^{(\ell - p_- m)(r - 1 - 2a)}
  \qm^{(\ell - p_+ m')(r' - 1 - 2a')}
  [x^1]\bigl(\CCplus{r,a}^{m}(x)\bigr)\,
  [y^1]\bigl(\CCplus{r',a'}^{m'}(y)\bigr)\\*
  {}\times\fp^{p_+ - 1 - m} \ep^{p_+ - 1 - m} \fm^{p_- - 1 - m'}
  \emi^{p_- - 1 - m'} K^{\ell - p_- m - p_+ m'}.
\end{multline*}

The rest of the proof (establishing the decomposition in terms of the
canonical central elements) is by direct evaluation of the action of
$\Radford^{\pm}(r,r')$, $\RadfordNESW(r,r')$, $\RadfordNWSE(r,r')$,
and $\RadfordUpUp(r,r')$ on the projective modules, as described
in~\bref{sec:coeffs}.  The nondiagonal part of the action on the
projective modules described in~\bref{app:proj-mod} is related to the
coefficients at the radical elements entering the decomposition with
respect to the canonical central elements.  The calculation yielding
\eqref{radfordNESW-decompose}, \eqref{radford-bdry-decomp1},
\eqref{radfordNWSE-decompose}, \eqref{radford-bdry-decomp2},
and~\eqref{phiUpUp-decompose} is very similar to the one given
in~\cite{[FGST]}, and we therefore omit it.

\end{document}